\documentclass[12pt,a4paper]{article}

\usepackage{amsmath}    
\usepackage{amssymb}
\usepackage{graphicx}   
\usepackage{verbatim}   
\usepackage{color}      

\usepackage{hyperref}   

\usepackage{enumerate}
\usepackage{mathrsfs}
\usepackage{empheq}
\usepackage{bbold}

\usepackage[final]{showlabels}

\usepackage{amsthm}

\newtheorem{lemma}{Lemma}[section]

\newtheorem{theorem}{Theorem}[section]

\newcommand{\del}{\partial}
\renewcommand{\theta}{\vartheta}
\renewcommand{\phi}{\varphi}
\newcommand{\vecc}[2]{\left ( \begin{array}{c}#1\\#2\\ \end{array}\right )}

\newcommand{\dd}{\mathrm{d}}

\newcommand{\id}{\mathbb{1}}

\renewcommand{\div}{\mathrm{div\,}}
\renewcommand{\vec}{\mathbf}

\newcommand{\sign}{\mathrm{sign}}

\newcommand{\const}{\mathrm{const}}

\newcommand{\ii}{\mathbb{i}}

\renewcommand{\title}{Truly multi-dimensional all-speed schemes for the \\Euler equations on Cartesian grids}

\newcommand{\authorOne}{Wasilij Barsukow\footnote{Institute for Mathematics, Zurich University, 8057 Zurich, Switzerland}}
\newcommand{\authorTwo}{}

\begin{document}

\begin{center} \Large
\title

\vspace{1cm}

\date{}
\normalsize

\authorOne, \authorTwo
\end{center}

\begin{abstract}

Finite volume schemes often have difficulties to resolve the low Mach number (incompressible) limit of the Euler equations. Incompressibility is only non-trivial in multiple spatial dimensions. Low Mach fixes, however generally are applied to the one-dimensional method and the method is then used in a dimensionally split way. This often reduces its stability. Here, it is suggested to keep the one-dimensional method as it is, and only to extend it to multiple dimensions in a particular, all-speed way. This strategy is found to lead to much more stable numerical methods. Apart from the conceptually pleasing property of modifying the scheme only when it becomes necessary, the multi-dimensional all-speed extension also does not include any free parameters or arbitrary functions, which generally are difficult to choose, or might be problem dependent. The strategy is exemplified on a Lagrange Projection method and on a relaxation solver.



\end{abstract}

\section{Introduction}

This paper focuses on the low Mach number limit of the Euler equations. It has been shown (e.g. \cite{ebin77,klainerman81,metivier01}) that, with appropriate initial data, in this limit the solutions fulfill the incompressible Euler equations. Many standard numerical methods for the compressible Euler equations have difficulties resolving low Mach number flow. Because incompressibility greatly simplifies in one spatial dimension (the velocity becomes spatially constant) the low Mach number limit can be considered a genuinely multi-dimensional feature. The numerical difficulties can thus be understood as being embedded into the larger context of multi-dimensionality, which amounts to the question whether multi-dimensional problems require (or, at least, profit from) an approach beyond dimensional splitting (i.e. the solution of one-dimensional problems in different directions). Before coming back to low Mach number flow, a discussion of other multi-dimensional problems therefore is due, which is followed by an overview of reasons why multi-dimensional methods have appeared in the literature so far.

The development of numerical methods for hyperbolic conservation laws for a long time focused on the ability of the numerical method to remain stable (e.g. \cite{courant28,vanleer06}), to converge to a weak solution of the problem (e.g. \cite{lax64}), to respect a discrete entropy inequality (e.g. \cite{tadmor2003entropy}), to resolve discontinuities sharply and to attain a high order of accuracy (e.g. \cite{vanleer79,colella84,jiang96,cockburn98}). All these aspects have mostly been studied in a one-dimensional setting, and the extension to multiple spatial dimensions using dimensional splitting yielded seemingly satisfactory results.

In multiple spatial dimensions, however, hyperbolic conservation laws have a much richer phenomenology. Vortices have complicated dynamics, interactions of shocks can be intricate in multiple dimensions, contact discontinuities can become unstable (slip lines) and thus change the solution dramatically with respect to one-dimensional counterparts. At the same time, multi-dimensional simulations require much more resources (in particular in terms of computational time) than one-dimensional ones. Thus, in multiple dimensions one faces the challenge of resolving many more features with much less resources. 

Truly multi-dimensional schemes so far arose in two domains. The strict application of Godunov's procedure of \emph{reconstruction}--\emph{evolution}--\emph{averaging} in multi-d requires the solution of the multi-dimensional Riemann problems arising at the vertices of the grid. This is possible for linear problems \cite{abgrall93,leveque97,gilquin96,barsukow17}, but is very difficult to accomplish in an exact manner for the Euler equations. A lot of effort (e.g. \cite{colella90,zheng12}) therefore has been devoted to understanding which essential properties of multi-dimensional interactions need to be taken into account and how an approximate solution of the multi-d Riemann problem should look like. In \cite{barsukow17} a limitation of this approach has been uncovered by showing that even the exact solution of a multi-dimensional problem does not lead to a low Mach number compliant scheme.

\begin{figure}
 \centering
 \includegraphics[width=0.45\textwidth]{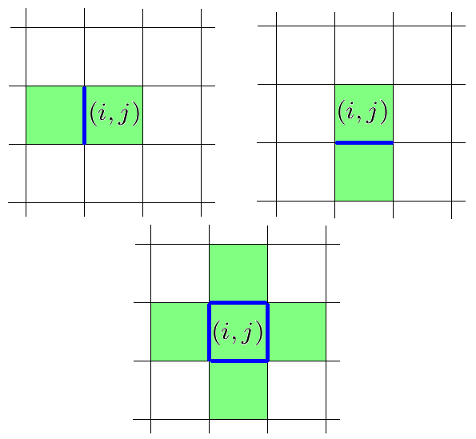} \hfill \includegraphics[width=0.45\textwidth]{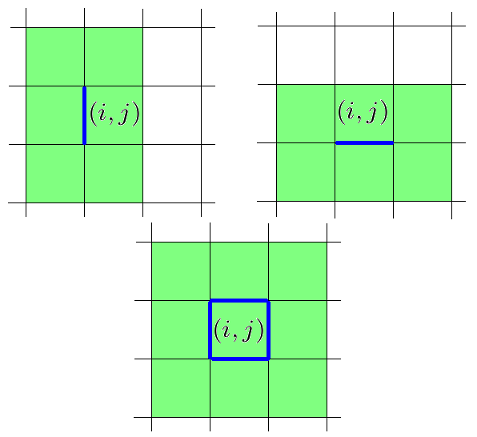}
 \caption{\emph{Left}: If the flux depends only on the values in the two cells adjacent to a cell interface, then the update formula uses a five-point stencil on Cartesian grids. \emph{Right}: A truly multi-dimensional scheme also uses values in cells which only share a corner.}
 \label{fig:stencils}
\end{figure}

Another case in which truly multi-dimensional schemes have appeared naturally is the quest for preserving discrete structures, such as involutions. On Cartesian grids, this has been extensively studied in the context of linear acoustics, which possesses an involution (vorticity) (e.g. \cite{morton01,jeltsch06,mishra09preprint,barsukow17a}). It turns out that including the corner-cells of a $3 \times 3$ Cartesian stencil (see Figure \ref{fig:stencils}) allows to easily construct numerical methods that preserve a discrete counterpart of the involution. Discrete preservation of an involution in particular means that its time evolution is exempt from any numerical diffusion.

{The difficulties of many numerical methods to resolve the low Mach number limit of the Euler equations have also been attributed to an excessive numerical diffusion. It has been shown in e.g. \cite{guillard04} that a (dimensionally split) Riemann solver is at the origin of the difficulties. Interestingly, it has been found in \cite{rieper09,dellacherierieper10,guillard09} that on triangular and tetrahedral computational grids the Roe solver is low Mach number compliant, while the HLL solver is not (\cite{rieper10}).} {Many numerical methods have been suggested which therefore, partly, replace the Riemann solvers by central differences (e.g. \cite{haack12,cordier12}), at the expense of having to use at least partially implicit (e.g. IMEX) time integrators. A time-implicit method, however, is costly, and whenever the solution involves both high and low Mach number flows the time steps need to be chosen according to the acoustic time scale anyway. In this work, low Mach number compliance shall be achieved while retaining a diffusion sufficient to allow explicit time integration.}

{So far, modifications of the diffusion that result in low Mach compliance were performed on the one-dimensional scheme.} In particular, {using asymptotic analysis}, terms (typically, discretizations of derivatives of velocity components $\del_x u$) were identified {which were not compliant with the low Mach number limit.} These terms were then multiplied by some factor $\mathcal O(\epsilon)$ (e.g. \cite{li08,thornber08,dellacherie10,rieper11,li13,chalons16,birken16,barsukow16,dellacherie16}). {They thus appeared at one order higher in the asymptotic analysis, and did not present a problem any longer.} Then this scheme was applied to multi-d in a dimensionally split manner. Such an approach involves the potential risk of spoiling stability properties of the scheme (see e.g. \cite{birken05,barsukow16} for examples of reduced CFL-conditions). The restriction to dimensionally split schemes comes at the expense of accepting free parameters and reduced stability.

In \cite{barsukow17a}, the low Mach number problem was investigated {on Cartesian grids} for the equations of linear acoustics. It has been found that low Mach number compliance leads to the same condition as the preservation of a discrete vorticity. 
From the viewpoint of low Mach compliance this has opened up a new avenue, as numerical methods initially devoted to vorticity preservation came into focus of low Mach number analysis. Indeed, one can draw a close parallel between the preservation of a discrete structure (such as vorticity) and low Mach number compliance. The low Mach number compliant numerical method needs to ensure that the numerical diffusion does not interfere with this limit. For standard numerical methods, the numerical diffusion becomes dominant in the low Mach limit and wipes out all the features of the flow. It is this aspect that this paper focuses on.

It is largely accepted that a certain amount of diffusion is unavoidable to achieve a stable numerical method. For systems of PDEs, and even more so when they are to be solved in multiple spatial dimensions, it is, however, unclear how little diffusion actually is necessary. Structure preserving methods (e.g. asymptotic preserving, involution preserving, well-balanced) can be understood as reducing or removing the numerical diffusion on certain setups or on the evolution of certain quantities. This requires a careful control over the numerical diffusion and this is where truly multi-dimensional numerical methods turn out to be useful. On Cartesian grids, dimensionally split methods, despite being slightly easier to implement, are perceived as a huge restriction when the numerical diffusion is to be tailored in a specific way. It is, for example, clear that allowing for a 9-point stencil in 2D leads to a much larger choice of finite difference discretizations than on a 5-point stencil (see also Figure \ref{fig:stencils}).
Another way to express this, is to say that while focusing on each direction individually, dimensionally split methods neglect the way contributions from different directions balance each other. 

Structure preserving (e.g. low Mach compliant) methods greatly improve the quality of the simulation. In order to reach the same results, a standard method would need to be run on a highly resolved numerical grid, which often is impossible in practice. Truly multi-dimensional methods which are structure preserving are much easier to construct than dimensionally split ones, and are endowed with a lot of favourable properties, while dimensionally split methods need to compromise on them.
It is important to emphasize that the reduction of the numerical diffusion and the improved quality of simulation do not entail any change in computational time or memory usage. This is particularly important in multiple spatial dimensions, where the computational resources are limited.

The new strategy, pursued in this paper, is to use a good one-dimensional scheme (endowed e.g. with known stability properties) and to \emph{extend} it to the multi-dimensional case on a Cartesian grid in a particular, all-speed manner. Instead of removing derivatives of velocity components, they are completed to make the divergence appear, knowing that in the low Mach limit the divergence is $\mathcal O(\epsilon)$. Thus, a term such as a discretization of a second derivative $\del_x^2 u$ of an individual velocity component $u$ would become $\del_x^2 u + \del_x \del_y v$ in the multi-dimensional case. The new scheme remains the same in one spatial dimension, where there is no problem with low Mach compliance (because the incompressible limit is trivial). In multi-d, however, it is all-speed without the need for free parameters. One can expect that the multi-dimensional scheme inherits many of the favorable properties of the one-dimensional scheme, in particular stability. The appearance of cross-derivatives is a sign that the resulting scheme needs to be truly multi-dimensional, and cannot be achieved with a dimensionally split approach.

In this paper the new approach of an all-speed extension to multiple dimensions is applied to two one-dimensional schemes which are known to have provable stability properties. They are chosen to represent two philosophies: a Lagrange Projection method and a relaxation solver. Both schemes are analyzed and extended to multiple dimensions in an all-speed fashion.

An all-speed version of the Lagrange Projection scheme has already been suggested in \cite{chalons16} using a low Mach fix applied to the one-dimensional method. In \cite{chalons16} it was also noticed that the same low Mach fix applied to a relaxation solver yields an unstable method. This paper suggests an alternative strategy, which provides stable all-speed methods in both cases, and on top of that is free of a number of inconvenient properties of standard low Mach fixes.

The paper is organized as follows: section \ref{sec:overview} presents a more detailed overview of existing low Mach number modifications, and in particular the multi-dimensional approach for linear acoustics. Taking it as inspiration, in section \ref{sec:lagrange} a one-dimensional Lagrange-Projection scheme is extended to multi-d; similarly in section \ref{sec:relaxation} a one-dimensional relaxation solver is extended to multi-d. Low Mach compliance of the resulting schemes is shown in both cases. Numerical examples are presented in section \ref{sec:numerical}.

\section{Overview of low Mach number modifications} \label{sec:overview}
\subsection{Linear acoustics}\label{sec:overviewacoustics}

The equations of linear acoustics arise as a linearization of the Euler equations around the state of constant density and pressure and vanishing velocity:
\begin{align}
 \del_t \vec v + \frac{\nabla p}{\epsilon^2} &= 0 \label{eq:acoustic2dv}\\
 \del_t p + c^2 \nabla \cdot \vec v &= 0  \label{eq:acoustic2dp}
\end{align}
where
\begin{align}
 p &:\mathbb R^+_0 \times \mathbb R^d \to \mathbb R \\
 \vec v &:\mathbb R^+_0 \times \mathbb R^d \to \mathbb R^d  & c &\in \mathbb R^+
\end{align}
Here and in the following, $d$ denotes the number of spatial dimensions, for simplicity $d=2$ is mostly assumed. Also, the notation $\vec v = (u,v)$ is used below. 

For a detailed derivation of the acoustic equations and their properties see \cite{barsukow17,barsukow17a}. In analogy to the rescaling of the Euler equations in the low Mach number limit $\epsilon \to 0$ (see e.g. \cite{klein95,barsukow16}) the acoustic equations \eqref{eq:acoustic2dv}--\eqref{eq:acoustic2dp} are endowed with a $\frac{1}{\epsilon^2}$-scaling of the pressure gradient. {(See section \ref{ssec:euler} for more details on the rescaling procedure for the Euler equations).} As $\epsilon \to 0$, the equations \eqref{eq:acoustic2dv}--\eqref{eq:acoustic2dp} formally reduce to
\begin{align}
 \nabla p &\in \mathcal O(\epsilon^2) & \nabla \cdot \vec v &\in \mathcal O(\epsilon) \label{eq:exactlimitequationsacoustics}
\end{align}

Apart from the low Mach number limit, the acoustic equations are of interest due to the presence of an involution. Indeed, applying the curl to \eqref{eq:acoustic2dv}, one obtains
\begin{align}
 \del_t (\nabla \times \vec v )& =0
\end{align}
By analogy with the Euler equations, $\nabla \times \vec v $ is referred to as \emph{vorticity}. Numerical methods that would keep a discretization of the vorticity stationary (\emph{vorticity preserving}) have been subject of various investigations, e.g. \cite{morton01,jeltsch06,mishra09preprint,lung14}. All these methods have in common that they are truly multi-dimensional.

In \cite{barsukow17a} it was shown that vorticity preservation is equivalent to low Mach number compliance for linear acoustics. Moreover, these properties are also equivalent to stationarity preservation: a method that is low Mach number compliant (or vorticity preserving) also possesses a discretization of the divergence $\nabla \cdot \vec v$, which together with a constant pressure characterizes its stationary states. (This is not true for any numerical method generally). The reader is referred to \cite{barsukow17a,barsukow18hypproceeding} for all details of these relationships. Although a dimensionally split method can have these properties, they often come at the expense of decreased stability and the introduction of ad-hoc parameters. Truly multi-dimensional methods achieve the same with greater ease. In the following a simple reason for the natural appearance of multi-dimensional schemes in the context of low Mach compliance is given. 

The intermediate state of the exact one-dimensional Riemann solver for linear acoustics is
\begin{align}
 p^*_{i+\frac12} &= \frac{p_i + p_{i+1}}{2} - c \epsilon \frac{u_{i+1} - u_i}{2} \label{eq:intermediate1dacousticp}\\
 u^*_{i+\frac12} &= \frac{u_i + u_{i+1}}{2} - \frac{1}{c \epsilon} \frac{p_{i+1} - p_i}{2} \label{eq:intermediate1dacousticu}
\end{align}
and the velocity equation on an equidistant one-dimensional grid becomes
\begin{align}
 \del_t u_i + \frac{p_{i+1} - p_{i-1}}{2 \Delta x \epsilon^2} - \frac{c}{\epsilon} \frac{u_{i+1} - 2 u_i + u_{i-1}}{2 \Delta x} &= 0
\end{align}

The dimensionally split scheme on a two-dimensional Cartesian grid then reads
\begin{align}
 \del_t u_{ij} + \frac{p_{i+1,j} - p_{i-1,j}}{2 \Delta x \epsilon^2} - \frac{c}{\epsilon} \frac{u_{i+1,j} - 2 u_{ij} + u_{i-1,j}}{2 \Delta x} &= 0\\
 \del_t v_{ij} + \frac{p_{i,j+1} - p_{i,j-1}}{2 \Delta y \epsilon^2} - \frac{c}{\epsilon} \frac{v_{i,j+1} - 2 v_{ij} + v_{i,j-1}}{2 \Delta y} &= 0\\
 \del_t p_{ij} + c^2 \left( \frac{u_{i+1,j} - u_{i-1,j}}{2 \Delta x } + \frac{v_{i,j+1} - v_{i,j-1}}{2 \Delta y } \right )\phantom{mmmm}& \nonumber \\- \frac{c}{\epsilon} \frac{p_{i+1,j} - 2 p_{ij} + p_{i-1,j}}{2 \Delta x} - \frac{c}{\epsilon} \frac{p_{i,j+1} - 2 p_{ij} + p_{i,j-1}}{2 \Delta y} &= 0
\end{align}

Expanding formally every quantity as power series in $\epsilon$, e.g.
\begin{align}
 p_{ij} = p^{(0)}_{ij} +p^{(1)}_{ij} \epsilon + p^{(2)}_{ij} \epsilon^2 + \mathcal O(\epsilon^3)
\end{align}
and collecting order by order yields, to lowest orders
\begin{align}
 \frac{p^{(0)}_{i+1,j} - p^{(0)}_{i-1,j}}{2 \Delta x} = \frac{p^{(0)}_{i,j+1} - p^{(0)}_{i,j-1}}{2 \Delta y} = 0\\
 \frac{p^{(1)}_{i+1,j} - p^{(1)}_{i-1,j}}{2 \Delta x} - c \frac{u^{(0)}_{i+1,j} - 2 u^{(0)}_{ij} + u^{(0)}_{i-1,j}}{2 \Delta x} &= 0 \label{eq:roeacousticp1x}\\
 \frac{p^{(1)}_{i,j+1} - p^{(1)}_{i,j-1}}{2 \Delta y} - c \frac{v^{(0)}_{i,j+1} - 2 v^{(0)}_{ij} + v^{(0)}_{i,j-1}}{2 \Delta y} &= 0 \label{eq:roeacousticp1y}
\end{align}

Comparing to \eqref{eq:exactlimitequationsacoustics}, one observes the $\mathcal O(1/\epsilon)$ equations \eqref{eq:roeacousticp1x}--\eqref{eq:roeacousticp1y} to include discrete second derivatives of the velocity. Thus, $p^{(1)}= \const$ is, in general, not a solution of \eqref{eq:roeacousticp1x}--\eqref{eq:roeacousticp1y}. This is identified in e.g. \cite{guillard99} as the source of the low Mach number artefacts of the scheme. {Note that \eqref{eq:roeacousticp1x}--\eqref{eq:roeacousticp1y} do not, in general, imply a discrete version of $\nabla p^{(1)} \neq 0$, and thus $\nabla p \in \mathcal O(\epsilon)$. They can also imply a discrete version of $\nabla p^{(1)} = 0$, and additionally a discrete version of $\del_x^2 u^{(0)} = 0 = \del_y^2 v^{(0)}$. Whereas it is a priori unclear which case will occur, both these outcomes are unacceptable.} In e.g. \cite{dellacherie10,barsukow16} it is suggested to multiply the second derivatives of the velocity by a factor $\mathcal O(\epsilon)$ such that they cease to appear in the $\mathcal O(1/\epsilon)$ equations.

The truly multi-dimensional numerical scheme (e.g. the one from \cite{barsukow17a}), on the other hand, can be associated with the following intermediate state of the pressure:
\begin{align}
 p^*_{i+\frac12,j} &= \frac{p_{i,j-1} + p_{i+1,j-1} + 2(p_{ij} + p_{i+1,j}) + p_{i,j-1} + p_{i+1,j-1}}{8} \nonumber\\
 &- c \epsilon \left( \frac{u_{i+1,j+1} - u_{i,j+1} + 2(u_{i+1,j} - u_{i,j}) + u_{i+1,j-1} - u_{i,j-1}}{8} \right . \\ &+ \left. \frac{\Delta x}{\Delta y} \frac{v_{i,j+1} - v_{i,j-1} + v_{i+1,j+1} - v_{i+1,j-1}}{8}    \right )
\end{align}
(One may note that this state does \emph{not} arise from an exact multi-dimensional Riemann problem.)

In order to cope with the lengthy expressions, the following notation is used in the following:
\begin{align}
 [q]_{i+\frac12} &:= q_{i+1} - q_i & \{q\}_{i+\frac12} &:= q_{i+1} + q_i\\
 [q]_{i\pm1} &:= q_{i+1} - q_{i-1}\\
 [[q]]_{i\pm\frac12} &:= [q]_{i+\frac12} - [q]_{i-\frac12} & \{ \{ q \}\}_{i\pm\frac12}&:= \{q\}_{i+\frac12} + \{q\}_{i-\frac12}\\
 &= q_{i+1} - 2 q_i + q_{i-1} & &= q_{i+1} + 2 q_i + q_{i-1}
\end{align}
The only nontrivial identity is
\begin{align}
 \{[q]\}_{i\pm\frac12} = [q]_{i+\frac12} + [q]_{i-\frac12} = [q]_{i\pm1}
\end{align}
For multiple dimensions the notation is combined, e.g.
\begin{align}
 [[q_i]]_{j\pm\frac12} &= q_{i,j+1} - 2 q_{ij} + q_{i,j-1}\\
 \{[q]_{i+\frac12}\}_{j+\frac12} &= q_{i+1,j+1} - q_{i,j+1} + q_{i+1,j} - q_{ij}\\
 [[q]_{i\pm1}]_{j\pm1} &= q_{i+1,j+1} - q_{i-1,j+1} - q_{i+1, j-1} + q_{i-1,j-1}
\end{align}

The intermediate state of the truly multi-dimensional scheme then can be rewritten as
\begin{align}
 p^*_{i+\frac12,j} &= \frac{\{\{  \{ p\}_{i+\frac12} \}\}_{j\pm\frac12}}{8} - c \epsilon \left( \frac{\{\{ [u]_{i+\frac12} \}\}_{j\pm\frac12}}{8} + \frac{\Delta x}{\Delta y} \frac{ [\{ v\}_{i+\frac12} ]_{j\pm1}}{8}    \right ) \label{eq:intermediatemultidacousticpx}\\
 p^*_{i,j+\frac12} &= \frac{\{\{  \{ p\}\}_{i\pm\frac12} \}_{j+\frac12}}{8} - c \epsilon \left( \frac{\Delta y}{\Delta x} \frac{ \{ [u]_{i\pm1} \}_{j+\frac12}}{8} +  \frac{ [\{ \{ v\}\}_{i\pm\frac12} ]_{j+\frac12}}{8}    \right )\\
 u^*_{i+\frac12,j} &= \frac{\{\{  \{ u\}_{i+\frac12} \}\}_{j\pm\frac12}}{8} - \frac{1}{c\epsilon} \frac{\{\{ [p]_{i+\frac12} \}\}_{j\pm\frac12}}{8} \label{eq:intermediatemultidacousticu}
 \end{align}

Observe that instead of the derivatives of individual velocity components now $\frac12 \{ \mathscr D_{i+\frac12} \}_{j\pm\frac12}$ and $\frac12 \{ \mathscr D \}_{i\pm\frac12,j+\frac12}$ appear in the intermediate pressure, where $\mathscr D$ is the discrete divergence
\begin{align}
 \mathscr D_{i+\frac12,j+\frac12} := \frac{ \{ [u]_{i+\frac12}\}_{j+\frac12}}{2 \Delta x} + \frac{  [\{v\}_{i+\frac12}]_{j+\frac12}}{2 \Delta y}
\end{align}
Observe also that the intermediate state \eqref{eq:intermediatemultidacousticpx}--\eqref{eq:intermediatemultidacousticu} reverts to the usual one (equations \eqref{eq:intermediate1dacousticp}--\eqref{eq:intermediate1dacousticu}) in a one-dimensional situation.

The numerical scheme then reads 
\begin{align}
 \del_t u_{ij} + \frac{\{\{ [p]_{i\pm1} \}\}_{j\pm\frac12}}{8 \Delta x \epsilon^2} - \frac{c}{\epsilon} \left( \frac{\{\{ [[u]]_{i\pm\frac12} \}\}_{j\pm\frac12} }{8 \Delta x} + \frac{ [[v]_{i\pm1}]_{j\pm1} }{8 \Delta y} \right )&= 0 \label{eq:multidschemeacousticu}\\
 \del_t v_{ij} + \frac{ [\{\{ p \}\}_{i\pm\frac12}]_{j\pm1}}{8 \Delta y \epsilon^2} - \frac{c}{\epsilon} \left( \frac{ [[u]_{i\pm1}]_{j\pm1} }{8 \Delta x}  +   \frac{[[   \{\{ v \}\}_{i\pm\frac12} ]]_{j\pm\frac12}}{8 \Delta y} \right )&= 0 \label{eq:multidschemeacousticv}\\
 \del_t p_{ij} + c^2 \left( \frac{\{\{ [u]_{i\pm1} \}\}_{j\pm\frac12}}{8 \Delta x}  + \frac{ [\{\{ v \}\}_{i\pm\frac12}]_{j\pm1}}{8 \Delta y} \right )\phantom{mmmm}& \nonumber \\- \frac{c}{\epsilon} \left( \frac{\{\{ [[p]]_{i\pm\frac12} \}\}_{j\pm\frac12} }{8 \Delta x} + \frac{[[   \{\{ p \}\}_{i\pm\frac12} ]]_{j\pm\frac12}}{8 \Delta y} \right ) &= 0 \label{eq:multidschemeacousticp}
\end{align}
and again collecting order by order yields
\begin{align}
 \frac{\{\{ [p^{(0)}]_{i\pm1} \}\}_{j\pm\frac12}}{8 \Delta x} = \frac{ [\{\{ p^{(0)} \}\}_{i\pm\frac12}]_{j\pm1}}{8 \Delta y} = 0\\
 \frac{\{\{ [p^{(1)}]_{i\pm1} \}\}_{j\pm\frac12}}{8 \Delta x} - \frac{c}{2} \frac{\{[ \mathscr D^{(0)} ]_{i\pm\frac12} \}_{j\pm\frac12} }{2}&= 0 \label{eq:fstorderacousticmultidpx}\\
 \frac{ [\{\{ p^{(1)} \}\}_{i\pm\frac12}]_{j\pm1}}{8 \Delta y} - \frac{c}{2} \frac{[\{ \mathscr D^{(0)} \}_{i\pm\frac12} ]_{j\pm\frac12} }{2}&= 0 \label{eq:fstorderacousticmultidpy}\\
 c^2 \frac{ \{\{ \mathscr D^{(0)}\}_{i\pm\frac12} \}_{j\pm\frac12} }{4} - c \left( \frac{\{\{ [[p^{(1)}]]_{i\pm\frac12} \}\}_{j\pm\frac12} }{8 \Delta x} + \frac{[[   \{\{ p^{(1)} \}\}_{i\pm\frac12} ]]_{j\pm\frac12}}{8 \Delta y} \right ) &= 0 \label{eq:fstorderacousticmultidu}
\end{align}
Here, the particular shape of the divergence operators and of the additional perpendicular averaging is crucial for the appearance of the same divergence operator in the different equations. One observes that \eqref{eq:fstorderacousticmultidpx}--\eqref{eq:fstorderacousticmultidu} is solved by
\begin{align}
 \frac{\{\{ [p^{(1)}]_{i\pm1} \}\}_{j\pm\frac12}}{8 \Delta x} = \frac{ [\{\{ p^{(1)} \}\}_{i\pm\frac12}]_{j\pm1}}{8 \Delta y} = 0\\
 \mathscr D_{i+\frac12,j+\frac12}^{(0)} = \frac{ \{ [u^{(0)}]_{i+\frac12}\}_{j+\frac12}}{2 \Delta x} + \frac{  [\{v^{(0)}\}_{i+\frac12}]_{j+\frac12}}{2 \Delta y} &= 0 
\end{align}
These are discretizations of \eqref{eq:exactlimitequationsacoustics}.

While previously the velocity derivatives needed to be multiplied by an extra factor of $\epsilon$ to achieve low Mach number compliance, now they appear as derivatives of the divergence, which is $\mathcal O(\epsilon)$ and in the low Mach number limit automatically contributes a factor of $\epsilon$. The extension of this approach to 3-d is entirely analogous.

The advantage of this approach is twofold: it does not require any free parameters or ad hoc factors, and the one-dimensional numerical method is unchanged. The good performance of the scheme has been demonstrated experimentally in \cite{barsukow17a,barsukow18thesis}, and a von Neumann stability analysis for this scheme is provided here for the first time in \ref{ap:stability}. The aim of the next chapters is to apply the same procedure to the Euler equations in order to obtain parameter-free all-speed schemes. The linearized versions of these schemes will basically reduce to the acoustic scheme described above, such that its stability analysis is relevant for these new schemes as well.

\subsection{Euler equations} \label{ssec:euler}

The compressible Euler equations read

{
\begin{align}
 \del_t \rho + \nabla \cdot (\rho \vec v) &= 0 \label{eq:eulerunscaled1}\\
 \del_t (\rho \vec v) + \nabla \cdot (\rho \vec v \otimes \vec v + p \id) &= 0\\
 \del_t e + \nabla \cdot(\vec v (e+p)) &= 0 \label{eq:eulerunscaled3}
\end{align}
with
\begin{align}
 \rho,e,p &: \mathbb R^+_0 \times \mathbb R^d \to \mathbb R^+\\
 \vec v &: \mathbb R^+_0 \times \mathbb R^d \to \mathbb R^d
\end{align}
and
\begin{align}
 e &= \frac{p}{\gamma-1} + \frac12 \rho |\vec v|^2 \qquad \gamma > 1
\end{align}}

In two dimensions, again the notation $\vec v = (u, v)^\text T$ is used. The conservative variables are denoted by $q = (\rho, \rho u, \rho v, e)^\text T$. In the numerical examples, $\gamma = 1.4$ is used. 

{It is convenient to ``rescale'' the Euler equations in order to exhibit the low Mach number limit more prominently. Mostly, this is done by choosing some ``reference quantities''. However, it remains unclear what consequences a particular choice of these quantities has, and the physical arguments easily get mixed up with mathematical necessity. Here, a derivation (first appeared in \cite{barsukow16}) is provided which postpones all physical arguments to the very end, thus showing that most of the ``rescaling'' is a mathematical consequence of the equations alone.

Consider a family of solutions $(\rho_\epsilon, \vec v_\epsilon, p_\epsilon)$ of \eqref{eq:eulerunscaled1}--\eqref{eq:eulerunscaled3} parametrized by $\epsilon \in \mathbb R^+$, so far without any relation to Mach number. Assume every quantity to be a power series in $\epsilon$ and make the leading order explicit:
\begin{align}
 \rho_\epsilon(t, \vec x) &= \epsilon^\mathfrak a \tilde \rho_\epsilon(t, \vec x) & \tilde \rho_\epsilon(t, \vec x) &= \tilde \rho^{(0)}(t,\vec  x) + \epsilon \tilde \rho^{(1)}(t, \vec x) + \mathcal O(\epsilon^2) \label{eq:rescalingrho}
\end{align}
and similarly for all other quantities, each with its own leading power of $\epsilon$:
\begin{align}
 \vec v_\epsilon(t, \vec x) &= \epsilon^\mathfrak b \tilde{ \vec v}_\epsilon(t, \vec x) & 
 p_\epsilon(t,\vec  x) &= \epsilon^\mathfrak c \tilde p_\epsilon(t, \vec x)   \label{eq:rescalingpressurevelocity}
\end{align}
Consider also the independent variables to depend on $\epsilon$:
\begin{align}
 \vec x &= \epsilon^\mathfrak d \tilde {\vec x} & t &= \epsilon^\mathfrak e \tilde t & \mathfrak a, \mathfrak b, \mathfrak c, \mathfrak d, \mathfrak e \in \mathbb Z  \label{eq:rescalingxt}
\end{align}
which can be realized as an $\epsilon$-dependent system size or an $\epsilon$-dependent final time.

Computing the Mach number yields
\begin{align}
 M_\epsilon := \frac{|\vec v_\epsilon|}{\sqrt{\gamma p_\epsilon / \rho_\epsilon}} = \epsilon^{\mathfrak b - \frac{\mathfrak c}{2} + \frac{\mathfrak a}{2}} \frac{|\tilde {\vec v}_\epsilon|}{\sqrt{\gamma \tilde p_\epsilon / \tilde \rho_\epsilon} }
\end{align}
At this point, $\epsilon$ shall be chosen to control the Mach number, which implies $\mathfrak b - \frac{\mathfrak c}{2} + \frac{\mathfrak a}{2} = 1$. Now $M_\epsilon \in \mathcal O(\epsilon)$. 

Inserting \eqref{eq:rescalingrho}--\eqref{eq:rescalingxt} into the Euler equations then yields 
\begin{align}
 \epsilon^{\mathfrak d - \mathfrak b - \mathfrak e} \del_{\tilde t} \tilde \rho_\epsilon + \tilde \nabla \cdot (\tilde \rho_\epsilon \tilde {\vec v_\epsilon}) &= 0 \label{eq:eulerintermd1}\\
 \epsilon^{\mathfrak d - \mathfrak b - \mathfrak e} \del_{\tilde t} (\tilde \rho_\epsilon \tilde{\vec v}_\epsilon) + \tilde \nabla \cdot (\tilde\rho_\epsilon \tilde {\vec v}_\epsilon \otimes \tilde {\vec v}_\epsilon + \tilde p_\epsilon/\epsilon^2 \id) &= 0\\
 \epsilon^{\mathfrak d - \mathfrak b - \mathfrak e} \del_{\tilde t} e_\epsilon + \tilde \nabla \cdot(\tilde {\vec v}_\epsilon (\tilde e_\epsilon+\tilde p_\epsilon)) &= 0 \label{eq:eulerintermd3}
\end{align}
and
\begin{align}
 \tilde e_\epsilon &= \frac{\tilde p_\epsilon}{\gamma-1} + \frac12 \epsilon^2 \tilde \rho_\epsilon |\tilde {\vec v}_\epsilon|^2  \label{eq:eulerintermd4}
\end{align}

It is well-known that the Euler equations contain two dimensionless numbers, the Mach number and the Strouhal number $\frac{|\vec x|}{|\vec v_\epsilon|t} = \epsilon^{\mathfrak d - \mathfrak b - \mathfrak e} \frac{|\tilde{\vec x}|}{|\tilde {\vec v}_\epsilon|\tilde t}$. By choosing $\mathfrak d - \mathfrak b - \mathfrak e = 0$ one ensures that the time and length scales match those of the fluid velocity. This is the only decision taken on physical grounds in the rescaling procedure. Note also that it is not necessary to define a length scale, a time scale and so on, because these details are, in fact, irrelevant for the outcome of the rescaling. An example of \eqref{eq:rescalingrho}--\eqref{eq:rescalingxt} with $\mathfrak c = -2$, $\mathfrak a = \mathfrak b = \mathfrak d = \mathfrak e = 0$ is given by the vortex test in section \ref{sec:numerical}.

The rescaling only serves the purpose of keeping track of the influence of a decreasing Mach number, and thus is only needed for the theoretical analysis. In the practical application, the numerical scheme is solving the usual Euler equations \eqref{eq:eulerunscaled1}--\eqref{eq:eulerunscaled3}.

For better readability, the tilde and the $\epsilon$-subscript of all quantities is dropped from now on, such that the rescaled Euler equations read
\begin{align}
 \del_t  \rho + \nabla \cdot ( \rho  {\vec v}) &= 0 \\
 \del_t ( \rho {\vec v}) +  \nabla \cdot (\rho  {\vec v} \otimes  {\vec v} +  p/\epsilon^2 \id) &= 0\\
 \del_t e +  \nabla \cdot( {\vec v} ( e+ p)) &= 0 
\end{align}
and
\begin{align}
  e &= \frac{ p}{\gamma-1} + \frac12 \epsilon^2  \rho | {\vec v}|^2  
\end{align}
}

As $\epsilon \to 0$, the limit is characterized by
\begin{align}
 \nabla p &\in \mathcal O(\epsilon^2) &  \nabla \cdot \vec v &\in \mathcal O(\epsilon)
\end{align}
which is very similar to the acoustic case. Also the formal asymptotic expansion of standard dimensionally split schemes (e.g. the Roe scheme) is very similar to the acoustic case with velocity derivatives appearing in the $\mathcal O(1/\epsilon)$-equations. More details on this can be found in \cite{klainerman81,metivier01,klein95,guillard99,barsukow16}. Existing low Mach number modifications focus on the one-dimensional scheme which is then used in a dimensionally split way (e.g. \cite{turkel87,li08,thornber08,dellacherie10,rieper11,li13,chalons16,birken16,barsukow16,dellacherie16}), or alternatively on removing the diffusion altogether and using an implicit time integration to stabilize the central derivatives.

\subsection{Time integration}

For explicit time integration methods, the time step depends on the maximal eigenvalue of the Jacobian. In the low Mach number limit of the Euler equations, the speed of sound surpasses by far the speed of the fluid. In order to resolve any noticeable fluid flow one thus needs many time steps, which is impractical. In those cases when no high speed phenomena and effects of compressibility (shocks, or sound waves) need to be resolved, an implicit time integrator can greatly accelerate the simulation, e.g. \cite{haack12,cordier12}. To balance the advantages of explicit and implicit time integration, IMEX schemes (e.g. \cite{haack12,bispen17,boscarino19,thomann19}) have been suggested. In all cases the implicit time integration allows to use central derivatives, such that there is no numerical diffusion to interfere with the low Mach number limit. In this paper, however, the focus lies on maintaining sufficient numerical diffusion for the usage of an explicit time integrator. The analysis therefore focuses on semi-discrete methods, and the numerical examples are obtained using forward Euler.

\section{Multi-dimensional modification of a Lagrange-Pro\-jec\-tion scheme} \label{sec:lagrange}

\subsection{Review of the one-dimensional method}

The first numerical method to be discussed here is a time-explicit scheme on Cartesian grids suggested in \cite{chalons16,despres17}. It splits the Euler equations into acoustics and pressureless Euler equations (advection and compression)
\begin{align}
 \del_t \rho &&&+ (u \del_x + v \del_y) \rho &&+ \rho (\del_x u + \del_y v) &&= 0 \label{eq:euler1}\\
 \del_t (\rho u) &&&+ (u \del_x + v \del_y) (\rho u) &&+ \rho u (\del_x u + \del_y v) &&+ \frac{\del_x p}{\epsilon^2} = 0\\
 \del_t (\rho v) &&&+ (u \del_x + v \del_y) (\rho v) &&+ \rho v (\del_x u + \del_y v) &&+ \frac{\del_y p}{\epsilon^2} = 0\\
 \del_t e &&&+ (u \del_x + v \del_y) e &&+ e (\del_x u + \del_y v) &&+ \del_x (up) + \del_y (vp)= 0 \label{eq:euler4}\\
 \nonumber &&&\text{\,\,\,\,\,\,advection}&&\text{\,\,\,\,\,\,compression}&&\text{nonlinear acoustics}
\end{align}

This splitting is tightly linked to Lagrange-Projection schemes, but is interesting for others reasons as well. The different treatment of acoustics and (nonlinear) advection parallels observations (e.g. \cite{roe17}) of the very different structure of these operators in multiple spatial dimensions and the possible need to treat them differently numerically. 

Although the scheme from \cite{chalons16} is an Eulerian scheme, it carries distinctive features of its Lagrangian origin. Here, the derivation in one spatial dimension is recalled briefly and the reader is referred to \cite{chalons16} for more details.

First, the acoustic system is solved by means of a 1D relaxation solver:
\begin{align}
 u^*_{i+\frac12} &= \frac{u_{i+1} + u_{i}}{2} - \frac1{2a \epsilon} (p_{i+1} - p_{i}) = \frac{ \{ u \}_{i+\frac12} }{2} - \frac1{2a \epsilon} [p]_{i+\frac12} \label{eq:1dlagprojrelaxacu}\\
 p^*_{i+\frac12} &= \frac{p_{i+1} + p_{i}}{2} - \frac{a \epsilon}{2} (u_{i+1} - u_{i}) = \frac{\{p \}_{i+\frac12} }{2} - \frac{a \epsilon}{2} [u]_{i+\frac12} \label{eq:1dlagprojrelaxacp}
\end{align}
with the subcharacteristic condition $a > \max (\rho c)$, where $c$ is the speed of sound.

The compressive terms are taken into account in a typically Lagrangian way defining
\begin{align}
 L_{i} &:= 1 + \frac{\Delta t}{\Delta x} [u^*]_{i\pm\frac12}
\end{align}
\begin{align}
 \rho^{n+1,-}_{i} &:= \frac{\rho^n_{i}}{L_{i} }  & (\rho u)^{n+1,-}_{i} &:= \frac{(\rho u)^n_{i} - \frac{\Delta t}{\Delta x} \frac{1}{\epsilon^2} [p^*]_{i\pm\frac12} }{L_{i} } \label{eq:intermediatecgkrho}\\
 (\rho v)^{n+1,-}_{i} &:= \frac{(\rho v)^n_{i}}{L_{i} } &
 e^{n+1,-}_{i} &:= \frac{e^n_{i} - \frac{\Delta t}{\Delta x} [u^* p^*]_{i\pm\frac12}}{L_{i} } \label{eq:intermediatecgke}
\end{align}
This intermediate step is not conservative, but the overall scheme is, as will be seen later.

Observe that the compressive terms are rather an ODE-type operator of the form
\begin{align}
 \del_t q + q \cdot D &= 0 \label{eq:odecompression}
\end{align}
where $D$ now is the placeholder for the divergence $\del_x u + \del_y v$. Solving \eqref{eq:odecompression} with the \emph{backward} Euler method yields
\begin{align}
 0 &= \frac{q^{n+1} - q^n}{\Delta t} + q^{n+1} \cdot D & 
 q^{n+1} &= \frac{q^{n}}{1 +  \Delta t  D}
\end{align}
This demonstrates how the denominators in the intermediate step \eqref{eq:intermediatecgkrho}--\eqref{eq:intermediatecgke} can be derived from operator splitting.

Finally, the advective terms are taken into account via a \emph{multiplicative} operator split, with $u^*$ the advecting speed and $q^{n+1,-}$ the advected quantities:
\begin{align}
 q_i^{n+1} &= q^{n+1,-}_{i} L_{i} - \frac{\Delta t}{\Delta x} [f^\text{adv}]_{i\pm\frac12} \label{eq:advstepcgk}\\
 f^\text{adv}_{i+\frac12} &:= \begin{cases}  u^*_{i+\frac12} q^{n+1,-}_{i}  & u^*_{i+\frac12} > 0 \\ u^*_{i+\frac12} q^{n+1,-}_{i+1}  & \text{else} \end{cases}
\end{align}
Observe that the first term on the right hand side of \eqref{eq:advstepcgk} involves $q^{n+1,-}_{i} L_{i}$, rather than $q^{n+1,-}_{i}$ alone, thus undoing the nonconservativity of the intermediate step. The advective flux $f^{\text{adv}} \sim u^* q^{n+1,-}_{i}$ uses the results of the intermediate step directly, but is also conservative. Overall, this gives a conservative scheme 
\begin{align}
 q_i^{n+1} = q_i^n - \frac{\Delta t}{\Delta x} [f]_{i\pm\frac12}
\end{align}
with the following flux (if $u^*_{i+\frac12} > 0$):
\begin{align}
 f_{i+\frac12} &:= \left ( \begin{array}{c} 
                            \displaystyle u^*_{i+\frac12} \frac{\rho^n_{i}}{1 + \frac{\Delta t}{\Delta x} [u^*]_{i\pm\frac12}} \\
                            \displaystyle u^*_{i+\frac12} \frac{(\rho u)^n_{i} - \frac{\Delta t}{\Delta x} \frac{1}{\epsilon^2} [p^*]_{i\pm\frac12} }{1 + \frac{\Delta t}{\Delta x} [u^*]_{i\pm\frac12} } + \frac{p^*_{i+\frac12}}{\epsilon^2}\\
                            \displaystyle u^*_{i+\frac12} \frac{(\rho v)^n_{i}}{1 + \frac{\Delta t}{\Delta x} [u^*]_{i\pm\frac12} }\\
                            \displaystyle u^*_{i+\frac12} \frac{e^n_{i} - \frac{\Delta t}{\Delta x} [u^* p^*]_{i\pm\frac12}}{1 + \frac{\Delta t}{\Delta x} [u^*]_{i\pm\frac12} } +
                            (u^* p^*)_{i+\frac12}
                           \end{array} \right ) \label{eq:cgk1dflux}
\end{align}
The case $u^*_{i+\frac12} < 0$ is analogous. This flux has some similarities to the flux obtained using a relaxation solver for the full Euler system. This is discussed in section \ref{sec:relaxation}.

In \cite{chalons16}, the existence of a discrete entropy inequality for this scheme in the one-dimensional case is demonstrated. 

In multiple spatial dimensions, the multiplication by $L_{ij}$ in \eqref{eq:advstepcgk} needs to make sure that the non-conservativity of the intermediate step is undone. There cannot be several $L_{ij}$, associated to different directions, but there should rather be only one, extended to include contributions from all directions, e.g.
\begin{align}
 L_{ij} &:= 1 + \frac{\Delta t}{\Delta x} [u^*]_{i\pm\frac12,j} + \frac{\Delta t}{\Delta y} [v^*]_{i,j\pm\frac12}
\end{align}
The resulting scheme for multiple spatial dimensions would then be conservative as well, but still dimensionally split and not low Mach number compliant. The aim of the next section is to find an all-speed, truly multi-dimensional extension of this scheme.

\subsection{Multi-dimensional all-speed modification} \label{ssec:multidlagpro}

Observe that in multiple spatial dimensions, starting from \eqref{eq:1dlagprojrelaxacu}--\eqref{eq:1dlagprojrelaxacp}, the dimensionally split acoustic solver would become
\begin{align}
 u^*_{i+\frac12,j} &=  \frac{ \{ u \}_{i+\frac12,j} }{2} - \frac1{2a \epsilon} [p]_{i+\frac12,j} &
 v^*_{i,j+\frac12} &=  \frac{ \{ v \}_{i,j+\frac12} }{2} - \frac1{2a\epsilon} [p]_{i,j+\frac12} \label{eq:ustarlagprodimsplit}\\
 p^*_{i+\frac12,j} &= \frac{\{p \}_{i+\frac12,j} }{2} - \frac{a\epsilon}{2} [u]_{i+\frac12,j}&
 p^*_{i,j+\frac12} &= \frac{\{p \}_{i,j+\frac12} }{2} - \frac{a\epsilon}{2} [v]_{i,j+\frac12} \label{eq:pstarlagprodimsplit}
\end{align}

Usual low Mach number modifications would require a factor of $\epsilon$ in front of $[u]_{i+\frac12,j}$ and $[v]_{i,j+\frac12}$ in \eqref{eq:pstarlagprodimsplit}. Here, it is proposed to replace these derivatives of individual velocity components by the divergence. The suggested multi-dimensional modification of the acoustic step therefore is
\begin{align}
 u^*_{i+\frac12,j} &=  \frac{\{\{ \{ u \}_{i+\frac12} \}\}_{j\pm\frac12} }{8} - \frac1{2a\epsilon} \frac{\{\{ [p]_{i+\frac12} \}\}_{j\pm\frac12}}{4} \label{eq:lagprojmultidfirst}\\
 p^*_{i+\frac12,j} &= \frac{\{\{\{p \}_{i+\frac12} \}\}_{j\pm\frac12} }{8} - \frac{a\epsilon}{2}\left( \frac{\{\{ [u]_{i+\frac12} \}\}_{j\pm\frac12}}{4} + \frac{\Delta x}{\Delta y} \frac{[\{v\}_{i+\frac12}]_{j\pm1} }{4}\right ) \label{eq:lagprojmultidfirstpxdir}\\
 v^*_{i,j+\frac12} &=  \frac{ \{\{ \{ v \}\}_{i\pm\frac12}\}_{j+\frac12} }{8} - \frac1{2a\epsilon} \frac{ [ \{\{ p \}\}_{i\pm\frac12}]_{j+\frac12}}{4}\\
 p^*_{i,j+\frac12} &= \frac{\{\{ \{p \}\}_{i\pm\frac12}\}_{j+\frac12} }{8} - \frac{a\epsilon}{2} \left( \frac{\Delta y}{\Delta x} \frac{\{[u]_{i\pm1} \}_{j+\frac12}}{4} +  \frac{[ \{\{v\}\}_{i\pm\frac12}]_{j+\frac12} }{4}\right )  \label{eq:lagprojmultidfirstp}
\end{align}
Observe that averaging in the perpendicular direction has been added. This is necessary in order for the different divergence discretizations to vanish simultaneously (see Theorem \ref{thm:cgkmultid} below). Of course, these finite difference operators are inspired by the acoustic scheme discussed in section \ref{sec:overviewacoustics}. The seemingly asymmetric divergence operators in \eqref{eq:lagprojmultidfirstpxdir} and \eqref{eq:lagprojmultidfirstp} become symmetric once the flux difference, i.e. $\frac{[p^*]_{i\pm\frac12,j}}{\Delta x}$ and $\frac{[p^*]_{i,j\pm\frac12}}{\Delta y}$, is computed.

In view of the role of the divergence in the compressive terms, the natural multi-dimensional extension of the intermediate step reads
\begin{align}
 L_{ij} &:= 1 + \frac{\Delta t}{\Delta x} [u^*]_{i\pm\frac12,j} + \frac{\Delta t}{\Delta y} [v^*]_{i,j\pm\frac12}\\
 \rho^{n+1,-}_{ij} &:= \frac{\rho^n_{ij}}{L_{ij} } \\
 (\rho u)^{n+1,-}_{ij} &:= \frac{(\rho u)^n_{ij} - \frac{\Delta t}{\Delta x} \frac{1}{\epsilon^2} [p^*]_{i\pm\frac12,j} }{L_{ij} }\\
 (\rho v)^{n+1,-}_{ij} &:= \frac{(\rho v)^n_{ij} - \frac{\Delta t}{\Delta y} \frac{1}{\epsilon^2} [p^*]_{i,j\pm\frac12}  }{L_{ij} }\\
 e^{n+1,-}_{ij} &:= \frac{e^n_{ij} - \frac{\Delta t}{\Delta x} [u^* p^*]_{i\pm\frac12,j} - \frac{\Delta t}{\Delta y} [v^* p^*]_{i,j\pm\frac12}}{L_{ij} } \label{eq:lagprojmultidlast}
\end{align}
The relation between the discrete divergence appearing in the acoustic step and the divergence $\frac{ [u^*]_{i\pm\frac12,j}}{\Delta x} + \frac{ [v^*]_{i,j\pm\frac12}}{\Delta y}$ appearing in the compressive step is clarified later.

Finally, for low Mach number compliance the advective step does not require a truly multi-dimensional extension, although many would be possible (e.g. \cite{colella90,leveque02}). Here, the original dimensionally split version is retained for simplicity.

Formal asymptotic analysis is employed to study the low Mach number limit $\epsilon \to 0$ of this scheme. 

Recall the definition
 \begin{align}
  \mathscr D^{(0)}_{i+\frac12,j+\frac12} := \frac{\{ [u^{(0)}]_{i+\frac12} \}_{j+\frac12} }{2 \Delta x} + \frac{[ \{v^{(0)} \}_{i+\frac12} ]_{j+\frac12} }{2 \Delta y}
 \end{align}

\begin{lemma} \label{lem:discrdiv}
 If $\mathscr D^{(0)}_{i+\frac12,j+\frac12}$ vanishes $\forall i,j$ then
 \begin{enumerate}[i)]
  \item \label{it:derivD}
  \begin{align}
   \left[ \frac{a^{(0)}_{\cdot,j} }{2} \left( \frac{ \{\{ [u^{(0)}]_{\cdot} \}\}_{j\pm\frac12} }{4 \Delta x}  +   \frac{[\{ v^{(0)}\}_{\cdot}]_{j\pm1}}{4 \Delta y}     \right ) \right ]_{i\pm\frac12} = 0 \quad \forall i,j
  \end{align}
  where $a^{(0)}_{i+\frac12,j}$ is any grid function defined on the edges,
  \item \label{it:avgD}
  \begin{align}
   &\frac{\{\{ [ u^{(0)} ]_{i\pm1} \}\}_{j\pm\frac12} }{8 \Delta x}   + \frac{ [ \{\{ v^{(0)} \}\} _{i\pm\frac12} ]_{j\pm1} }{8 \Delta y}  =\\
   &\left[ \frac{\{\{ \{ u^{(0)} \}_{\cdot} \}\}_{j\pm\frac12} }{8 \Delta x}  \right ]_{i\pm\frac12,j} + \left[ \frac{ [ \{\{ v^{(0)} \}\} _{i\pm\frac12} ] }{8 \Delta y} \right ]_{j\pm\frac12} = 0 \quad \forall i,j
  \end{align}
 \end{enumerate}
\end{lemma}
\begin{proof}
 \begin{enumerate}[i)]
 \item 
 \begin{align}
  &\left[ \frac{a^{(0)}_{\cdot,j} }{2} \left( \frac{ \{\{ [u^{(0)}]_{\cdot} \}\}_{j\pm\frac12} }{4 \Delta x}  +   \frac{[\{ v^{(0)}\}_{\cdot}]_{j\pm1}}{4 \Delta y}     \right ) \right ]_{i\pm\frac12} = \left[ \frac{a^{(0)}_{\cdot,j} }{2} \cdot \frac12  \{   \mathscr D^{(0)}  \}_{\cdot, j\pm\frac12}  \right ]_{i\pm\frac12}\\
  &\phantom{mmmm}= \frac{a^{(0)}_{i+\frac12,j} }{2} \cdot \frac12  \{   \mathscr D^{(0)}_{i+\frac12}  \}_{ j\pm\frac12}  - \frac{a^{(0)}_{i-\frac12,j} }{2} \cdot \frac12  \{   \mathscr D^{(0)}_{i-\frac12}  \}_{j\pm\frac12} 
 \end{align}

 \item 
 \begin{align}
  \left[ \frac{\{\{ \{ u^{(0)} \}_{\cdot} \}\}_{j\pm\frac12} }{8 \Delta x}  \right ]_{i\pm\frac12,j} + \left[ \frac{ [ \{\{ v^{(0)} \}\} _{i\pm\frac12} ] }{8 \Delta y} \right ]_{j\pm\frac12} = \frac14 \left\{ \left\{ \mathscr D^{(0)} \right \}_{i\pm\frac12} \right \}_{j\pm\frac12}
 \end{align}

 \end{enumerate}
\end{proof}

The same statement can also be proved using the discrete Fourier transform and a discrete Leibniz rule, see \cite{barsukow18thesis}.

Recall that the low Mach number artefacts identified in \eqref{eq:roeacousticp1x}--\eqref{eq:roeacousticp1y}, and similarly those for the Euler equations (see e.g. \cite{guillard99}), are related to $p^{(1)} = \const$ not being a solution of the discrete limit equations. In other words, the numerical diffusion restricts the discrete limit equations unnecessarily (as has also been elucidated in \cite{barsukow17a} for linear acoustics). For the multi-dimensional extension of the Lagrange Projection scheme one can now show the following statement:

\begin{theorem} \label{thm:cgkmultid}
 Consider the numerical scheme \eqref{eq:lagprojmultidfirst}--\eqref{eq:lagprojmultidlast} for the Euler equations. It is low Mach compliant in the sense that
 \begin{align}
 p^{(0)} &= \const & p^{(1)} &= \const  \label{eq:lagprolimit1}
 \end{align}
 and (discrete counterpart to $\div \vec v^{(0)} = 0$)
 \begin{align}
   \mathscr D^{(0)}_{i+\frac12,j+\frac12} = \frac{\{ [u^{(0)}]_{i+\frac12} \}_{j+\frac12} }{2 \Delta x} + \frac{[ \{v^{(0)} \}_{i+\frac12} ]_{j+\frac12} }{2 \Delta y}  = 0 \quad \forall i,j \label{eq:lagprolimit2}
 \end{align}
 are solutions of the discrete limit equations as $\epsilon \to 0$.
\end{theorem}

\begin{proof}
{$p^{(0)} = \const$ implies $(u^*_{i+\frac12,j})^{(-1)} = 0$.} The $\mathcal O(1/\epsilon^2)$ terms are only present in the momentum equation:
\begin{align}
 [p^*]^{(0)}_{i\pm\frac12,j} &+ (u^*_{i+\frac12})^{(0)} \begin{cases}  \frac{ - \frac{\Delta t}{\Delta x}  [p^*]^{(0)}_{i\pm\frac12,j} }{1 + \frac{\Delta t}{\Delta x} [u^*]^{(0)}_{i\pm\frac12,j} + \frac{\Delta t}{\Delta y} [v^*]^{(0)}_{i,j\pm\frac12}} & u^*_{i+\frac12} > 0 \\ 
		\frac{ - \frac{\Delta t}{\Delta x}  [p^*]^{(0)}_{i+1\pm\frac12,j} }{1 + \frac{\Delta t}{\Delta x} [u^*]^{(0)}_{i+1\pm\frac12,j} + \frac{\Delta t}{\Delta y} [v^*]^{(0)}_{i+1,j\pm\frac12}} & \text{else}
                                                        \end{cases} \nonumber
 \\&- (u^*_{i-\frac12})^{(0)} \begin{cases}  \frac{ - \frac{\Delta t}{\Delta x}  [p^*]^{(0)}_{i-1\pm\frac12,j} }{1 + \frac{\Delta t}{\Delta x} [u^*]^{(0)}_{i-1\pm\frac12,j} + \frac{\Delta t}{\Delta y} [v^*]^{(0)}_{i-1,j\pm\frac12}} & u^*_{i-\frac12} > 0 \\ 
		\frac{ - \frac{\Delta t}{\Delta x}  [p^*]^{(0)}_{i\pm\frac12,j} }{1 + \frac{\Delta t}{\Delta x} [u^*]^{(0)}_{i\pm\frac12,j}  + \frac{\Delta t}{\Delta y} [v^*]^{(0)}_{i,j\pm\frac12}} & \text{else}
                                                        \end{cases}  \nonumber\\&+ \text{perpendicular direction} = 0 \label{eq:order2cgkmultid}
\end{align}

This is indeed solved by $p^{(0)} = \const$ with \eqref{eq:lagprojmultidfirst}--\eqref{eq:lagprojmultidfirstp} 
\begin{align}
 (p^*_{i+\frac12,j})^{(0)} &= \frac{\{\{\{p^{(0)} \}_{i+\frac12} \}\}_{j\pm\frac12} }{8} = p^{(0)}\\
 (p^*_{i,j+\frac12})^{(0)} &= \frac{\{\{ \{p^{(0)} \}\}_{i\pm\frac12}\}_{j+\frac12,j} }{8} = p^{(0)} 
\end{align}

{Using $p^{(0)} = \const$, the remaining} $\mathcal O(1/\epsilon)$ terms are
\begin{align}
 [p^*]^{(1)}_{i\pm\frac12,j} &+ (u^*_{i+\frac12})^{(0)} \begin{cases}  \frac{ - \frac{\Delta t}{\Delta x}  [p^*]^{(1)}_{i\pm\frac12,j} }{1 + \frac{\Delta t}{\Delta x} [u^*]^{(0)}_{i\pm\frac12,j} + \frac{\Delta t}{\Delta y} [v^*]^{(0)}_{i,j\pm\frac12}} & u^*_{i+\frac12} > 0 \\ 
		\frac{ - \frac{\Delta t}{\Delta x}  [p^*]^{(1)}_{i+1\pm\frac12,j} }{1 + \frac{\Delta t}{\Delta x} [u^*]^{(0)}_{i+1\pm\frac12,j} + \frac{\Delta t}{\Delta y} [v^*]^{(0)}_{i+1,j\pm\frac12}} & \text{else}
                                                        \end{cases} \nonumber\\
 &- (u^*_{i-\frac12})^{(0)} \begin{cases}  \frac{ - \frac{\Delta t}{\Delta x}  [p^*]^{(1)}_{i-1\pm\frac12,j} }{1 + \frac{\Delta t}{\Delta x} [u^*]^{(0)}_{i-1\pm\frac12,j} + \frac{\Delta t}{\Delta y} [v^*]^{(0)}_{i-1,j\pm\frac12}} & u^*_{i-\frac12} > 0 \\ 
		\frac{ - \frac{\Delta t}{\Delta x}  [p^*]^{(1)}_{i\pm\frac12,j} }{1 + \frac{\Delta t}{\Delta x} [u^*]^{(0)}_{i\pm\frac12,j}  + \frac{\Delta t}{\Delta y} [v^*]^{(0)}_{i,j\pm\frac12}} & \text{else}
                                                        \end{cases}  \nonumber\\&+ \text{perpendicular direction} =0 
\end{align}

One has from \eqref{eq:lagprojmultidfirst}--\eqref{eq:lagprojmultidfirstp} 
\begin{align}
 (p^*_{i+\frac12,j})^{(1)} &= \frac{\{\{\{p^{(1)} \}_{i+\frac12} \}\}_{j\pm\frac12} }{8} - \frac{a^{(0)}}{2}\left( \frac{\{\{ [u^{(0)}]_{i+\frac12} \}\}_{j\pm\frac12}}{4} + \frac{\Delta x}{\Delta y} \frac{[\{v^{(0)}\}_{i+\frac12}]_{j\pm1} }{4}\right )\\
 (p^*_{i,j+\frac12})^{(1)} &= \frac{\{\{ \{p^{(1)} \}\}_{i\pm\frac12}\}_{j+\frac12,j} }{8} - \frac{a^{(0)}}{2} \left( \frac{\Delta y}{\Delta x} \frac{\{[u^{(0)}]_{i\pm1} \}_{j+\frac12}}{4} +  \frac{[ \{\{v^{(0)}\}\}_{i\pm\frac12}]_{j+\frac12} }{4}\right )
\end{align}

Computing $\frac{[p^*]_{i\pm\frac12,j}}{\Delta x}$, $\frac{[p^*]_{i,j\pm\frac12}}{\Delta y}$, the $\mathcal O(1/\epsilon)$ terms vanish due to Lemma \ref{lem:discrdiv}\ref{it:derivD}.

The $O(1)$ energy equation, using $\del_t e^{(0)} = 0$ yields (without loss of generality for $u^*_{i+\frac12,j}$, $v^*_{i,j+\frac12}$ uniformly positive)
\begin{align*}
 0= (u^*_{i+\frac12,j})^{(0)} \frac{(e^n_{ij})^{(0)} - \frac{\Delta t}{\Delta x} [u^* p^*]^{(0)}_{i\pm\frac12,j} - \frac{\Delta t}{\Delta y} [v^* p^*]^{(0)}_{i,j\pm\frac12}}{1 + \frac{\Delta t}{\Delta x} [u^*]^{(0)}_{i\pm\frac12,j} + \frac{\Delta t}{\Delta y} [v^*]^{(0)}_{i,j\pm\frac12}} +
                            (u^* p^*)^{(0)}_{i+\frac12,j}
 \\-(u^*_{i-\frac12,j})^{(0)} \frac{(e^n_{i-1,j})^{(0)} - \frac{\Delta t}{\Delta x} [u^* p^*]^{(0)}_{i-1\pm\frac12,j}- \frac{\Delta t}{\Delta y} [v^* p^*]^{(0)}_{i-1,j\pm\frac12} }{1 + \frac{\Delta t}{\Delta x} [u^*]^{(0)}_{i-1\pm\frac12,j} + \frac{\Delta t}{\Delta y} [v^*]^{(0)}_{i-1,j\pm\frac12}} -
                            (u^* p^*)^{(0)}_{i-\frac12,j}\\
  +(v^*_{i,j+\frac12})^{(0)} \frac{(e^n_{ij})^{(0)} - \frac{\Delta t}{\Delta x} [u^* p^*]^{(0)}_{i\pm\frac12,j} - \frac{\Delta t}{\Delta y} [v^* p^*]^{(0)}_{i,j\pm\frac12}}{1 + \frac{\Delta t}{\Delta x} [u^*]^{(0)}_{i\pm\frac12,j} + \frac{\Delta t}{\Delta y} [v^*]^{(0)}_{i,j\pm\frac12}} +
                            (v^* p^*)^{(0)}_{i,j+\frac12}
 \\-(v^*_{i,j-\frac12})^{(0)} \frac{(e^n_{i,j-1})^{(0)} - \frac{\Delta t}{\Delta x} [u^* p^*]^{(0)}_{i\pm\frac12,j-1}- \frac{\Delta t}{\Delta y} [v^* p^*]^{(0)}_{i,j-1\pm\frac12} }{1 + \frac{\Delta t}{\Delta x} [u^*]^{(0)}_{i\pm\frac12,j-1} + \frac{\Delta t}{\Delta y} [v^*]^{(0)}_{i,j-1\pm\frac12}} -
                            (v^* p^*)^{(0)}_{i,j-1+\frac12}
\end{align*}

Here, a seemingly different divergence appears. However,
\begin{align}
 \frac{[u^*]^{(0)}_{i\pm\frac12,j}}{\Delta x} + \frac{[v^*]^{(0)}_{i,j\pm\frac12}}{\Delta y} &= \frac{\{\{ [ u^{(0)} ]_{i\pm1} \}\}_{j\pm\frac12} }{8 \Delta x} + \frac{ [\{ \{ v^{(0)} \}\}_{i\pm\frac12}]_{j\pm1} }{8 \Delta y} 
 \nonumber\\&- \frac1{2a^{(0)}} \frac{\{\{ [p^{(1)}]_{i+\frac12} \}\}_{j\pm\frac12}}{4 \Delta x}
 - \frac1{2a^{(0)}} \frac{ [ \{\{ p^{(1)} \}\}_{i\pm\frac12}]_{j+\frac12}}{4 \Delta y}
\end{align}
and thus it vanishes by Lemma \ref{lem:discrdiv}\ref{it:avgD}.

\end{proof}

Note that, although the individual steps of the proof are simple, it crucially depends on the presence of the \emph{correct} discrete divergence in the definition of $p^*_{i+\frac12}$, which satisfies the property of Lemma \ref{lem:discrdiv}. This builds on previous efforts (e.g. \cite{morton01,sidilkover02,jeltsch06,barsukow17a}) to construct this special divergence discretization. In \cite{barsukow17a}, for example, it is proven that a simple central divergence discretization on a 5-point stencil does not admit any statement analogous to Lemma \ref{lem:discrdiv}.

The existence of the discrete counterparts to the limit equations seems the crucial aspect of low Mach number compliance. This is elucidated in the linear case in \cite{barsukow17a}. It is shown there that a scheme that is not low Mach compliant fails to provide a discretization for \emph{all} the limit equations. 
The numerical diffusion restricts the discrete limit equations, such that \eqref{eq:roeacousticp1x}-\eqref{eq:roeacousticp1y}, for example, only implies $p^{(1)} = \const$ for spatially linear velocities. The above theorem shows that the multi-dimensional extension of the Lagrange Projection scheme is free of such a restriction and allows for a constant pressure for all (discrete) divergence-free velocities. This is understood to be the essence of a low Mach compliant scheme. It might in principle additionally be of interest whether \eqref{eq:lagprolimit1}--\eqref{eq:lagprolimit2} are the unique solutions of the limit equations. This does not seem to be related to the question of low Mach compliance. Any additional solutions might point to the existence of spurious modes in numerical simulations -- checkerboards, for example. Such spurious modes are unrelated to the low Mach problem, as they can appear in many contexts, and also for standard low Mach fixes. Such a uniqueness result, however, seems fairly difficult to establish and in the numerical experiments no spurious modes have been observed.

In \cite{chalons16} it is shown that the one-dimensional Lagrange Projection scheme satisfies a discrete entropy inequality as long as the low Mach fix is \emph{not} applied. As the multi-dimensional scheme proposed here reverts to the one-dimensional scheme in one-dimensional situations, one can expect such properties to carry over to the all-speed case more easily than with a low Mach number fix. However, it does not seem possible to immediately generalize the proof of the discrete entropy inequality to the multi-dimensional extension proposed here. To study this, and whether the existence of a discrete entropy inequality entails further conditions on the multi-dimensional all-speed extension, is subject of future work.

A von Neumann stability analysis in multiple spatial dimensions suffers from many technical difficulties, but in view of the fact that the linearization of the proposed method essentially becomes the acoustic scheme of section \ref{sec:overviewacoustics}, the stability result of \ref{ap:stability} can be considered to carry over. The stability of the scheme is confirmed experimentally. It remains stable for CFL numbers up to 1, rather than $\frac12$, which would be the case for a dimensionally split method in two spatial dimensions.

An extension to three dimensions can be performed analogously.

\section{Multi-dimensional modification of a relaxation solver}  \label{sec:relaxation}

\subsection{Review of the one-dimensional method} \label{ssec:relax1d}

Following \cite{chalons10,girardin14}, define the following numerical method, derived as a relaxation solver
\begin{align}
 Q\left( \frac{x}{t}; q_{i}, q_{i+1} \right) := \begin{cases}
    q_i & \frac{x}{t} \leq u^*_{i+\frac12} - c^*_{i+\frac12} \\
    q^*_{i+\frac12,\text L} & u^*_{i+\frac12} - c^*_{i+\frac12}   \leq \frac{x}{t} < u^*_{i+\frac12}  \\
    q^*_{i+\frac12,\text R} & u^*_{i+\frac12} \leq \frac{x}{t} < u^*_{i+\frac12} + c^*_{i+\frac12} \\
    q_{i+1} & \frac{x}{t} \geq u^*_{i+\frac12} + c^*_{i+\frac12} 
                                                       \end{cases} \label{eq:relaxsolverall}
\end{align}

with the intermediate states

\begin{align}
 u^*_{i+\frac12} &= \frac{\{ u \}_{i+\frac12}}{2}  - \frac{1}{2a \epsilon} [p]_{i+\frac12} \label{eq:relax1dustar} \\
 p^*_{i+\frac12} &= \frac{ \{ p \}_{i+\frac12} }{2} - \frac{a \epsilon}{2} [u]_{i+\frac12} \label{eq:relax1dpstar}\\
 \rho^*_{i+\frac12,\text L} &= \frac{\rho_i}{1 + \frac{\rho_i \epsilon}{a} (u^*_{i+\frac12} - u_i)} = 
 \frac{\rho_i}{1 + \rho_i \epsilon \frac{[u]_{i+\frac12}}{2a} - \rho_i\frac{[p]_{i+\frac12}}{2 a^2}} \label{eq:relax1drhoL}\\
 \rho^*_{i+\frac12,\text R} &= \frac{\rho_{i+1}}{1 + \frac{\rho_{i+1} \epsilon}{a} (u_{i+1} - u^*_{i+\frac12}) } =
 \frac{\rho_{i+1}}{1 + \rho_{i+1} \epsilon \frac{[u]_{i+\frac12}}{2 a} + \rho_{i+1}\frac{[p]_{i+\frac12} }{2 a^2} }\\
 \frac{e^*_{i+\frac12,\text L} }{ \rho^*_{i+\frac12,\text L}} &= \frac{e_i}{\rho_i} + \epsilon \frac{p_i u_i - p^*_{i+\frac12} u^*_{i+\frac12}}{a} \\
 \frac{e^*_{i+\frac12,\text R}}{\rho^*_{i+\frac12,\text R}} &= \frac{e_{i+1}}{\rho_{i+1}} + \epsilon \frac{ p^*_{i+\frac12} u^*_{i+\frac12} - p_{i+1} u_{i+1}}{a} \\
 v^*_{i+\frac12, \text{L}} &= v_i\\
 v^*_{i+\frac12, \text{R}} &= v_{i+1} \label{eq:relax1dlast}
\end{align}

For details of the derivation, see \cite{chalons10} and for an introduction to relaxation solvers see e.g. \cite{jin95,bouchut04}. The scalings of the intermediate states have been obtained by performing the rescaling procedure of Section \ref{ssec:euler}.

If $u^*_{i+\frac12} > 0$ then the flux associated with the intermediate state is
\begin{align}
 f_{i+\frac12} = \left ( \begin{array}{c}
                          u^*_{i+\frac12} \frac{\rho_i}{1 + \rho_i \epsilon \frac{[u]_{i+\frac12}}{2a} - \rho_i\frac{[p]_{i+\frac12}}{2 a^2}}\\
                          u^*_{i+\frac12} \left( \frac{\{ u \}_{i+\frac12}}{2}  - \frac{1}{2a \epsilon} [p]_{i+\frac12} \right ) \frac{\rho_i}{1 + \rho_i \epsilon \frac{[u]_{i+\frac12}}{2a} - \rho_i\frac{[p]_{i+\frac12}}{2 a^2}} + \frac{p^*_{i+\frac12}}{\epsilon^2}\\
                          u^*_{i+\frac12} v_i \frac{\rho_i}{1 + \rho_i \epsilon \frac{[u]_{i+\frac12}}{2a} - \rho_i\frac{[p]_{i+\frac12}}{2 a^2}}\\
                          u^*_{i+\frac12} \frac{\rho_i}{1 + \rho_i \epsilon \frac{[u]_{i+\frac12}}{2a} - \rho_i\frac{[p]_{i+\frac12}}{2 a^2}}  \left(  \frac{e_i}{\rho_i} +  \epsilon \frac{p_i u_i - p^*_{i+\frac12} u^*_{i+\frac12}}{a}  \right )    + u^*_{i+\frac12} p^*_{i+\frac12}
                         \end{array} \right ) \label{eq:relax1dflux}
\end{align}

In \cite{chalons10} it is shown that this solver satisfies a discrete entropy inequality.

The similarities between the two schemes go far beyond the observation that the Lagrange-Projection scheme (with the corresponding flux given in \eqref{eq:cgk1dflux}) uses a relaxation solver for the acoustic part. The way \eqref{eq:relax1dflux} treats the compressive and even the advective terms has striking parallels to the results obtained via Lagrange Projection, although the derivations are entirely different.

However, there are also clear differences. The Lagrange Projection scheme has a wider stencil because it uses a \emph{multiplicative} operator split (for an introduction to operator splitting see \cite{holden10}). A more technical observation is that the Lagrange Projection scheme is non-polynomial in $\Delta t$. One observes that the denominators (which account for the compressive terms in \eqref{eq:euler1}--\eqref{eq:euler4}) are different and in \eqref{eq:relax1dflux} do not involve $\Delta t$.

The relaxation scheme \eqref{eq:relaxsolverall} seems easier to implement than the two-step Lagrange Projection scheme \eqref{eq:intermediatecgkrho}--\eqref{eq:intermediatecgke}, \eqref{eq:advstepcgk}. Maybe this is what led \cite{girardin14} to also consider a scheme very similar to \eqref{eq:relaxsolverall}, and to study its possible low Mach number modification. However, it was found therein that modifying the one-dimensional scheme \eqref{eq:relaxsolverall} in the usual way of a low Mach fix, i.e. by multiplying the second derivatives of the velocity by a function $\mathcal O(\epsilon)$, spoiled the stability of the one-dimensional scheme already. Therefore this relaxation method has not been considered further in \cite{girardin14} and the success of the simple low Mach fix when applied to the Lagrange Projection method was attributed to the additional diffusion associated with the predictor step.

The multi-dimensional all-speed modification of the relaxation solver \eqref{eq:relaxsolverall} presented in this paper leaves the one-dimensional scheme untouched, and results in a stable method, as is seen next. 
Here the superiority of the approach of multi-dimensional modification becomes apparent.

\subsection{Multi-dimensional all-speed modification} \label{ssec:relaxmultid}

Replace again equations \eqref{eq:relax1dustar}--\eqref{eq:relax1dpstar} by
\begin{align}
 p^*_{i+\frac12,j} &= \frac{ \{\{ \{ p \}_{i+\frac12} \}\}_{j\pm\frac12}}{8} - \frac{a \epsilon}{2} \left( \frac{ \{\{ [u]_{i+\frac12} \}\}_{j\pm\frac12} }{4}  +  \frac{\Delta x}{\Delta y} \frac{[\{ v\}_{i+\frac12}]_{j\pm1}}{4}     \right ) \label{eq:relax2dpstar}\\
 u^*_{i+\frac12,j} &= \frac{\{\{ \{ u \}_{i+\frac12}\}\}_{j\pm\frac12}}{8}  - \frac{1}{2a \epsilon} \{\{ [p]_{i+\frac12} \}\}_{j\pm\frac12} \label{eq:relax2dustar}
\end{align}
Instead of \eqref{eq:relax1drhoL} it seems natural to use
\begin{align}
 \rho^*_{i+\frac12,j,\text L} &= 
 \frac{\rho_{ij}}{1 + \frac{\rho_{ij} \epsilon}{2a} \left( \frac{\{\{[u]_{i+\frac12}\}\}_{j\pm\frac12}}{4} + \frac{\Delta x}{\Delta y} \frac{[\{  v\}_{i+\frac12}]_{j\pm1}}{4} \right ) - \frac{\rho_{ij}}{ 2 a^2 } \frac{\{\{ [p]_{i+\frac12} \}\}_{j\pm\frac12}}{4}} \label{eq:relax2drhoL}
\end{align}
and similarly for $\rho^*_{i+\frac12,j,\text R}$. The other terms can remain unchanged.

\begin{theorem}
 Consider the numerical scheme \eqref{eq:relaxsolverall}--\eqref{eq:relax1dlast} with the modifications \eqref{eq:relax2dpstar}--\eqref{eq:relax2drhoL}. It is low Mach compliant in the sense that
 \begin{align}
 p^{(0)} &= \const & p^{(1)} &= \const \label{eq:pressureconstrelax2d}
 \end{align}
 and (discrete counterpart to $\div \vec v^{(0)} = 0$)
 \begin{align}
   \mathscr D^{(0)}_{i+\frac12,j+\frac12} = \frac{\{ [u^{(0)}]_{i+\frac12} \}_{j+\frac12} }{2 \Delta x} + \frac{[ \{v^{(0)} \}_{i+\frac12} ]_{j+\frac12} }{2 \Delta y}  = 0 \quad \forall i,j
 \end{align}
 are solutions of the discrete limit equations as $\epsilon \to 0$.
\end{theorem}

\begin{proof}
 Collect the terms $\mathcal O(1/\epsilon^2)$ in the momentum equation:
 \begin{align}
   [p^*]^{(0)}_{i\pm\frac12} = \frac{ \{\{ [ p^{(0)} ]_{i\pm1} \}\}_{j\pm\frac12}}{8} = 0
 \end{align}
 Similarly, the term $\mathcal O(1/\epsilon)$:
 \begin{align}
  \frac{[p^*]^{(1)}_{i\pm\frac12}}{\Delta x} &= \frac{ \{\{ [ p^{(1)} ]_{i\pm1} \}\}_{j\pm\frac12}}{8 \Delta x} - \left[ \frac{a^{(0)}_{\cdot,j} }{2} \left( \frac{ \{\{ [u^{(0)}]_{\cdot} \}\}_{j\pm\frac12} }{4 \Delta x}  +   \frac{[\{ v^{(0)}\}_{\cdot}]_{j\pm1}}{4 \Delta y}     \right ) \right ]_{i\pm\frac12}
  \end{align}
  This expression vanishes according to Lemma \ref{lem:discrdiv}\ref{it:derivD}.
  
  Consider now the $\mathcal O(1)$ energy equation, assuming without loss of generality $u^*_{i+\frac12} > 0$ $\forall i$:
  \begin{align*}
   (u^*_{i+\frac12,j})^{(0)} \rho_i^{(0)} \frac{e_i^{(0)}}{\rho_i^{(0)}} + (u^*_{i+\frac12})^{(0)} (p^*_{i+\frac12})^{(0)} - (u^*_{i-\frac12,j})^{(0)} \rho_{i-1}^{(0)} \frac{e_{i-1}^{(0)}}{\rho_{i-1}^{(0)}} + (u^*_{i-\frac12})^{(0)} (p^*_{i-\frac12})^{(0)} \\ + \text{terms in $y$ direction} =
   \left( \frac{ [u^*]^{(0)}_{i\pm\frac12,j}}{\Delta x} + \frac{[v^*]^{(0)}_{i,j\pm\frac12}}{\Delta y} \right ) \left( e^{(0)} + p^{(0)} \right ) 
  \end{align*}
 This term vanishes by Lemma \ref{lem:discrdiv}\ref{it:avgD} and \eqref{eq:pressureconstrelax2d}.
\end{proof}

The one-dimensional relaxation solver \eqref{eq:relax1dflux} has been endowed with a simple low Mach fix in \cite{girardin14}, but the resulting scheme was not found stable even in one dimension. The multi-dimensional all regime modification presented here does not change properties of the scheme in one dimension, and it has also been found stable in multiple dimensions. The all-speed multi-dimensional extension of a one-dimensional scheme, as presented in this paper, thus seems to be more attractive than a low Mach fix applied to the one-dimensional scheme.
An extension to three dimensions can be performed analogously.

{For the truly multi-dimensional discrete operators, boundary conditions (e.g. periodic, zero-gradient, inflow) can be straightforwardly implemented via ghost cells. They do not require more ghost cells than the dimensionally split versions, because the additionally involved corner cells already exist in any implementation of a Cartesian grid with ghost cells. Also, the truly multi-dimensional schemes are conservative, such that wherever boundaries are implemented by specifying the flux, nothing changes (e.g. wall boundaries).}

\section{Numerical examples} \label{sec:numerical}

\subsection{Incompressible vortex}

\begin{figure}[h]
\centering
\includegraphics[width=0.48\textwidth]{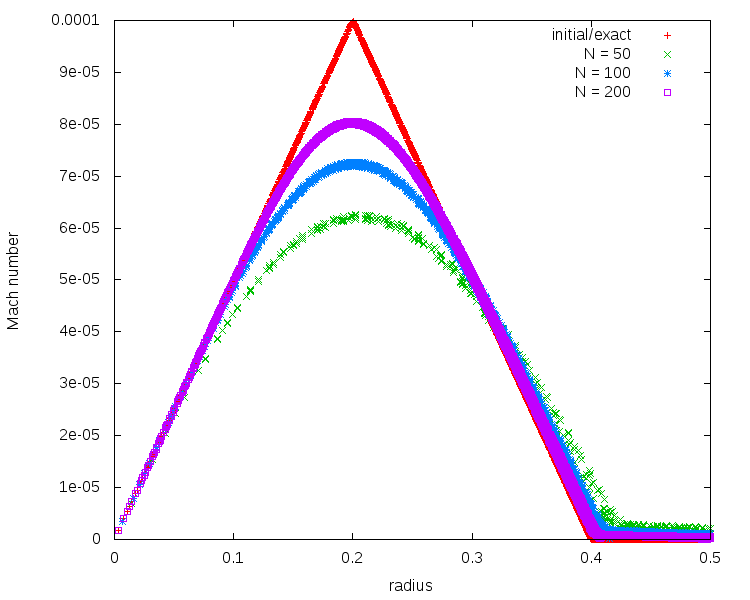}\hfill\includegraphics[width=0.48\textwidth]{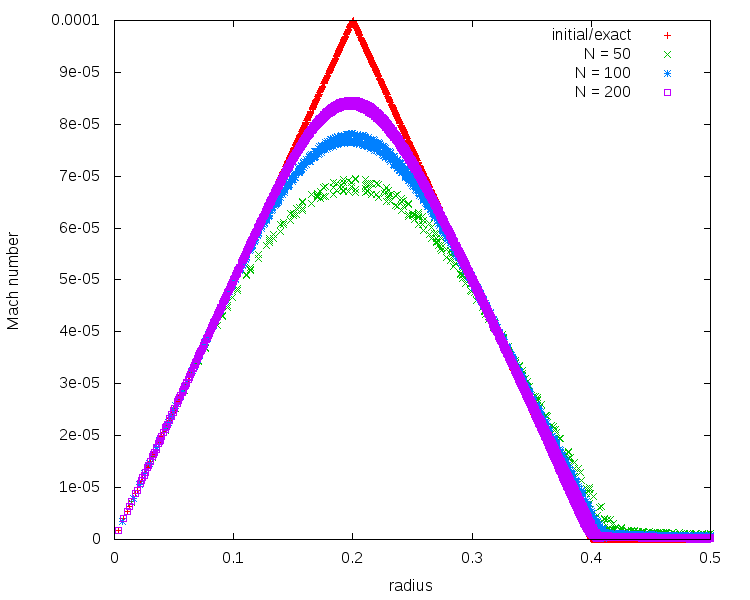}
\caption{{Stationary Gresho vortex \cite{gresho90,barsukow16} solved on the domain $[0,1]^2$ using CFL $= 0.9$ and periodic boundaries, until $t=1$. \emph{Left}: Multi-d extension of the Lagrange-Projection method from section \ref{ssec:multidlagpro}. \emph{Right}: Multi-d extension of a relaxation solver from section \ref{ssec:relaxmultid}. The setup uses a fixed $\epsilon = 10^{-4}$ and is solved on different grids ($50\times50$, $100\times100$, $200\times200$). A radial plot of the Mach number is shown. The figure demonstrates how the numerical solution converges to the exact one upon grid refinement.}}
\label{fig:greshoradialconvergence}
\end{figure}

\begin{figure}
 \centering
 \includegraphics[width=0.3\textwidth]{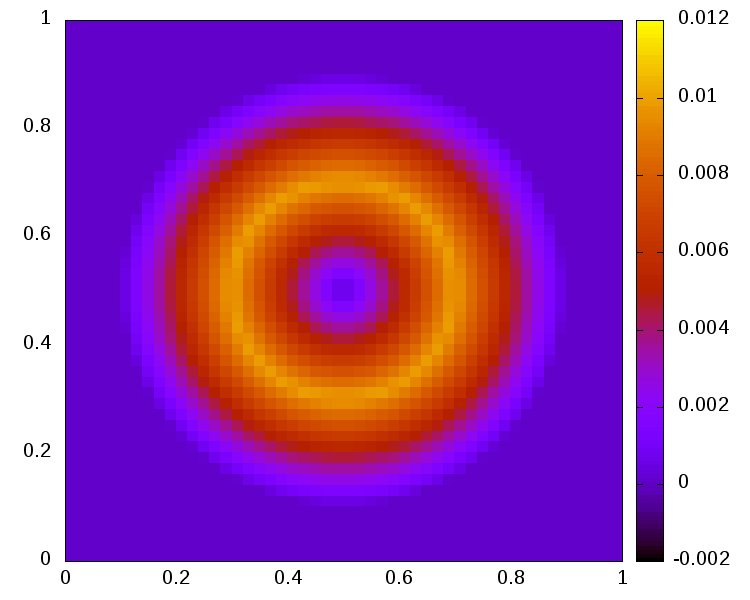} \hfill \includegraphics[width=0.3\textwidth]{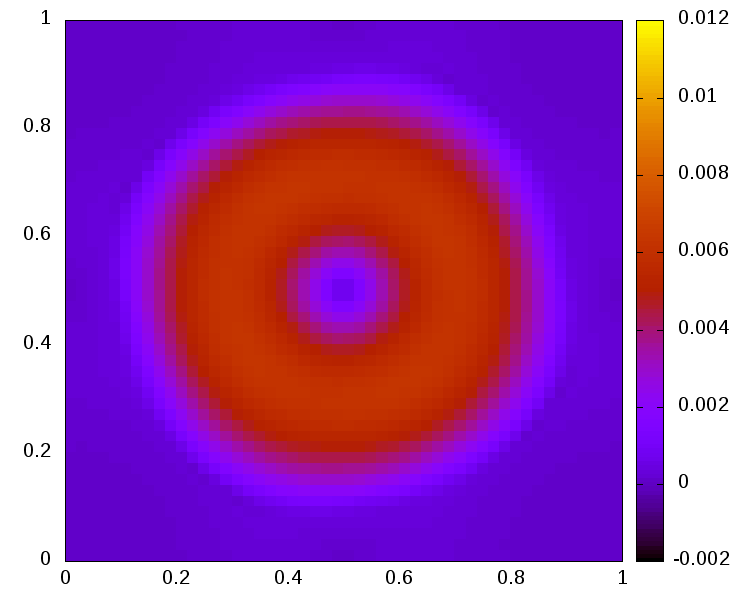} \hfill \includegraphics[width=0.3\textwidth]{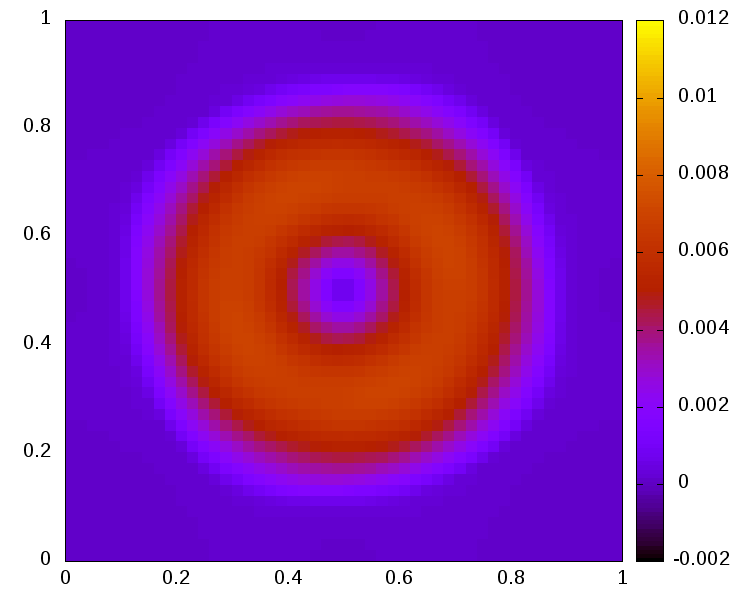}\\
 \includegraphics[width=0.3\textwidth]{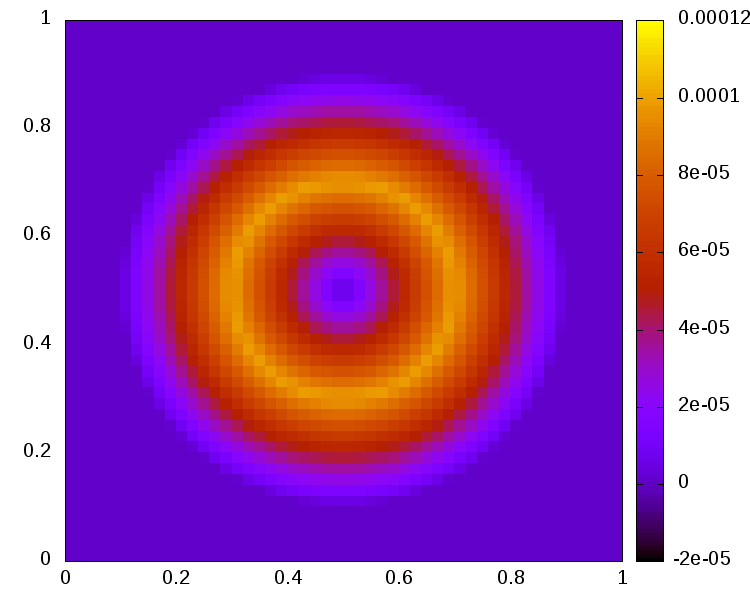} \hfill \includegraphics[width=0.3\textwidth]{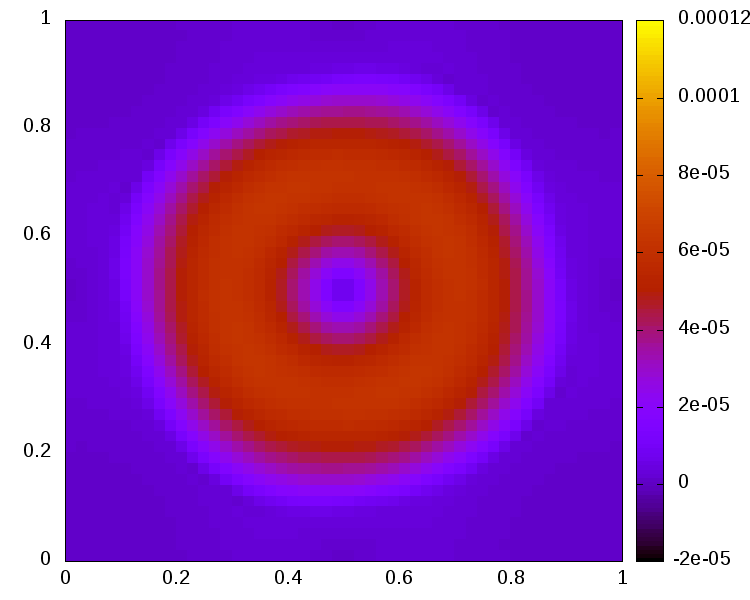}  \hfill  \includegraphics[width=0.3\textwidth]{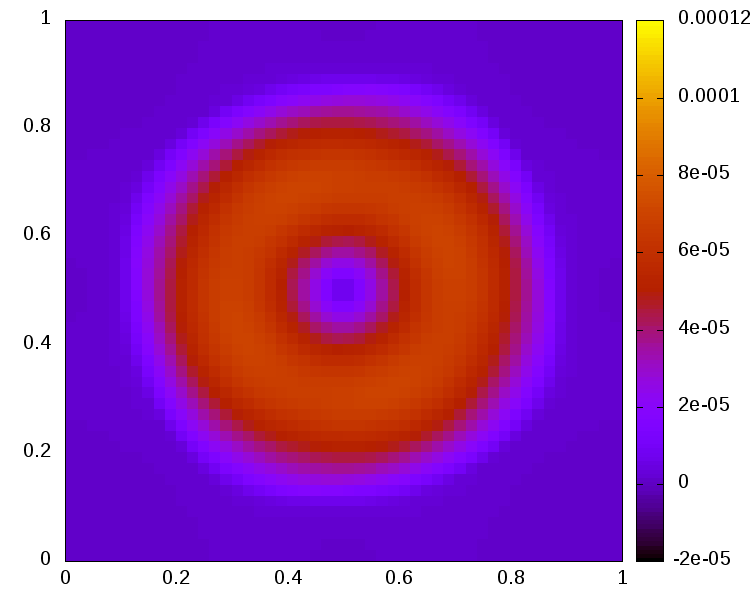}\\
 \includegraphics[width=0.3\textwidth]{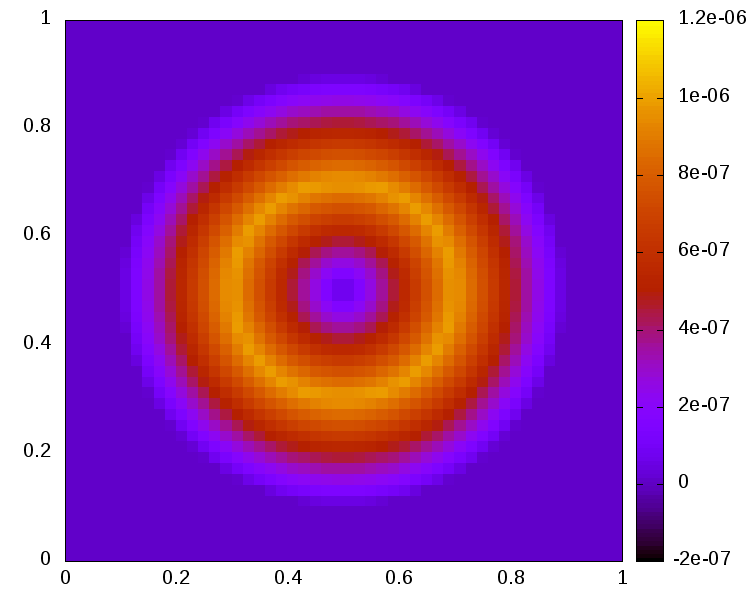} \hfill \includegraphics[width=0.3\textwidth]{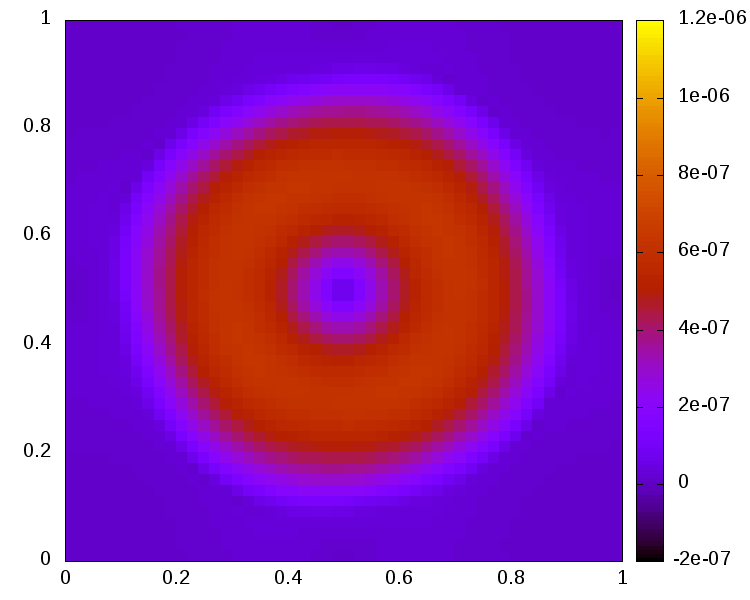} \hfill \includegraphics[width=0.3\textwidth]{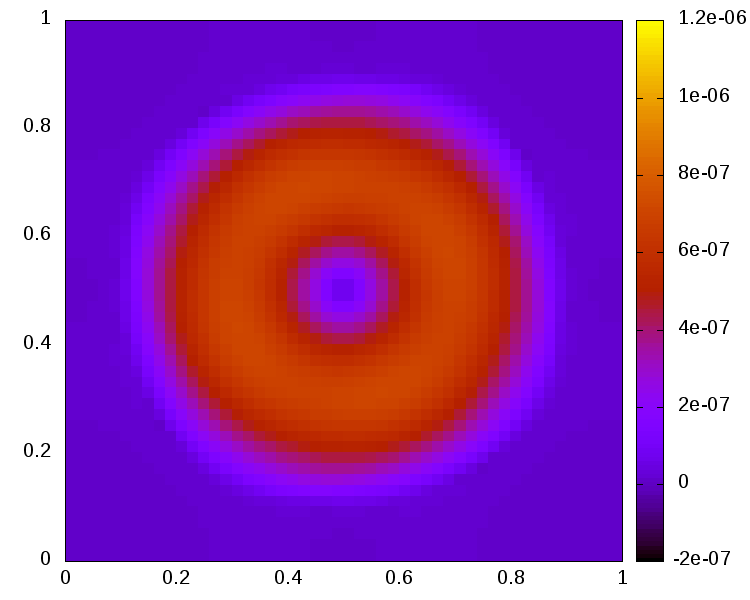}
 \caption{Setup as in Fig. \ref{fig:greshoradialconvergence}, solved on a $50 \times 50$ grid. Color coded is the local Mach number. \emph{Left}: Initial (and exact). \emph{Center}: Multi-d extension of the Lagrange-Projection method from section \ref{ssec:multidlagpro}. \emph{Right}: Multi-d extension of a relaxation solver from section \ref{ssec:relaxmultid}. \emph{Top to bottom}: Numerical solutions at $t=1$ for $\epsilon = 10^{-2}, 10^{-4}, 10^{-6}$.}
 \label{fig:greshomultid}
\end{figure}

The multi-dimensional all-speed extensions of the one-dimensional schemes discussed above shall be exemplified in this section on a number of setups. Figure \ref{fig:greshomultid} shows the stationary Gresho vortex (see \cite{gresho90,barsukow16}). By modifying the background pressure $p_0 := \frac{1}{\gamma \epsilon^2} - \frac12$, the sound speed is changed and thus the Mach number of the vortex. In particular, $p_0$ is defined such that the maximum local Mach number of the setup is just $\epsilon$. The rest of the setup is as follows:
\begin{align}
  \rho &= 1  \qquad  v_\phi = \begin{cases}
                         	5 r & r < 0.2  \\ 2-5 r & 0.2 \leq  r < 0.4 \\ 0 &\text{else} \end{cases}
  \\p &=  \begin{cases} p_0 + 12.5 r^2 & r < 0.2 \\ p_0 + 4 \ln(5r) + 4 - 20r + 12.5 r^2 & 0.2 \leq  r < 0.4 \\ p_0 + 4 \ln 2 - 2 & \text{else}
                        \end{cases}
\end{align}
{Results demonstrating the convergence of the numerical solution upon grid refinement for a fixed $\epsilon$ are shown in Figure \ref{fig:greshoradialconvergence}. The multi-step nature of the Lagrange Projection scheme might be responsible for the slightly larger diffusion observed for this method.}

\begin{figure}[h]
\centering
\includegraphics[width=0.3\textwidth]{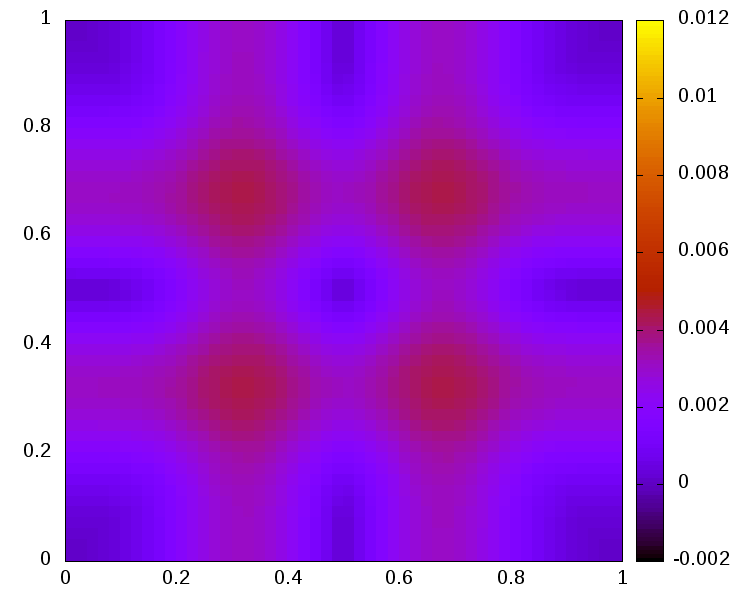}
\caption{Setup as in Fig. \ref{fig:greshomultid}, solved using the non-low Mach, dimensionally split relaxation solver from section \ref{ssec:relax1d} using CFL$=0.45$ on a $50\times50$ grid. Here, $\epsilon = 0.01$.}
\label{fig:gresho1d}
\end{figure}

Figure \ref{fig:greshomultid} shows setups with $ \epsilon = 10^{-2}$, $10^{-4}$, $10^{-6}$ on the same grid. In all cases, the results after $t=1$ ($\sim$ one revolution period) look identically, and the vortex is still recognizable. The numerical diffusion seems independent of $\epsilon$ in the limit. For comparison, the same setup solved with the Lagrange Projection scheme without any low Mach fix is shown in Figure \ref{fig:gresho1d}. As described in \cite{girardin14}, the relaxation solver does not admit a stable low Mach fix, so no comparison can be provided. For results on the Lagrange Projection scheme endowed with a one-dimensional low Mach fix, see \cite{chalons16}.

\begin{figure}[h]
\centering
\includegraphics[width=0.45\textwidth]{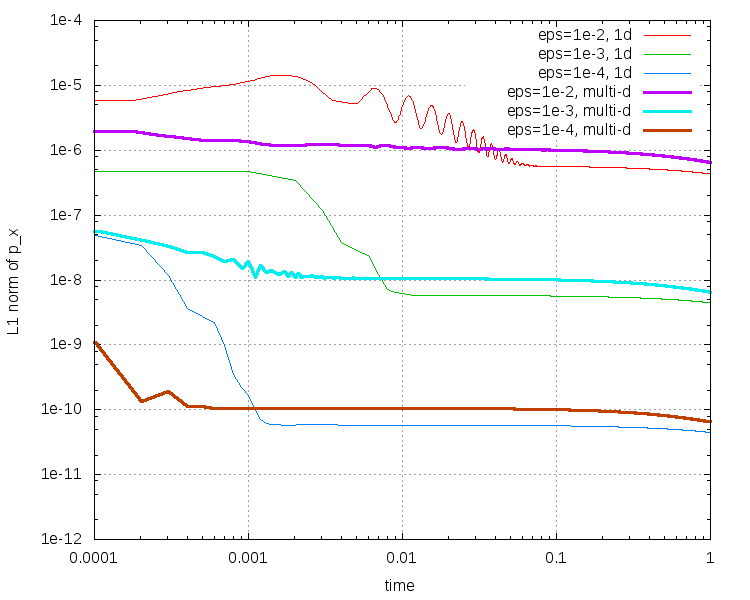}\hfill
\includegraphics[width=0.45\textwidth]{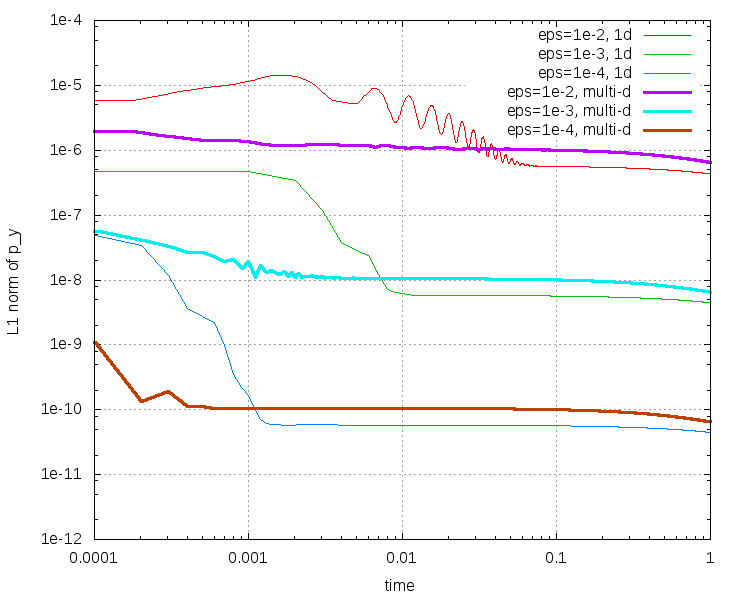}\\
\includegraphics[width=0.45\textwidth]{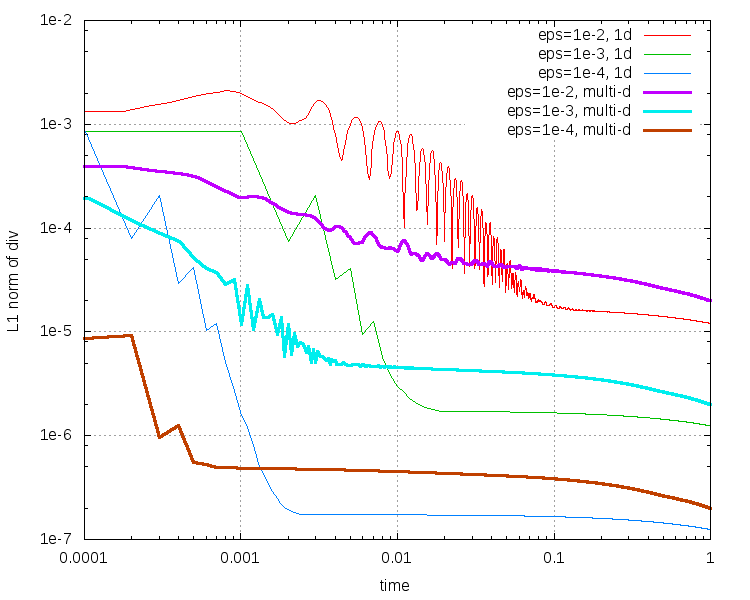}\hfill
\includegraphics[width=0.45\textwidth]{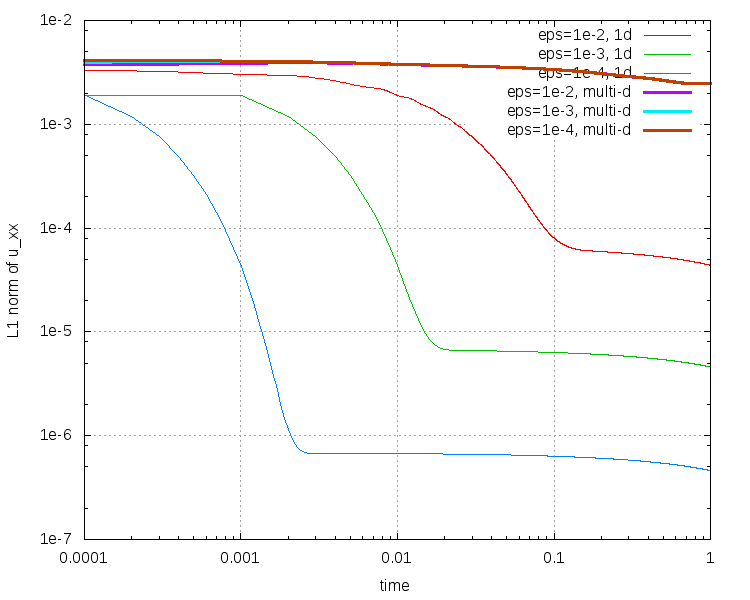}
\caption{{Setup as in Fig. \ref{fig:greshomultid}, using CFL$=0.45$ on a $50\times50$ grid with periodic boundaries for $\epsilon = 10^{-2}, 10^{-3}, 10^{-4}$. \emph{Top}: $\ell_1$-norm of the discrete pressure gradient, i.e. $\frac{1}{N}\sum_{ij}|[\epsilon^2 p]_{i\pm1,j}/2|$ (\emph{left}), $\frac{1}{N}\sum_{ij}|[\epsilon^2 p]_{i,j\pm1}/2|$, (\emph{right}). \emph{Bottom}: $\ell_1$-norm of the discrete divergence, i.e. $\Delta x \frac{1}{N}\sum_{ij}|\mathscr D_{i+\frac12,j+\frac12}|$ (\emph{left}), $\ell_1$-norm of the second derivative of $u$, i.e. $\frac{1}{N}\sum_{ij}|[[u]]_{i\pm\frac12,j}|$. Here, $N$ is the number of cells in the grid. Thick lines show results using the multi-d extension of the Lagrange-Projection method from section \ref{ssec:multidlagpro}, thin lines show its dimensionally split, non-low Mach counterpart.}}
\label{fig:gresholowmachanalysissemilag}
\end{figure}

\clearpage

To demonstrate the correct limit of the numerical solution, Figure \ref{fig:gresholowmachanalysissemilag} shows the relevant quantities as functions of time for different values of $\epsilon$. Note that in this setup, $p = \frac{\const}{\epsilon^2}  + \mathcal O(1)$. At the same time, by definition \eqref{eq:rescalingpressurevelocity} $p = \epsilon^{\mathfrak c} (p^{(0)} + \epsilon p^{(1)} + \mathcal O(\epsilon^2))$ with $\mathfrak c = -2$. This is why in Figure \ref{fig:gresholowmachanalysissemilag}, $\epsilon^2 \nabla p = \nabla p^{(1)} + \mathcal O(\epsilon)$ is plotted. In the literature, $\nabla p$ sometimes is instead divided by the leading constant $\frac{\const}{\epsilon^2}$, which amounts to the same. 

\begin{figure}[h]
\centering
\includegraphics[width=0.75\textwidth]{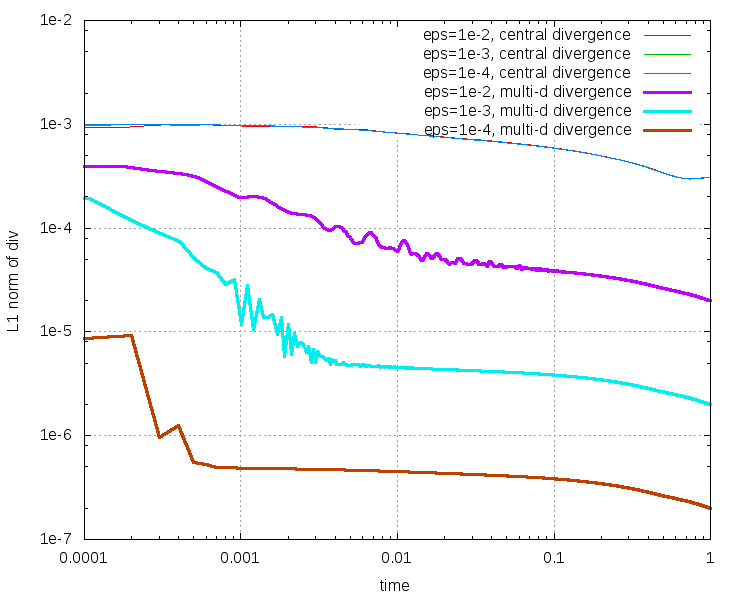}
\caption{{Setup as in Fig. \ref{fig:gresholowmachanalysissemilag} showing results obtained using the multi-d extension of a relaxation solver from section \ref{ssec:relaxmultid}. Thick lines show the evolution of the multi-dimensional divergence $\Delta x \frac{1}{N}\sum_{ij}|\mathscr D_{i+\frac12,j+\frac12}|$, thin lines show the evolution of the central divergence $\frac{1}{N}\sum_{ij}|[u]_{i\pm1,j}/2 + [v]_{i,j\pm1}/2|$.}}
\label{fig:gresholowmachanalysissemilagdivcomparison}
\end{figure}

Concentrate first on the thick lines showing the performance of the multi-dimensional extension of the Lagrange-Projection method from section \ref{ssec:multidlagpro}. As expected, the discrete gradient $\nabla p \in \mathcal O(\epsilon^2)$, and the discrete divergence $\mathscr D \in \mathcal O(\epsilon)$. The discrete second derivative of $u$ does not show dependence on $\epsilon$, as the low Mach number limit of the Euler equations does not enforce any conditions on second derivatives of velocity components. Note that the time evolution of some other discretization of the divergence does not, in general, show the same behaviour, and in particular is not necessarily $\mathcal O(\epsilon)$ (see Fig. \ref{fig:gresholowmachanalysissemilagdivcomparison}). This is clear, because the difference of two discretizations of the divergence must be a discretization of at least a second derivative by consistency, which does not decay to zero as $\epsilon \to 0$.

The thin lines in the plots show the behaviour of the non-low Mach dimensionally split Lagrange-Projection method. Interestingly, its behaviour is more complex. Concentrate first on short times, i.e. on times $t \in \mathcal O(\epsilon)$. For short times, one observes the colloquially well-known ``failure'' of this scheme, characterized by a discrete pressure gradient $\nabla p \in O(\epsilon)$ instead of $\mathcal O(\epsilon^2)$, and the discrete divergence $\mathscr D \in \mathcal O(1)$ instead of $\mathcal O(\epsilon)$. 

However, for long times, it is obvious that these ``failures'' disappear, i.e. the gradient of pressure \emph{does} become $\mathcal O(\epsilon^2)$ and the divergence \emph{does} decay as $\mathcal O(\epsilon)$. So what is wrong about this solution? Consider equation \eqref{eq:roeacousticp1x} that characterizes the limit states of a non-low Mach scheme\footnote{Equation \eqref{eq:roeacousticp1x} refers to the acoustic equations, but the results are very similar for the scheme given by \eqref{eq:cgk1dflux} for the Euler equations.}. As was mentioned, this equation does not imply $\nabla p^{(1)} \neq 0$. However, if $\nabla p^{(1)} = 0$, then the discrete second derivative of $u$ must vanish as well. And indeed, this is clearly seen in Figure \ref{fig:gresholowmachanalysissemilag} (\emph{bottom right}), where (discretely) $\del_x^2 u \in \mathcal O(\epsilon)$. This is as unacceptable as $\nabla p^{(1)} \neq 0$. So in fact, the behaviour of a non-low Mach scheme is, for this simulation, two-fold: for short times, the limit conditions $\nabla p \in \mathcal O(\epsilon^2)$, $\div \vec v \in \mathcal O(\epsilon)$ are violated, and for long times the limit conditions are fulfilled, but at the cost of additional relations (such as $\del_x^2 u \in \mathcal O(\epsilon)$) which have nothing to do with the limit. These latter are the reason for the visible deterioration of the solution shown in Figure \ref{fig:gresho1d}. Note that in the linear case, this behaviour has already been explained in \cite{barsukow17a}. The existence of a transition time $\mathcal O(\epsilon)$ is not surprising for a singular limit. In Section \ref{ap:asymptotic}, as an example, this time scale is derived for a toy problem.

\begin{figure}[h]
\centering
\includegraphics[width=0.45\textwidth]{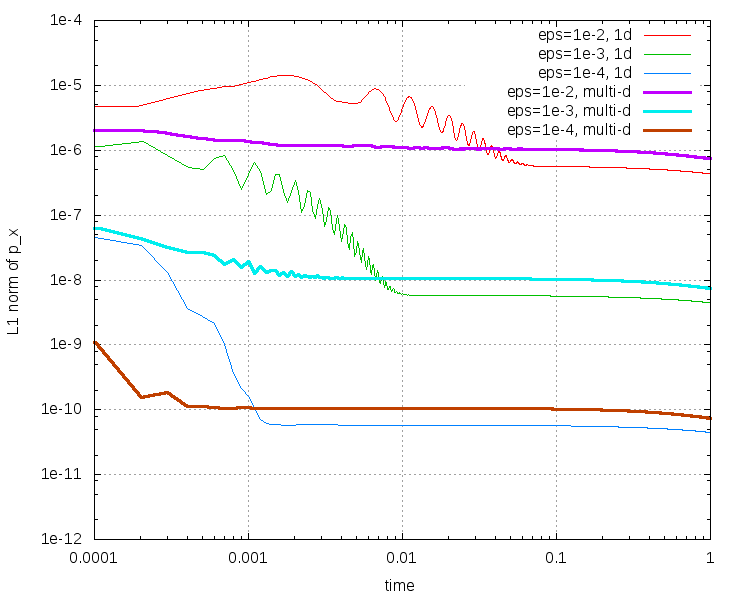}\hfill
\includegraphics[width=0.45\textwidth]{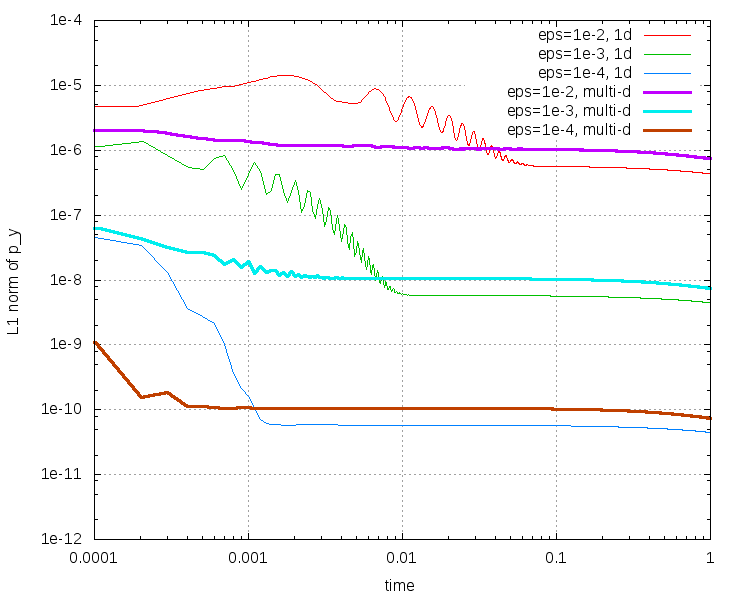}\\
\includegraphics[width=0.45\textwidth]{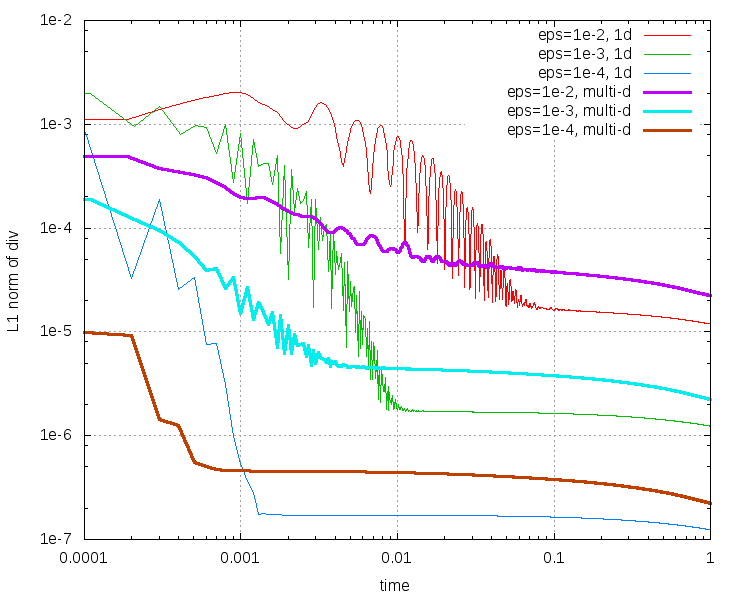}\hfill
\includegraphics[width=0.45\textwidth]{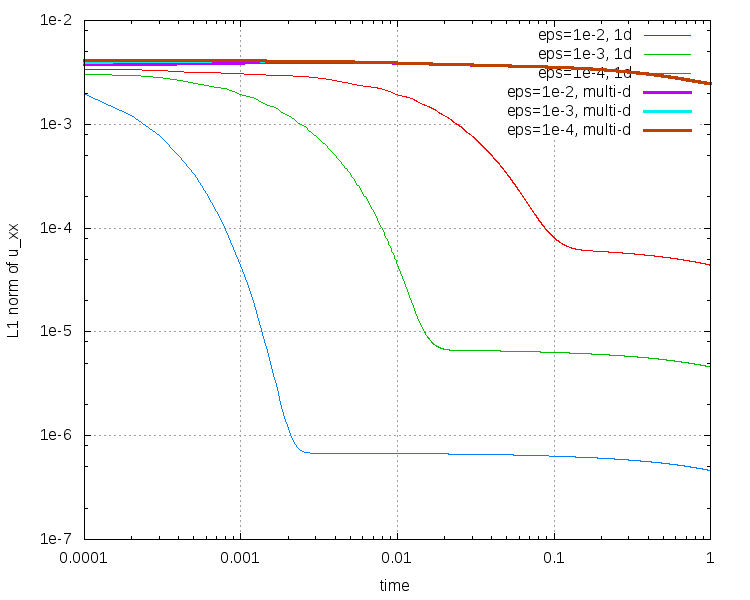}
\caption{{Setup as in Fig. \ref{fig:gresholowmachanalysissemilag}. Thick lines show results using the multi-d extension of a relaxation solver from section \ref{ssec:relaxmultid}, thin lines show its dimensionally split, non-low Mach counterpart.}}
\label{fig:gresholowmachanalysisonestep}
\end{figure}

Figure \ref{fig:gresholowmachanalysisonestep} demonstrates the analogous situation for the multi-d extension of a relaxation solver from section \ref{ssec:relaxmultid}.

\clearpage

\subsection{Incompressible vortex and a sound wave}

\begin{figure}[h]
\centering
\includegraphics[width=0.23\textwidth]{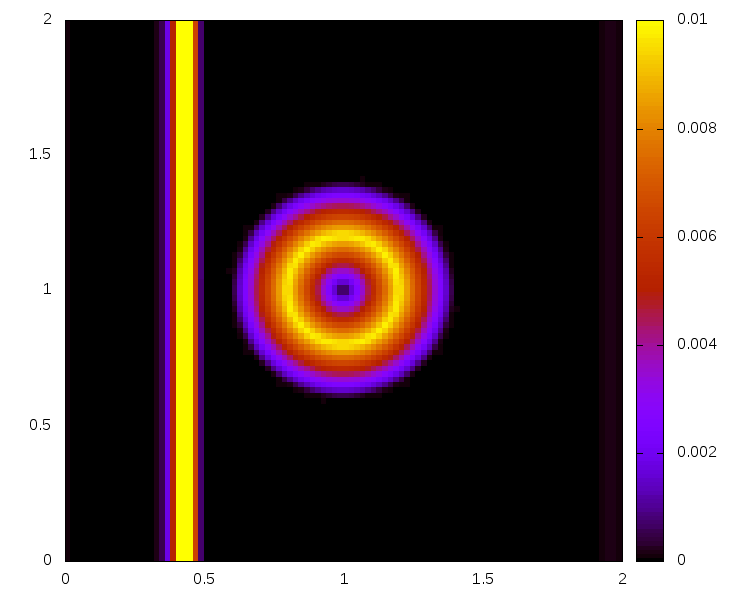}\hfill
\includegraphics[width=0.23\textwidth]{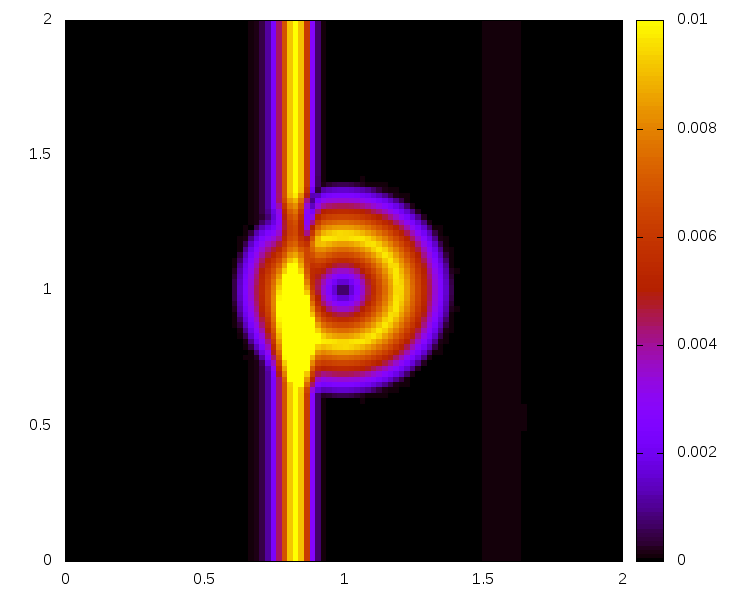}\hfill
\includegraphics[width=0.23\textwidth]{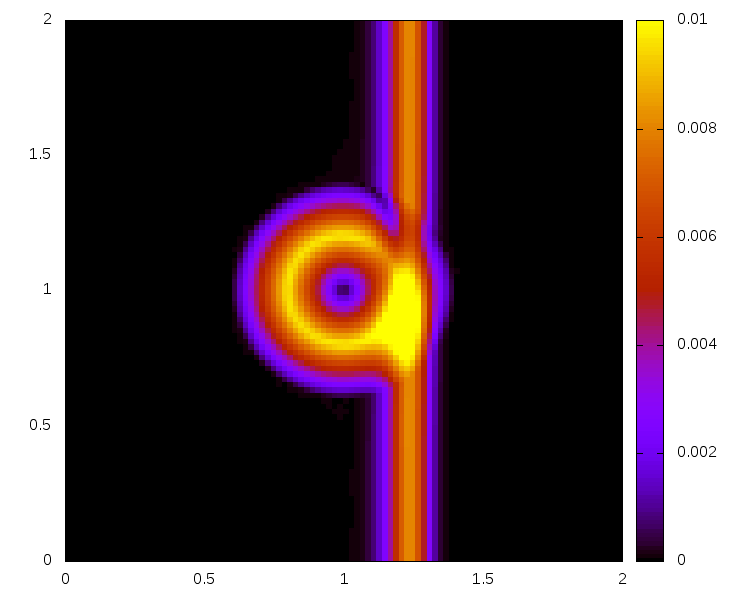}\hfill
\includegraphics[width=0.23\textwidth]{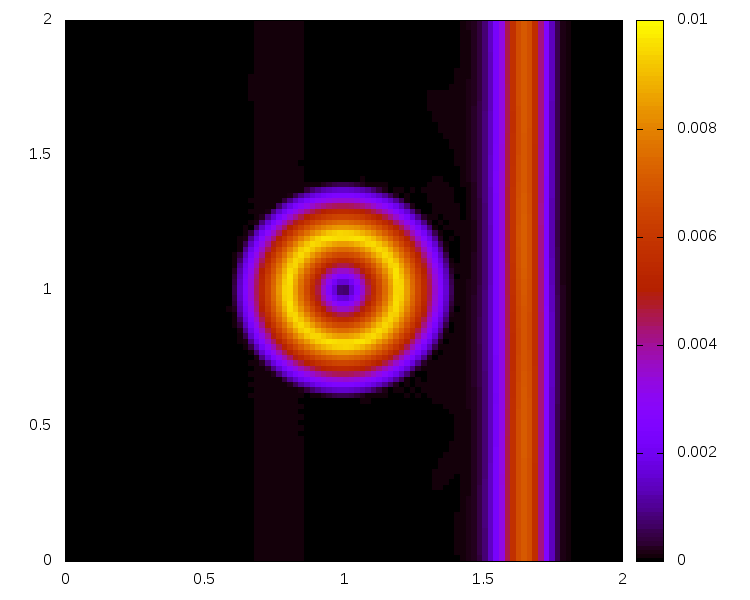}\\
\includegraphics[width=0.23\textwidth]{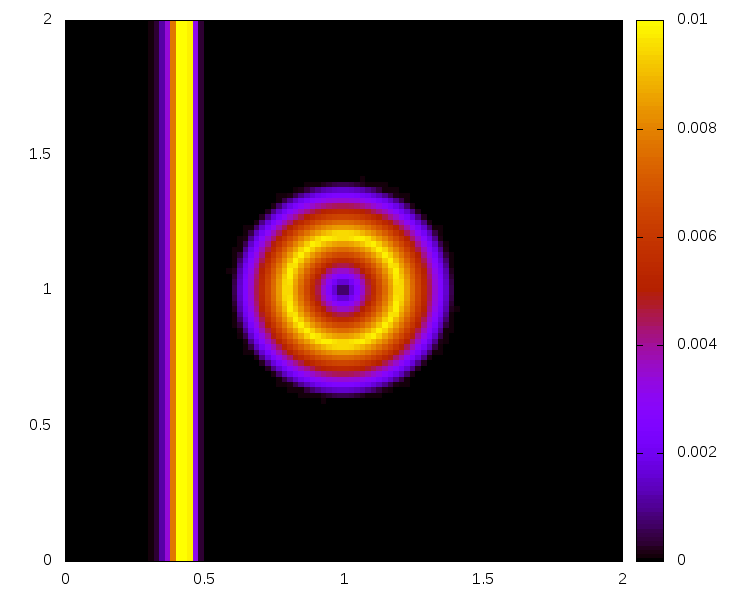}\hfill
\includegraphics[width=0.23\textwidth]{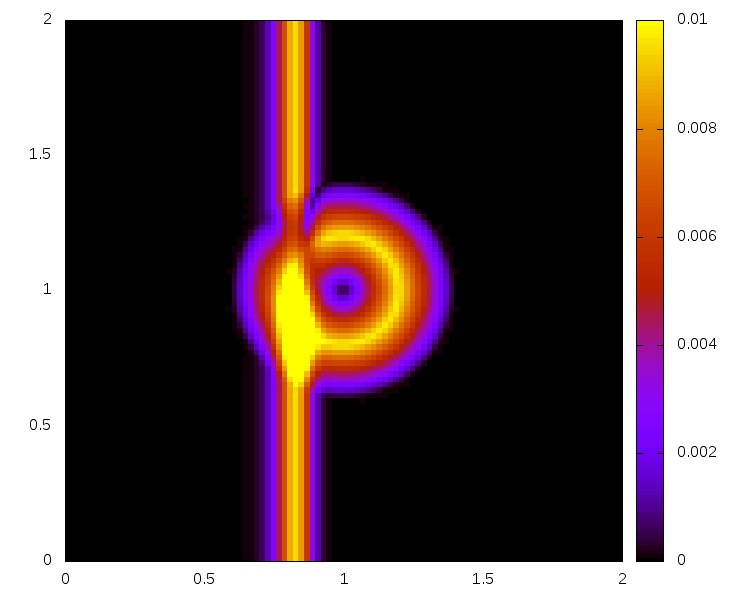}\hfill
\includegraphics[width=0.23\textwidth]{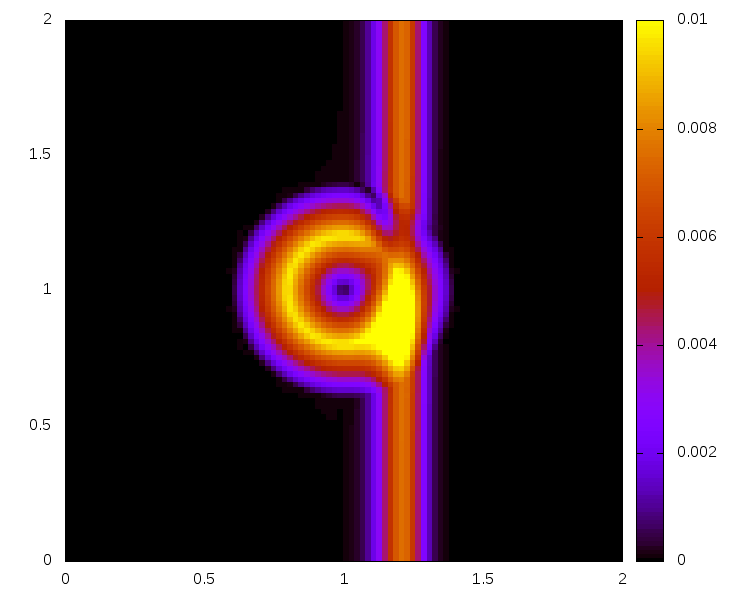}\hfill
\includegraphics[width=0.23\textwidth]{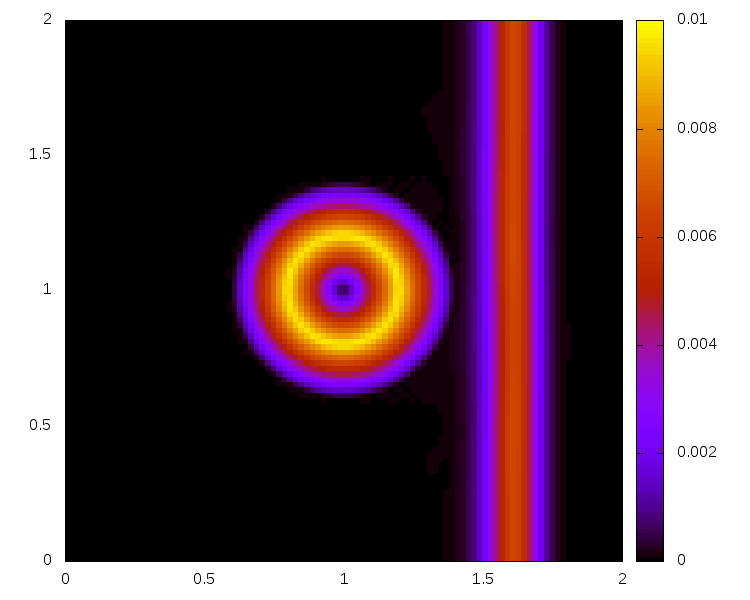}\\
\includegraphics[width=0.23\textwidth]{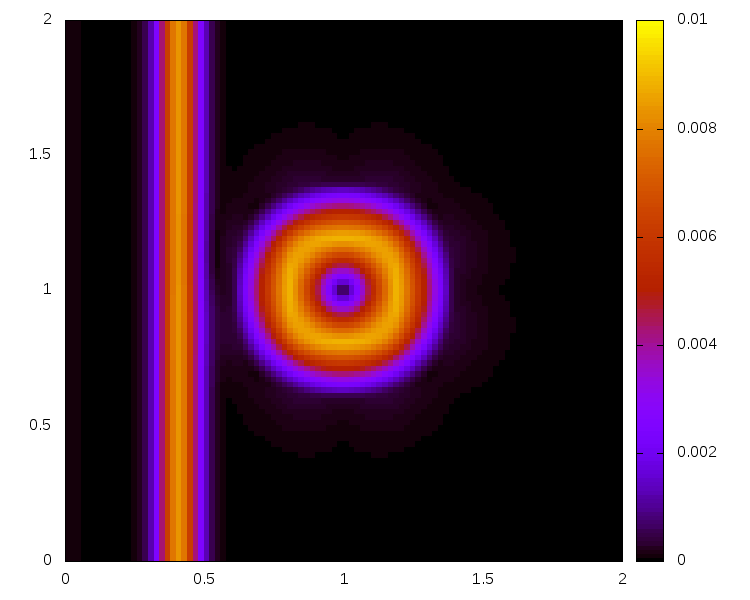}\hfill
\includegraphics[width=0.23\textwidth]{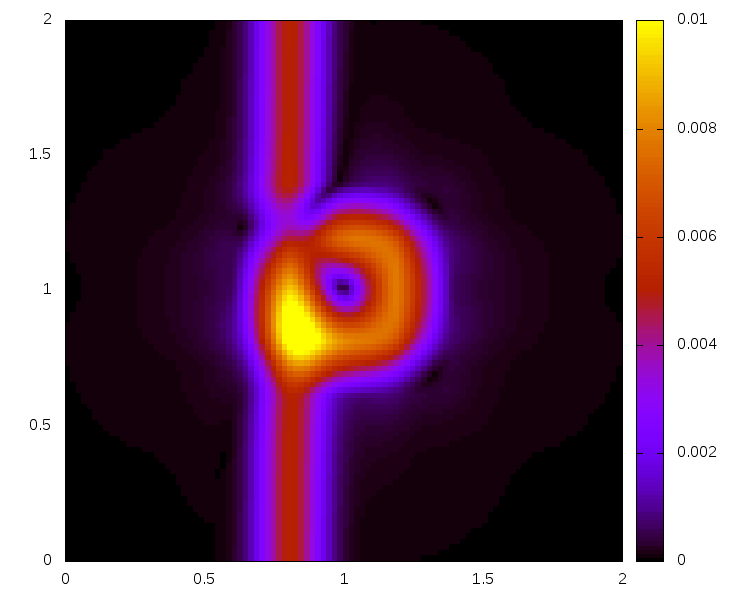}\hfill
\includegraphics[width=0.23\textwidth]{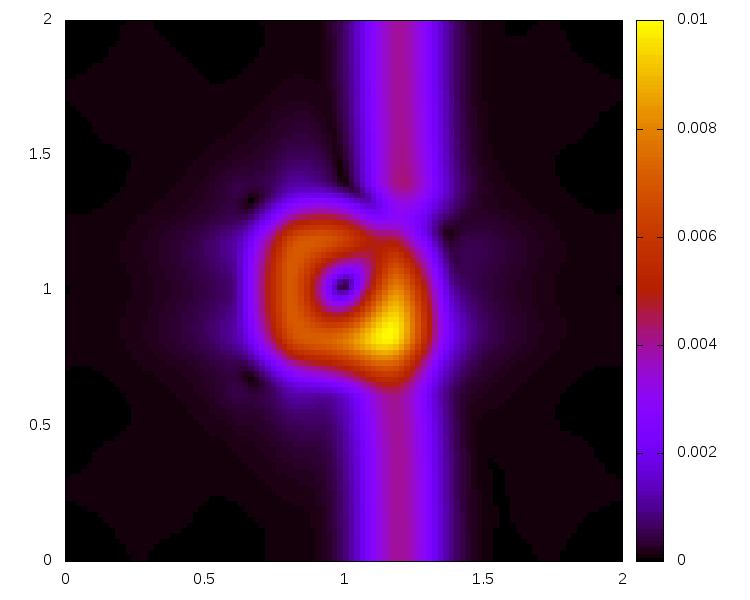}\hfill
\includegraphics[width=0.23\textwidth]{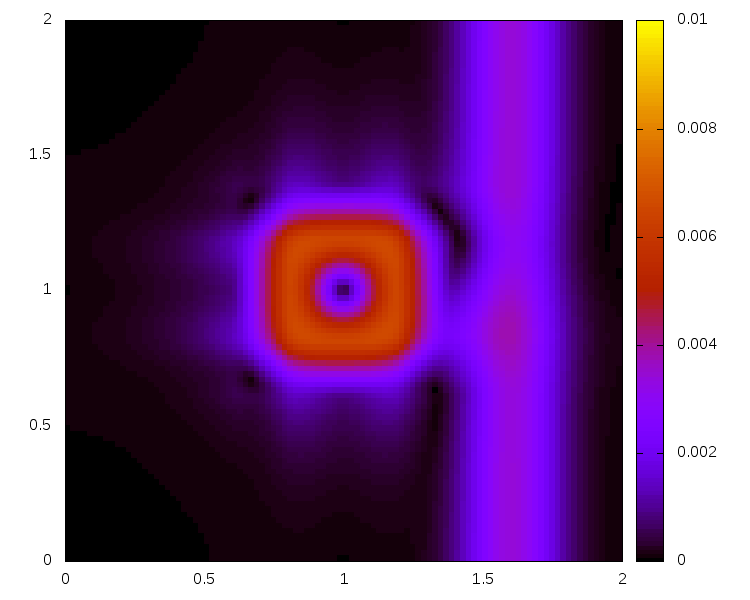}
\caption{{Vortex setup as in Fig. \ref{fig:greshomultid}, but with zero-gradient boundary conditions and additionally a sound wave traveling from left to right. \emph{Top}: Multi-d extension of the Lagrange-Projection method from section \ref{ssec:multidlagpro} using CFL$=0.9$. \emph{Center}: Multi-d extension of a relaxation solver from section \ref{ssec:relaxmultid} using CFL$=0.9$. \emph{Bottom}: Non-low Mach dimensionally split relaxation solver from section \ref{ssec:relax1d} using CFL$=0.45$. Here, $\epsilon = 10^{-2}$ and the domain $[0,2]^2$ is covered by a computational grid of $100\times100$ cells. \emph{Left to right}: $t = 0.002, 0.006, 0.01, 0.014$. Colour coded is the Mach number. The figure demonstrates the ability of the multi-dimensionally extended schemes to simultaneously evolve features of the acoustic and the incompressible regimes.}}
\label{fig:greshoradialconvergencesound}
\end{figure}

In the low Mach number limit, acoustic and incompressible phenomena decouple. To assess the ability of the schemes to cope with both, in Figure \ref{fig:greshoradialconvergencesound}, a right-moving sound wave is superposed on the initial condition of the vortex (in a way similar to \cite{miczek15}). The sound wave is given by
\begin{align}
 p(x) &= 300 \exp\left( -  \left(\frac{x -0.2}{0.02}\right)^2\right) &
 \rho(x) &= \frac{p(x)}{c_\infty^2} &
 u(x) &= \frac{p(x)}{\rho_\infty c_\infty}
\end{align}
where $c_\infty$ and $\rho_\infty$ are the speed of sound and the density outside the (compactly supported) vortex. One observes that the multi-dimensional methods do not have trouble to preserve the vortex even if a sound wave passes through it, while the non-low Mach dimensionally split method dissipates the vortex turning it into a square shape. Recall also that the dimensionally split scheme has a reduced stability range with CFL$< 0.5$. This adds further diffusion which is most visible in the poor resolution of the sound wave.

\subsection{Radial shock}

\begin{figure}
 \centering
 \includegraphics[width=0.48\textwidth]{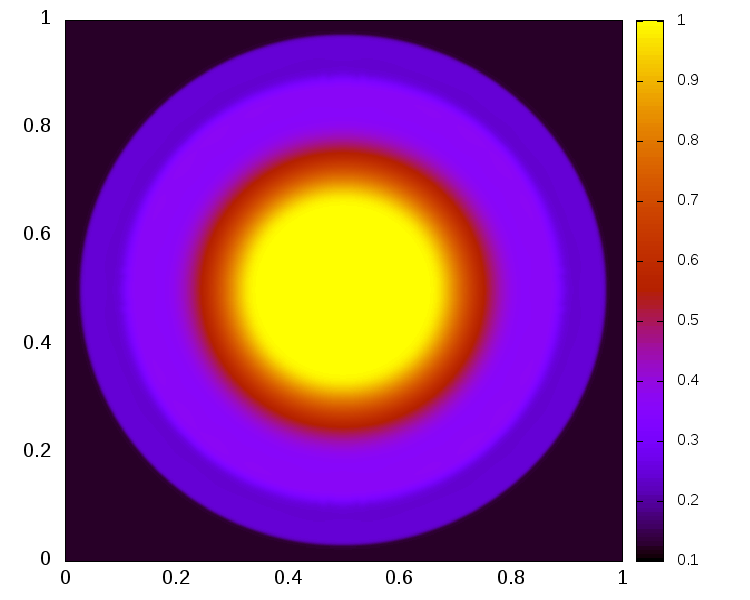} \includegraphics[width=0.48\textwidth]{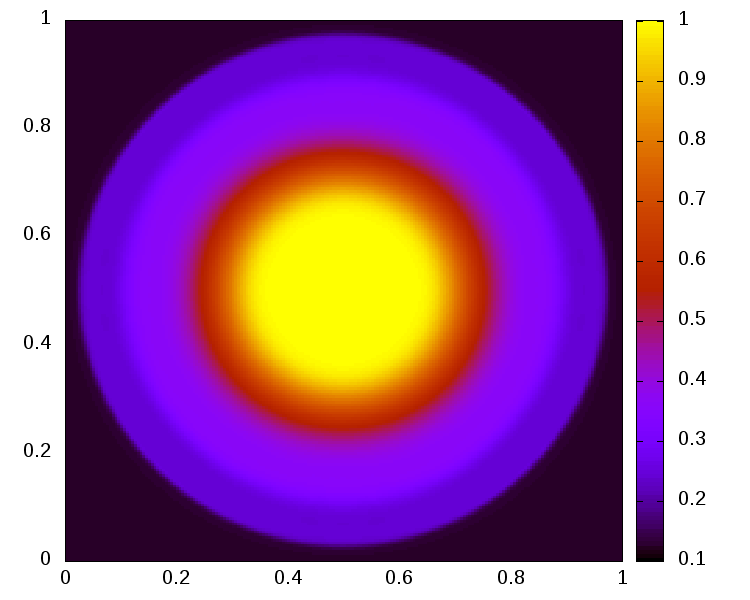} \\
 \includegraphics[width=0.48\textwidth]{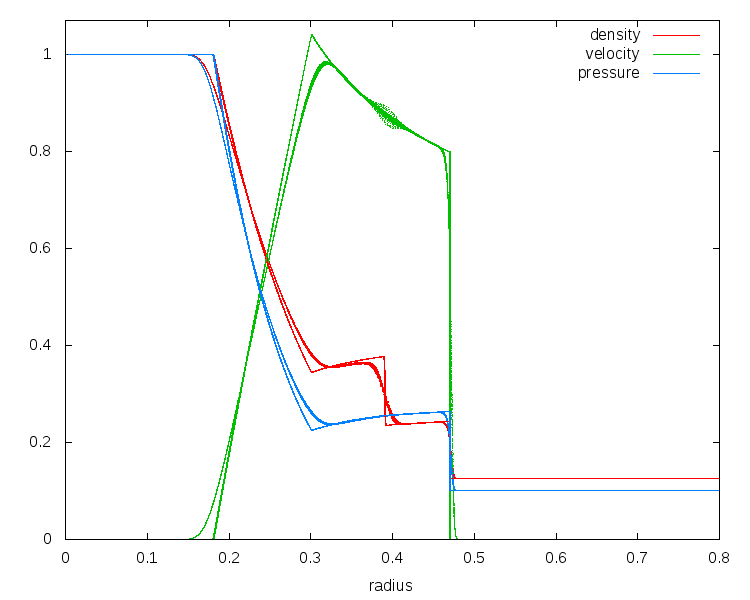} \includegraphics[width=0.48\textwidth]{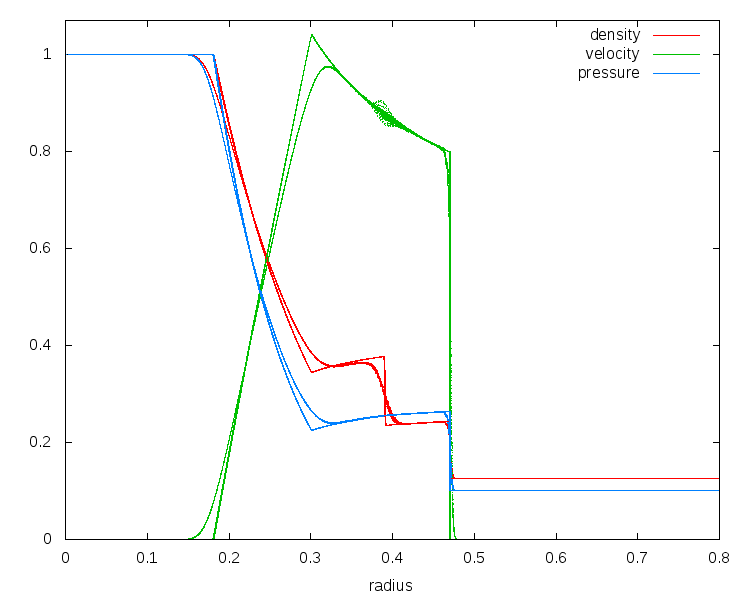} 
 \caption{The radial Sod shock setup solved using the multi-dimensional all-speed extension on a $500\times 500$ grid with CFL $=$ 0.9. Results are shown at $t = 0.1$. {The discontinuity is initially located at a radius of $0.3$.} \emph{Left}: Multi-d modification of a Lagrange-Projection method from section \ref{ssec:multidlagpro}. \emph{Right}: Multi-d modification of a relaxation solver from section \ref{ssec:relaxmultid}. \emph{Top}: Color coded is density. \emph{Bottom}: Radial scatter plots of density, velocity, pressure{, the solid lines show a reference solution obtained by solving just the radial equations on a fine grid}.}
 \label{fig:sodshock}
\end{figure}

Figure \ref{fig:sodshock} demonstrates the all-speed property of the schemes for a spherical version of Sod's shock tube. {Also, the ability of the schemes to maintain the symmetry of the solution is shown by means of a radial plot. The numerical error is most visible in the resolution of the contact discontinuity due to the lack of self-steepening for this wave. Here, the asymmetry of the solution is largest as well, as can be seen by increased scatter. This might be related to the fact that the multi-dimensional modifications suggested in this paper are focusing solely on the behaviour in the low Mach number regime, which mostly does not constrain the choice of the advection operator. Thus, future work (coupling with suggestions such as \cite{leveque97} or an extension to higher order) might improve the results further.}

\subsection{Kelvin-Helmholtz instability}

\begin{figure}
 \centering
 \includegraphics[width=0.24\textwidth]{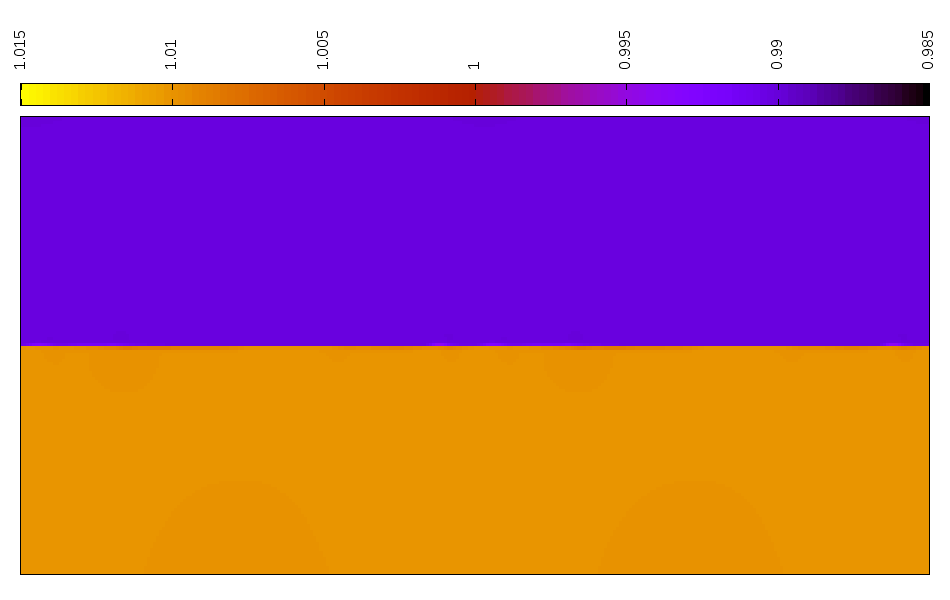} \includegraphics[width=0.24\textwidth]{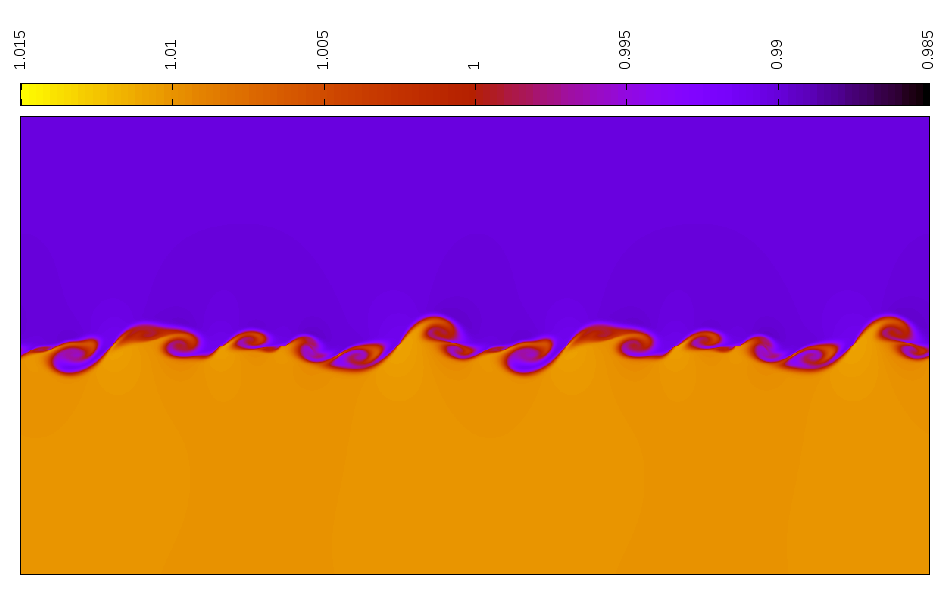} \includegraphics[width=0.24\textwidth]{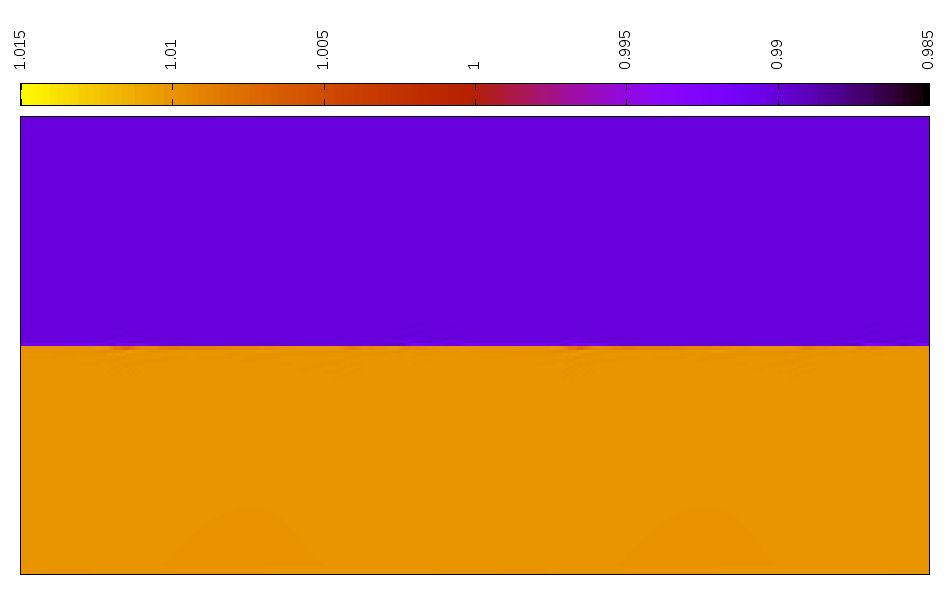} \includegraphics[width=0.24\textwidth]{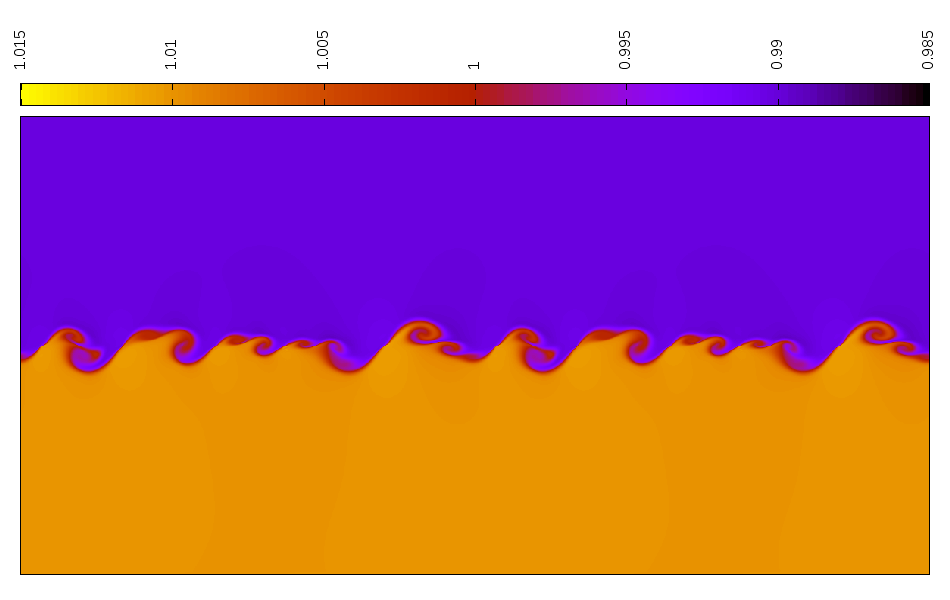}\\
 \includegraphics[width=0.24\textwidth]{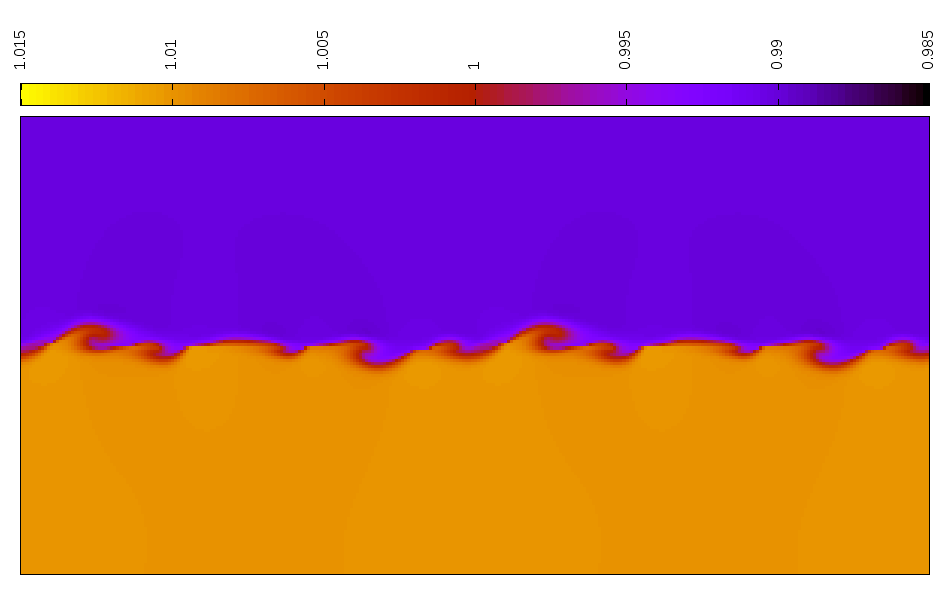} \includegraphics[width=0.24\textwidth]{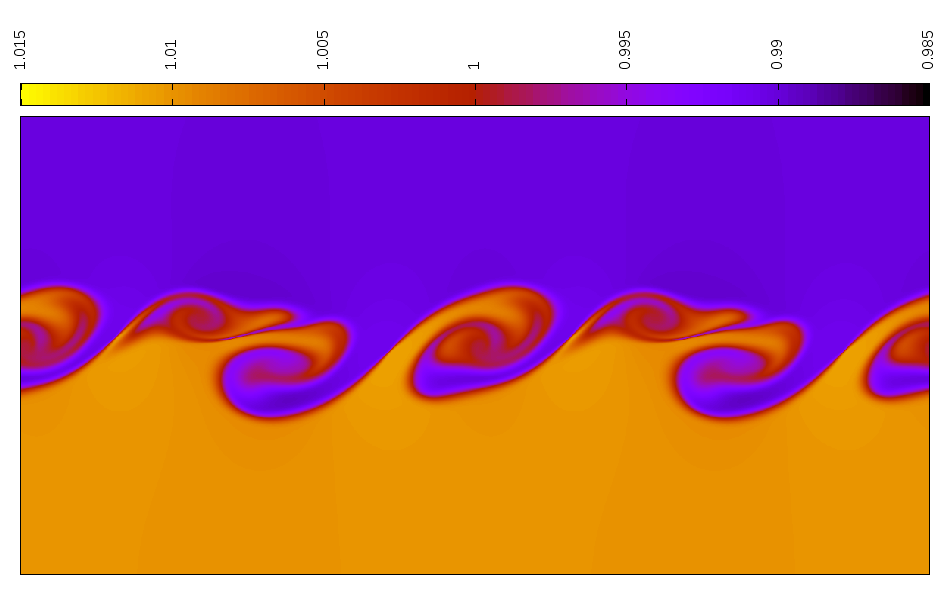} \includegraphics[width=0.24\textwidth]{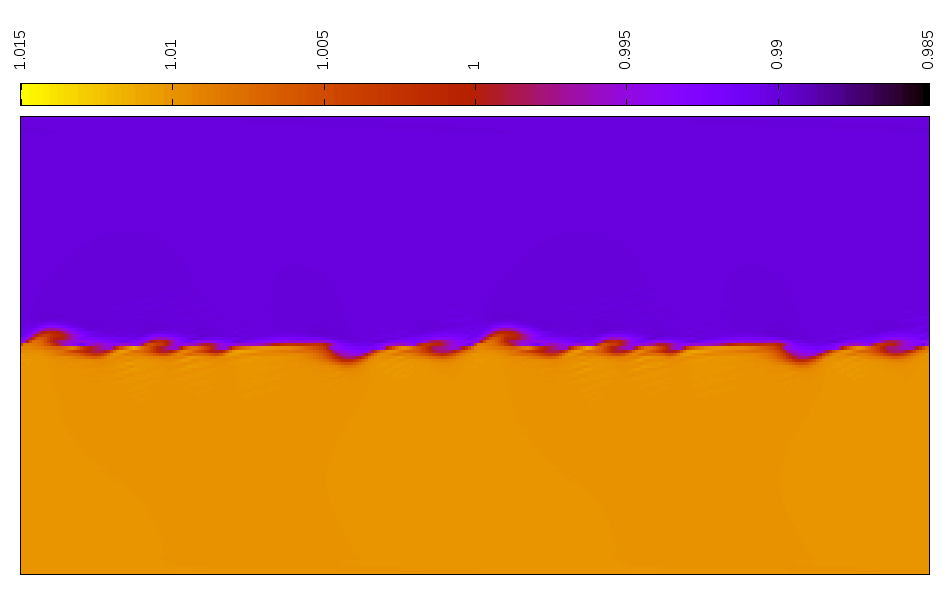} \includegraphics[width=0.24\textwidth]{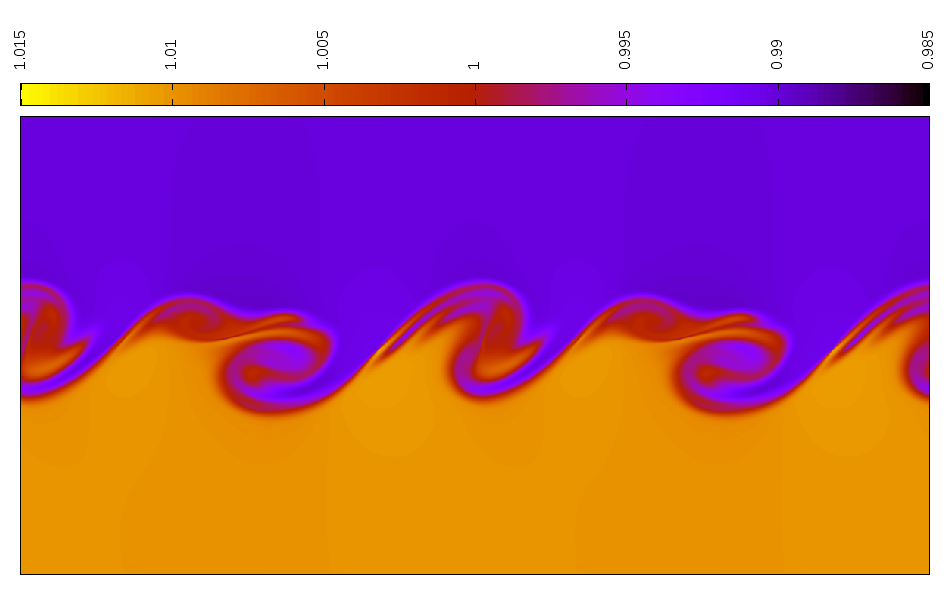}\\
 \includegraphics[width=0.24\textwidth]{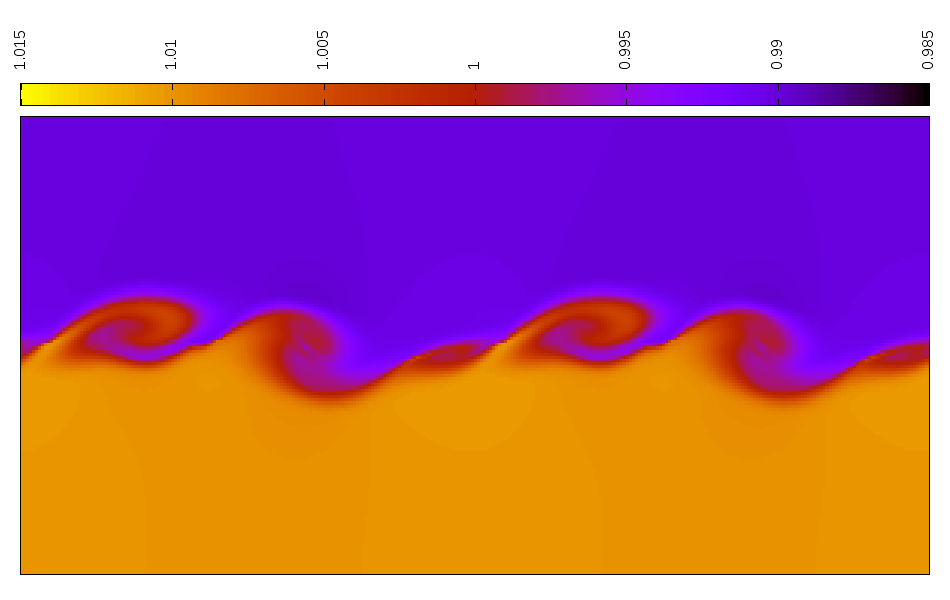} \includegraphics[width=0.24\textwidth]{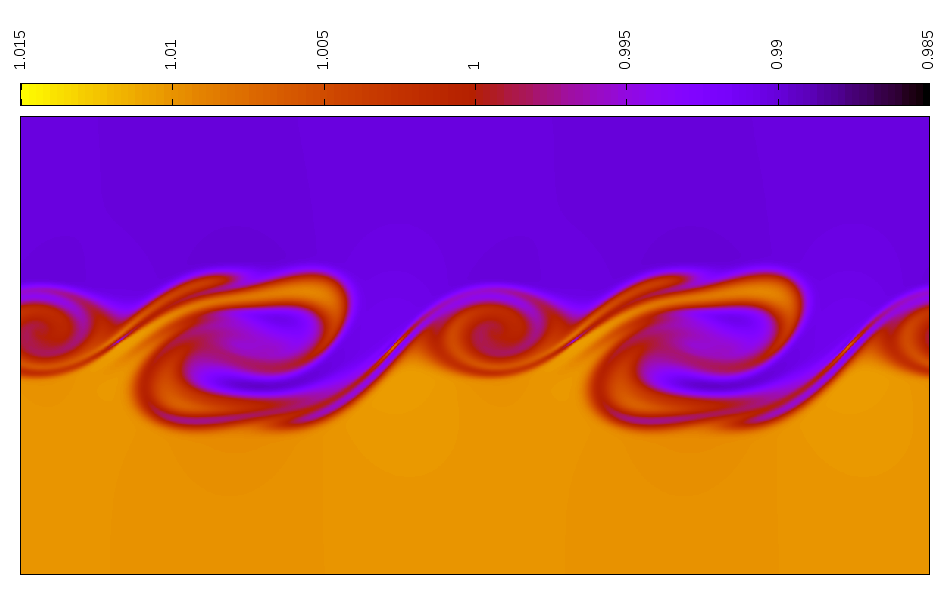} \includegraphics[width=0.24\textwidth]{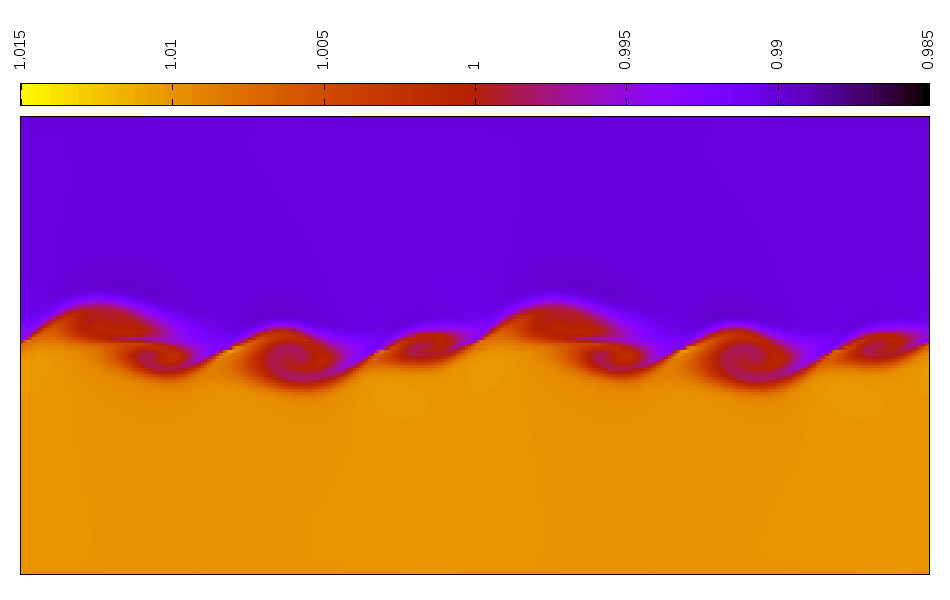} \includegraphics[width=0.24\textwidth]{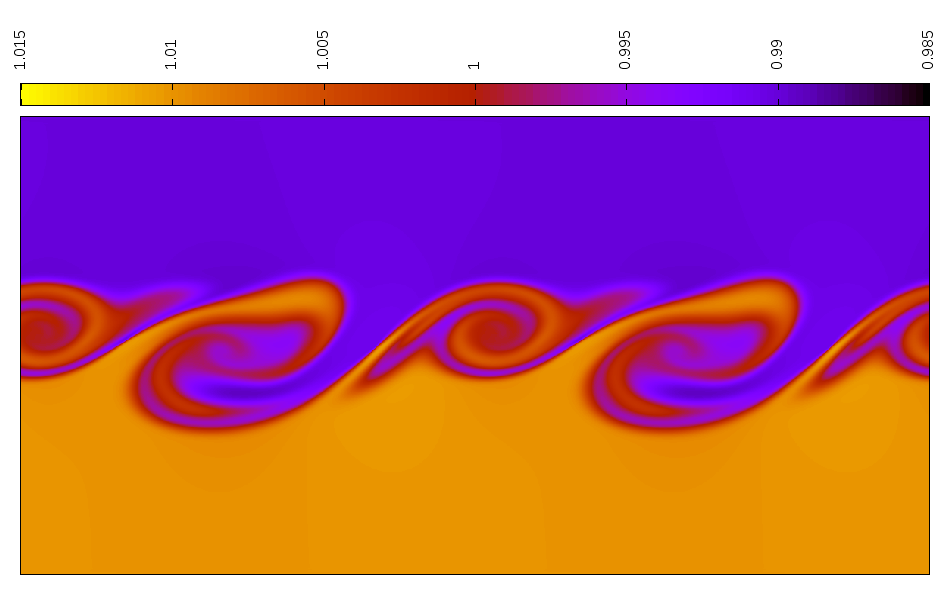}\\
 \includegraphics[width=0.24\textwidth]{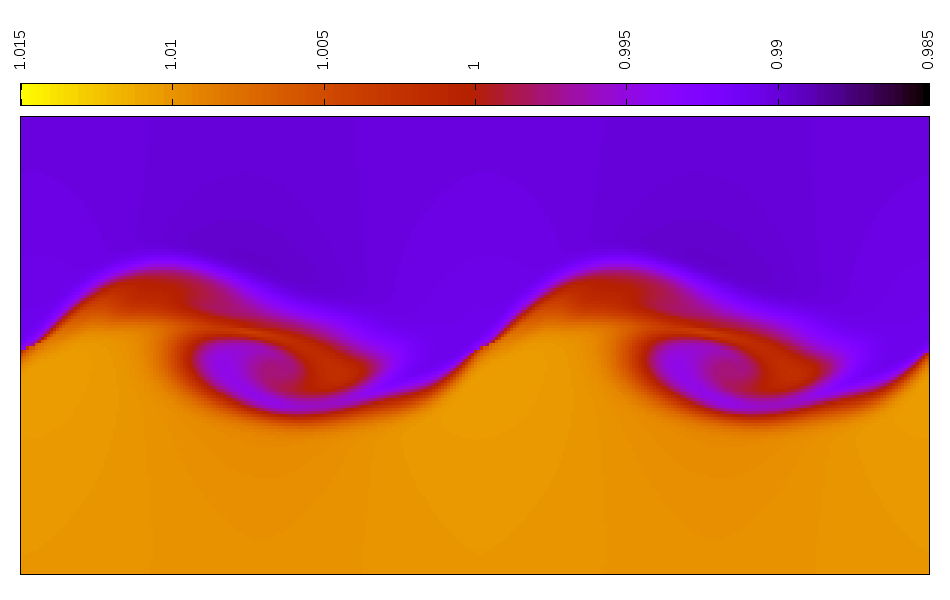} \includegraphics[width=0.24\textwidth]{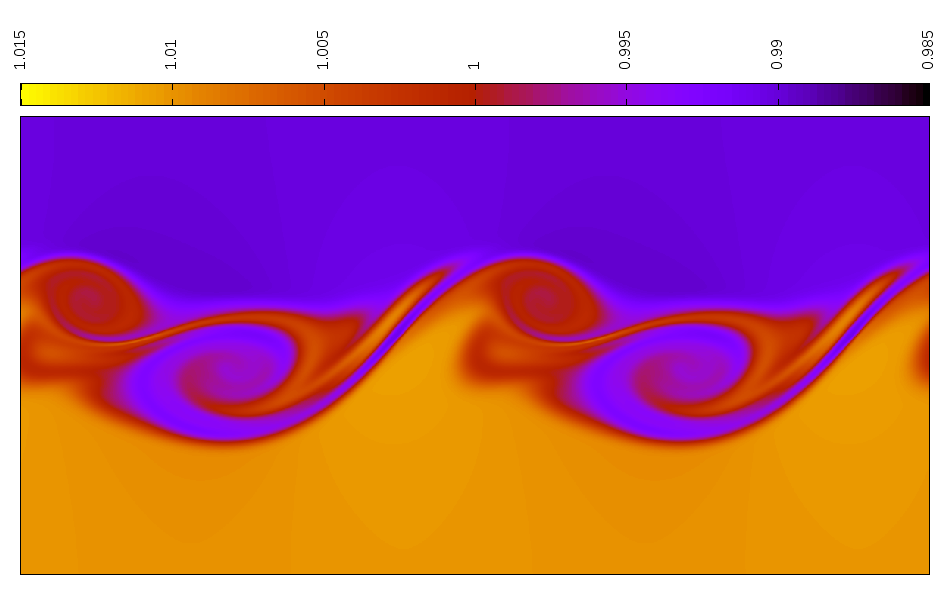} \includegraphics[width=0.24\textwidth]{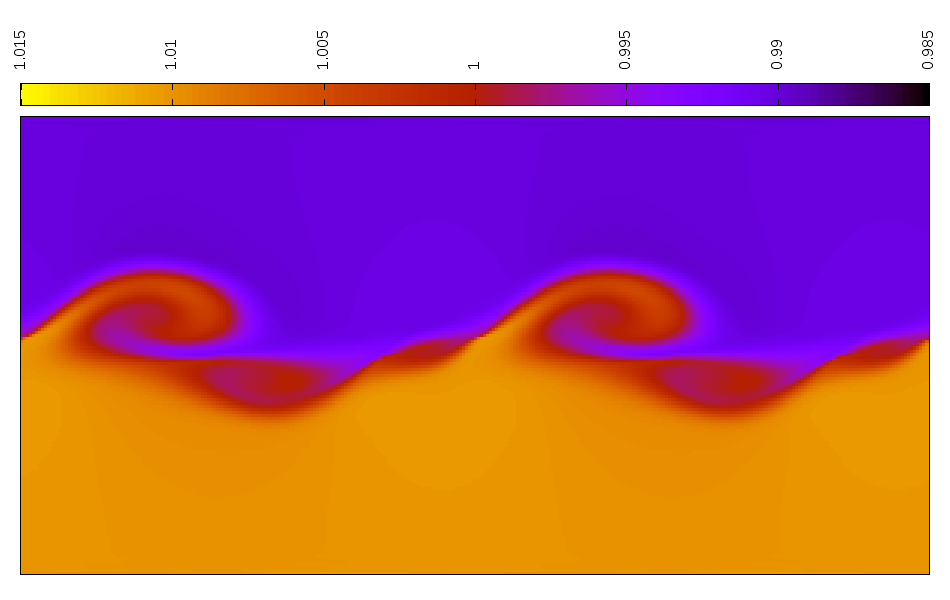} \includegraphics[width=0.24\textwidth]{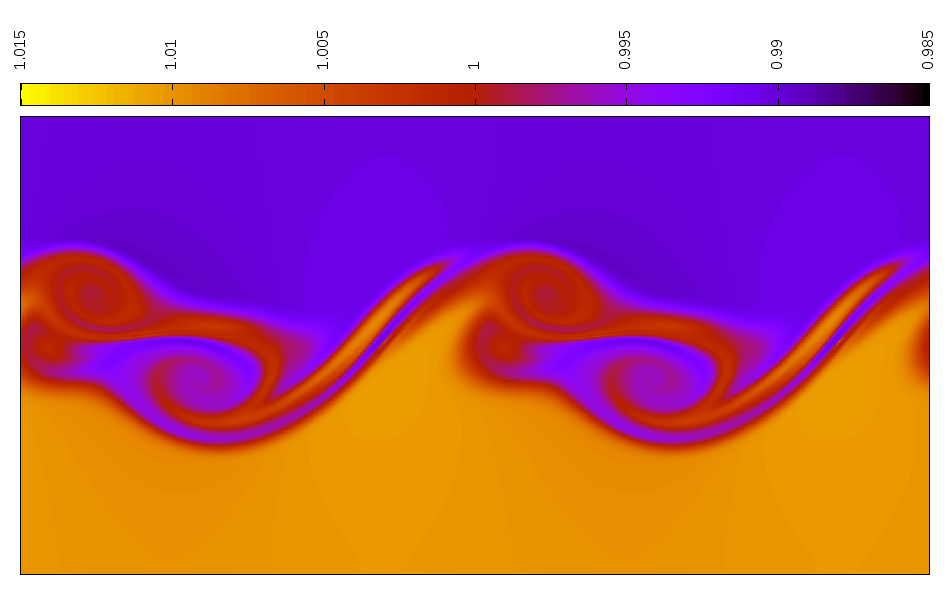}\\
 \includegraphics[width=0.24\textwidth]{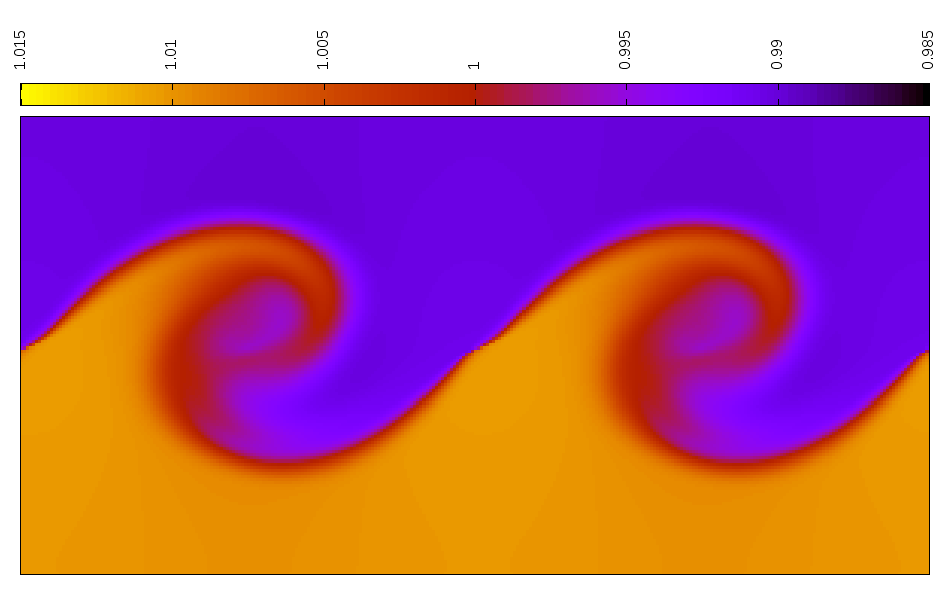} \includegraphics[width=0.24\textwidth]{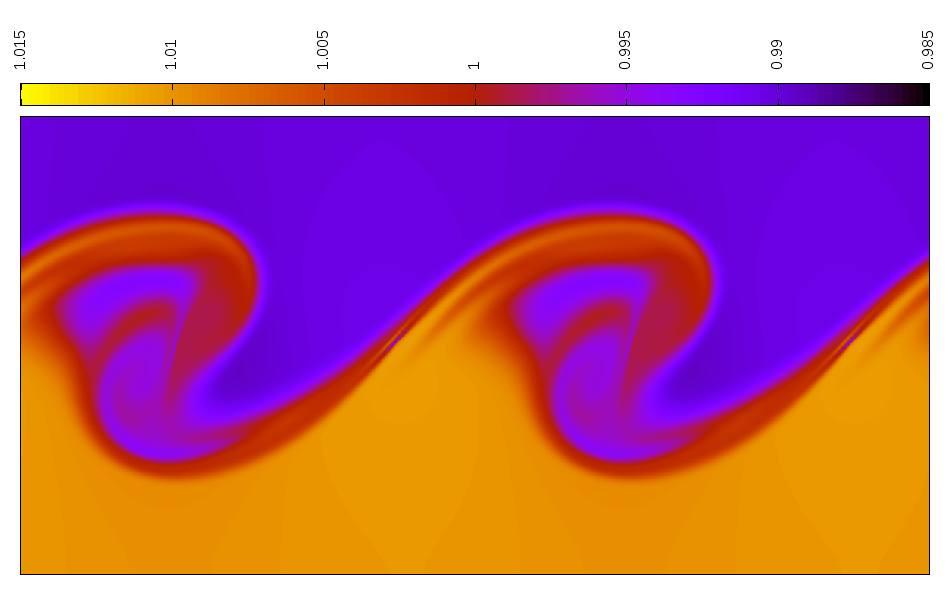} \includegraphics[width=0.24\textwidth]{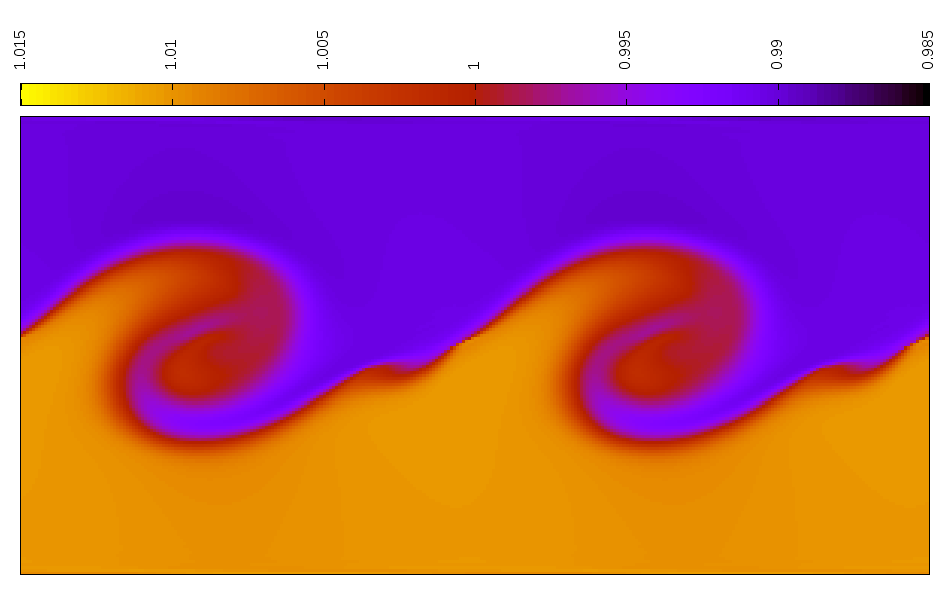} \includegraphics[width=0.24\textwidth]{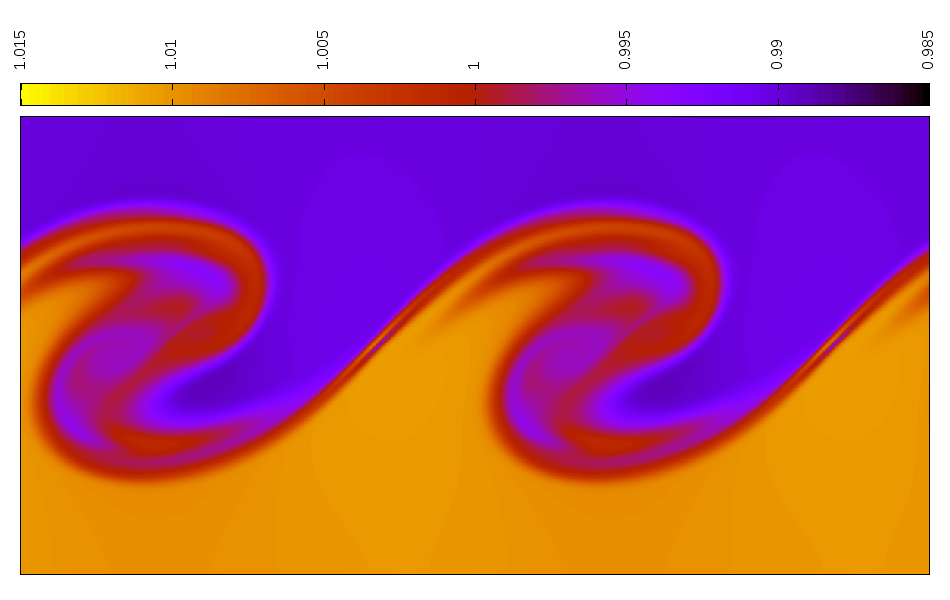}\\
 \includegraphics[width=0.24\textwidth]{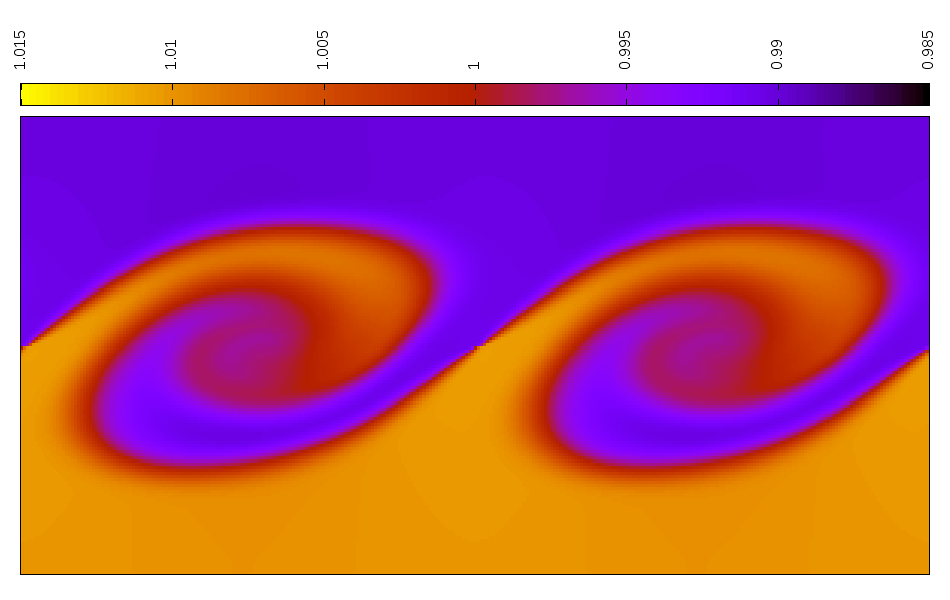} \includegraphics[width=0.24\textwidth]{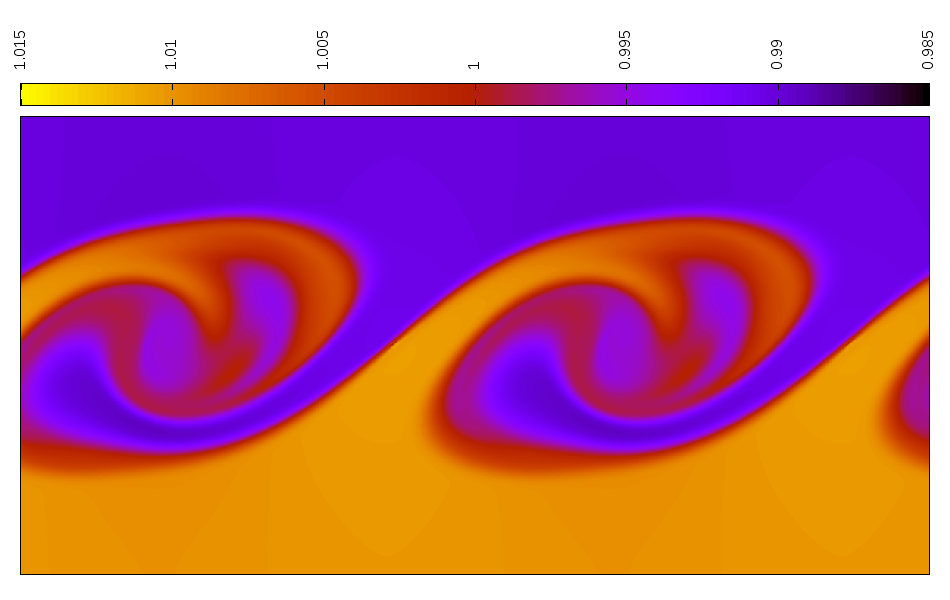} \includegraphics[width=0.24\textwidth]{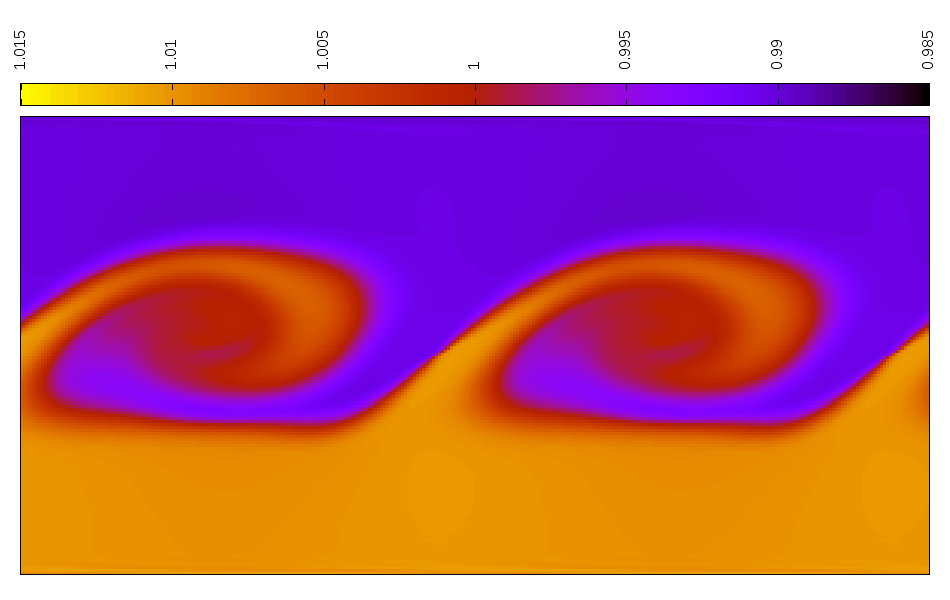} \includegraphics[width=0.24\textwidth]{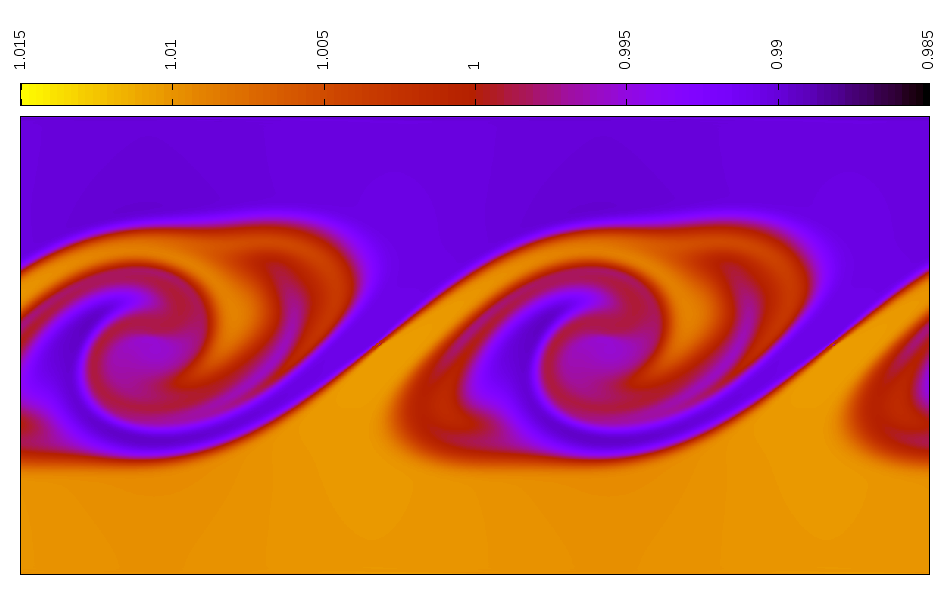}\\
 \includegraphics[width=0.24\textwidth]{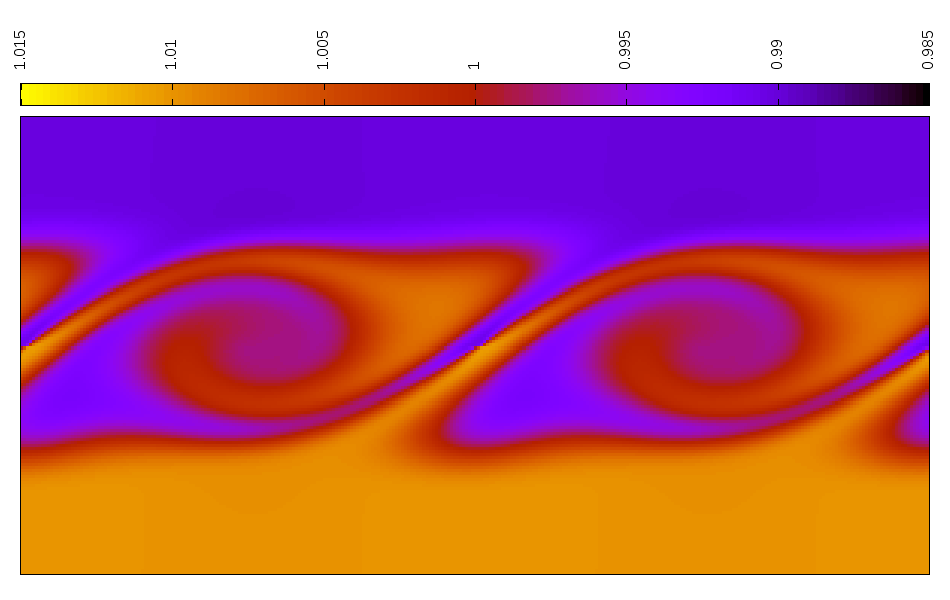} \includegraphics[width=0.24\textwidth]{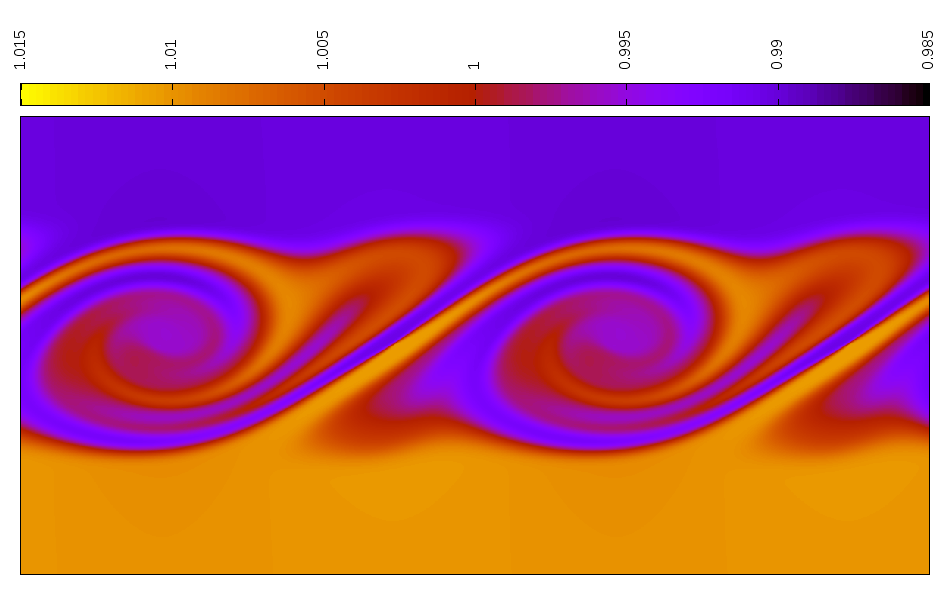} \includegraphics[width=0.24\textwidth]{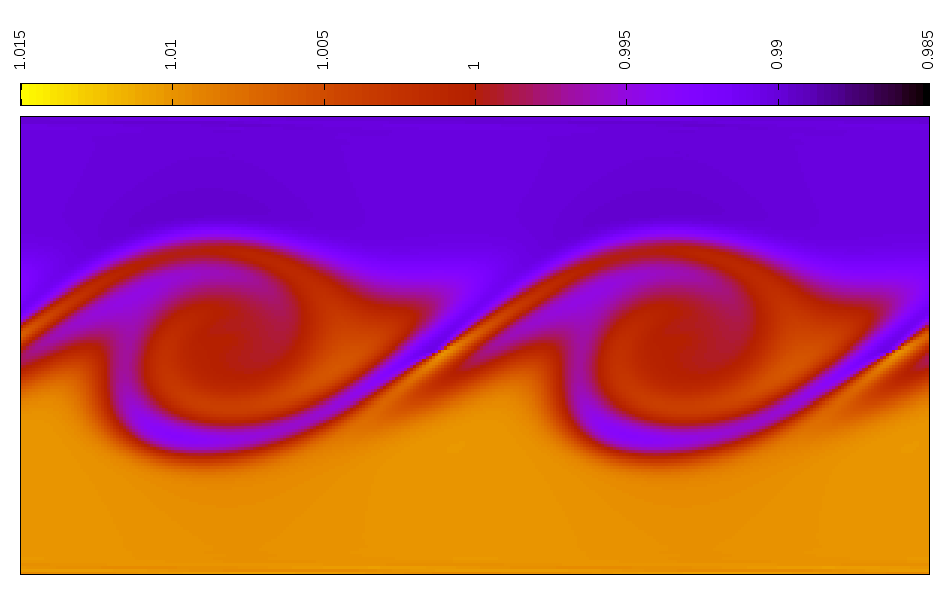} \includegraphics[width=0.24\textwidth]{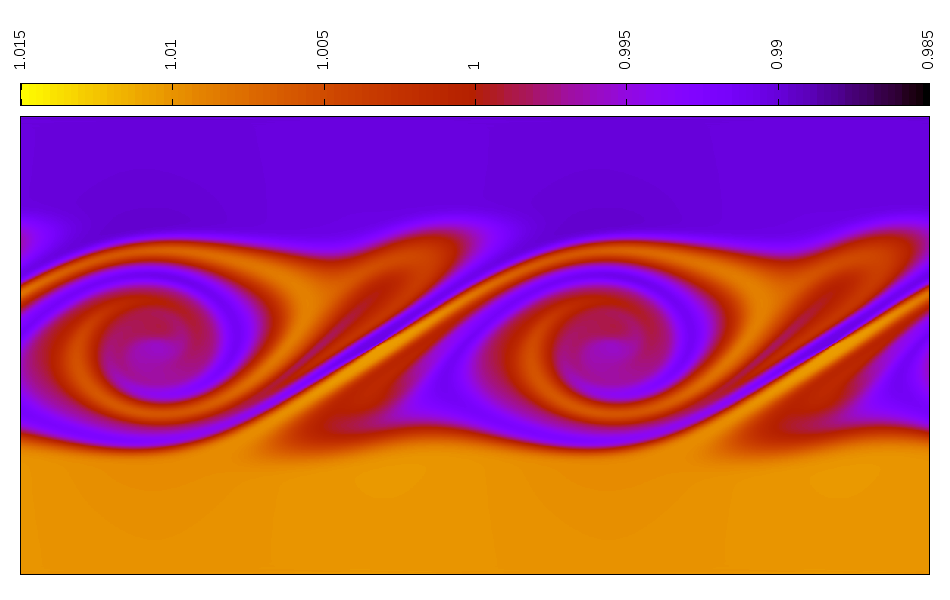}\\
 \includegraphics[width=0.24\textwidth]{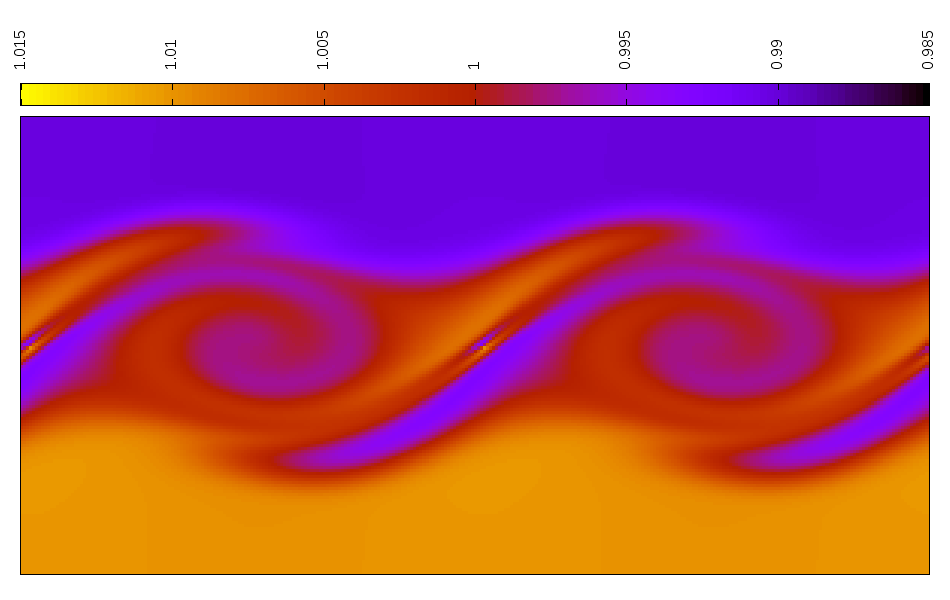} \includegraphics[width=0.24\textwidth]{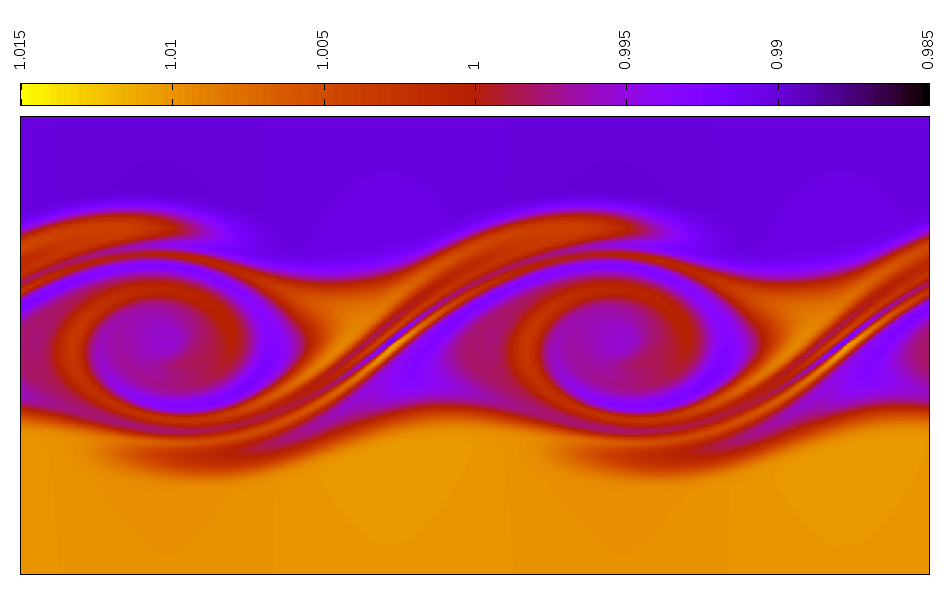} \includegraphics[width=0.24\textwidth]{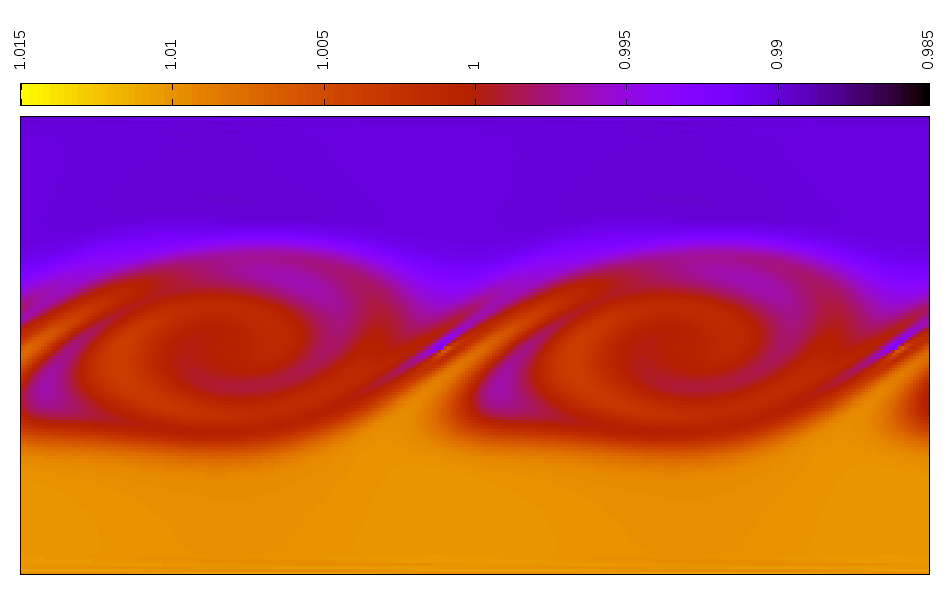} \includegraphics[width=0.24\textwidth]{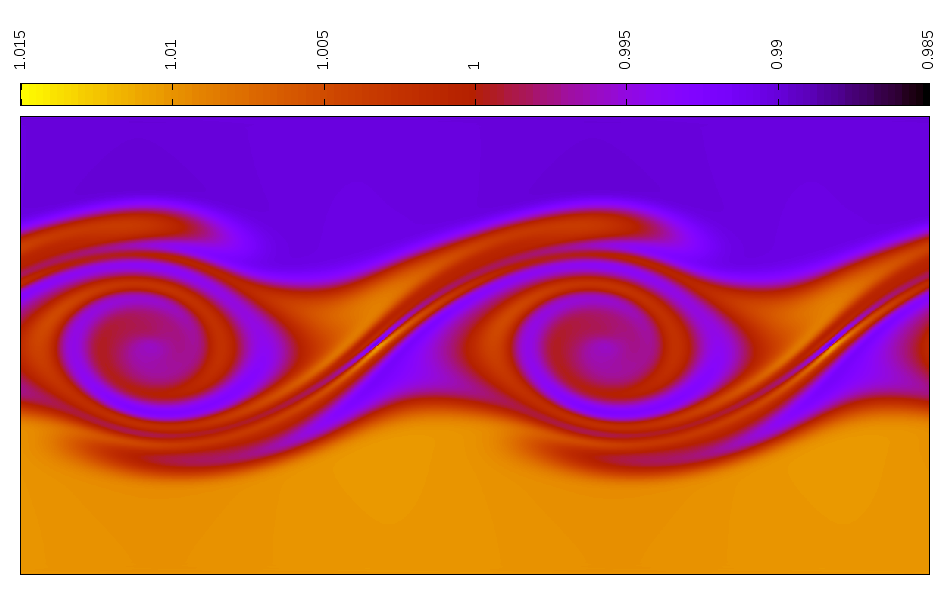}
 \caption{Kelvin-Helmholtz instability solved using the multi-dimensional all-speed extension on a domain $[0,2] \times [0,1]$ using CFL$=0.9$. Color-coded is the density, which is chosen 1\% larger for the left-moving fluid. The interface is initially perturbed with a perpendicular velocity of maximum magnitude of $10^{-3}$ and a wave length of 0.25. \emph{From left to right}: Multi-d modification of the Lagrange-Projection method (section \ref{ssec:multidlagpro}) on a grid of $300\times 150$, same on $1200 \times 600$, multi-d modification of a relaxation solver (section \ref{ssec:relaxmultid}) on $300\times 150$, same on $1200 \times 600$. \emph{From top to bottom}: Times $t = 4, 8, 12, \ldots, 32$.}
 \label{fig:KHmultid}
\end{figure}

\begin{figure}
 \centering
 \includegraphics[width=0.38\textwidth]{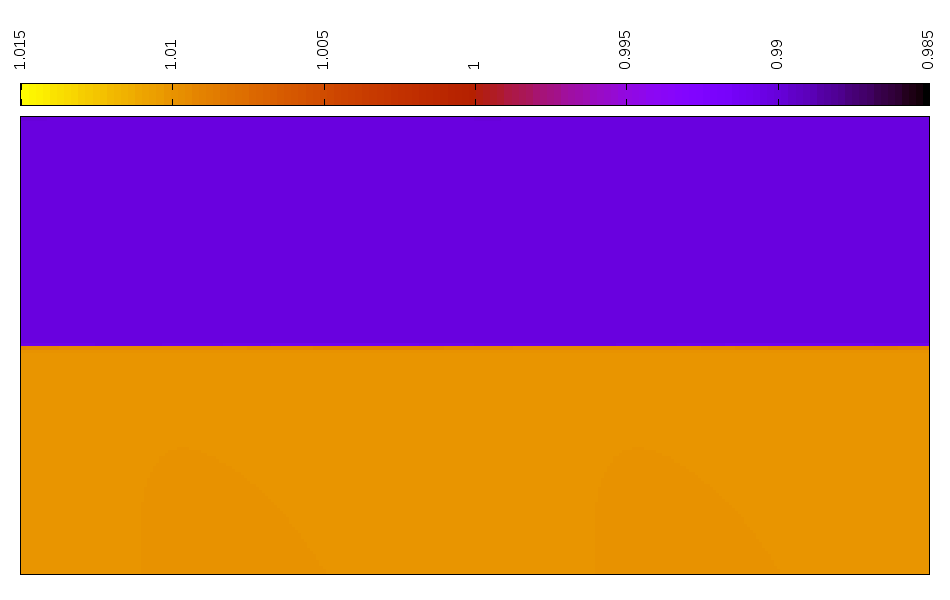} \includegraphics[width=0.38\textwidth]{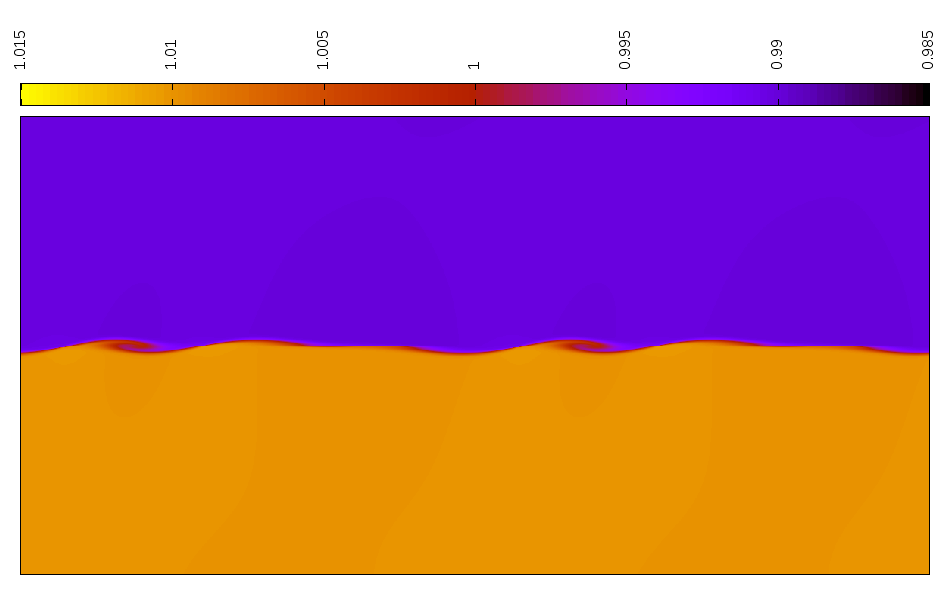} \\
 \includegraphics[width=0.38\textwidth]{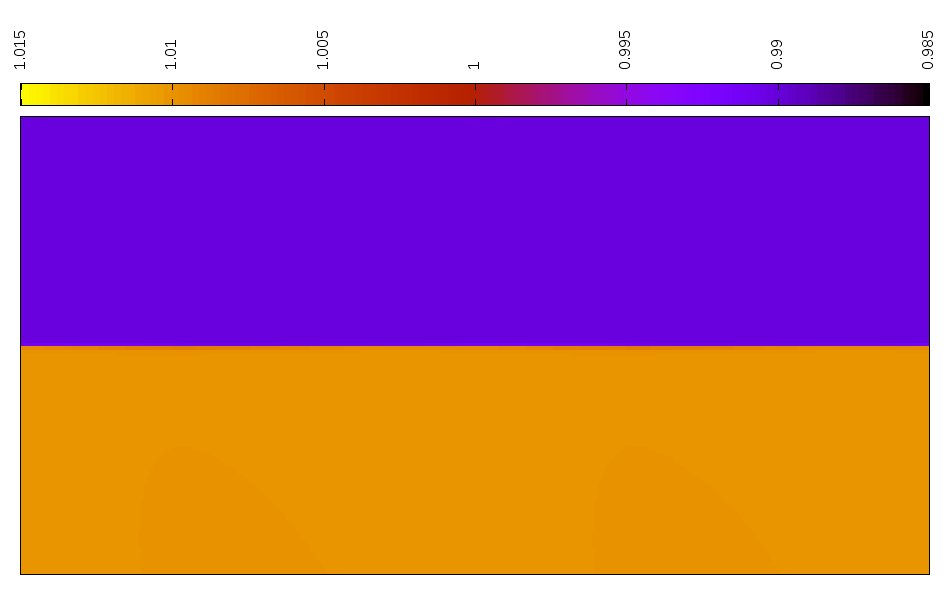} \includegraphics[width=0.38\textwidth]{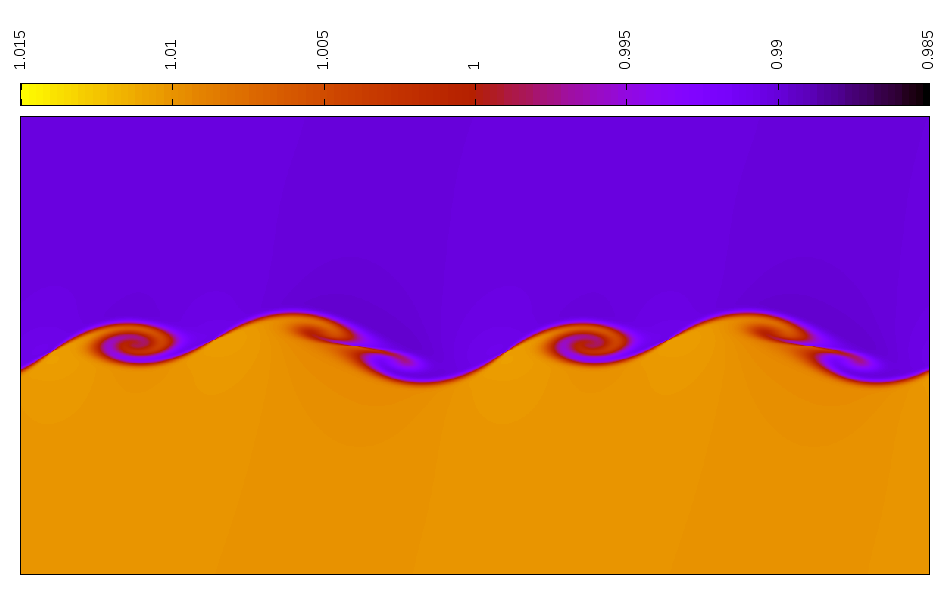} \\
 \includegraphics[width=0.38\textwidth]{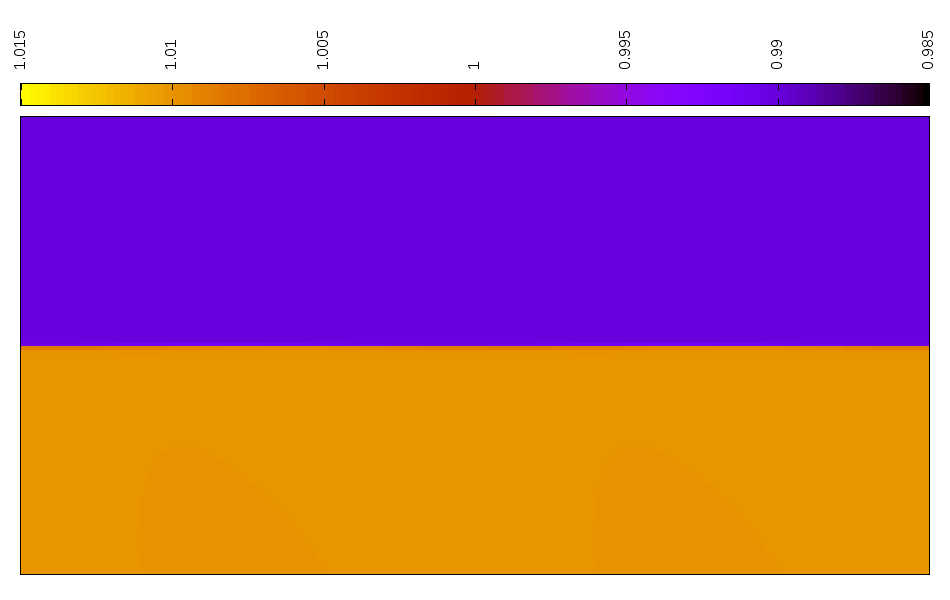} \includegraphics[width=0.38\textwidth]{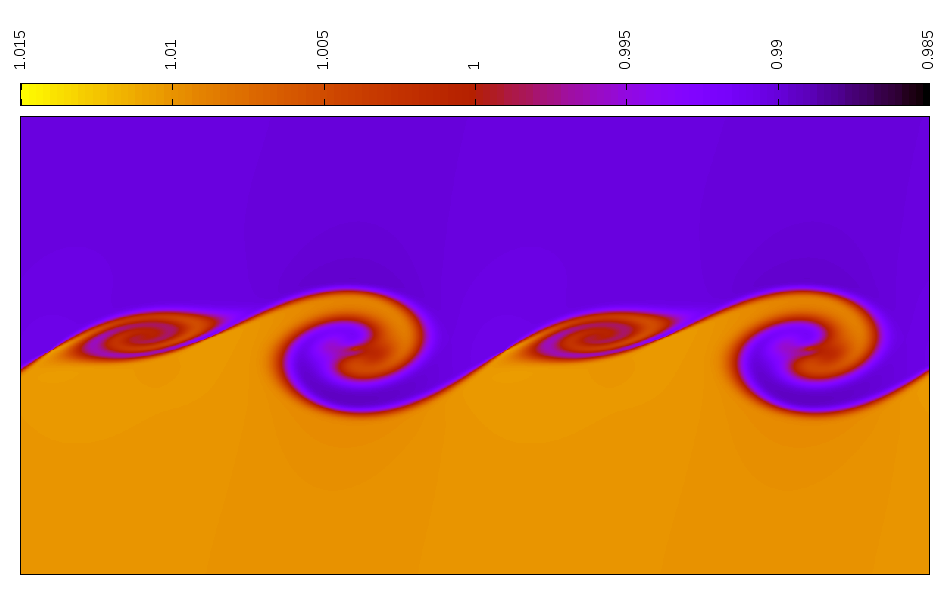} 
 \caption{Same setup as Figure \ref{fig:KHmultid}, solved using the non-low Mach dimensionally split Lagrange Projection method. \emph{Left}: Grid of $300\times 150$. \emph{Right}: Grid of $1200 \times 600$. \emph{From top to bottom}: Times $t = 16, 20, 24$. The interface in the low resolution simulation (left) does not get unstable.}
 \label{fig:KH1d}
\end{figure}

Finally, to assess the properties of the presented numerical methods in a more complex setup, a Kelvin-Helmholtz instability is shown in Figure \ref{fig:KHmultid}. One observes that the interface becomes unstable even for the lowest resolutions considered. There is a transitional phase of vortex leapfrogging and mergers. Generally, the two schemes show similar evolution. These results can be contrasted with the behaviour of the dimensionally split, non-low Mach Lagrange Projection method in Figure \ref{fig:KH1d}. There, the interface is not becoming unstable until the setup is run at high resolution, the simulation results are very different and lack the intricate structures that the other schemes are able to resolve. The artificial stabilization of the interface and lack of small scale flow details can be attributed to an excessive numerical viscosity. The all-speed schemes do not suffer from this problem. Also, they are able to display the intricate vortex dynamics on a low-resolution grid already, which in this case saves a factor of 64 of computational time with respect to the high resolution runs.

This example shows that the quality of a simulation can be dramatically improved when using a numerical method with a well-chosen numerical diffusion: sufficient for stability but not erasing too many details. The all-speed property can be contrasted with an increase of the order of accuracy of a method, which is a different way of reducing the numerical diffusion. Note however that switching to an all-speed scheme does not entail any noteworthy increase in computational cost, or any fundamental change of the algorithmic structure -- be it with usual low Mach fixes, or the multi-dimensional extension presented here. The advantage of the multi-dimensional all-speed extension (apart from the conceptually pleasing property of modifying the scheme only when it becomes necessary, i.e. retaining the one-dimensional scheme) seems to be the stability of the resulting methods. This is very important, as low Mach number modifications deal with precisely this balance of an amount of diffusion necessary to stabilize the method, but otherwise as small as possible. The multi-dimensional all-speed extension also does not include any free parameters or arbitrary functions, which generally are difficult to choose, or might be problem dependent.

\section{Summary and outlook}

Numerical methods for nonlinear systems of conservation laws face a number of challenges in multiple spatial dimensions. They need to resolve many more phenomena than in one dimension at a relatively lower spatial resolution (because for practical reasons it is often not possible to just square or cube the number of cells in the computational grid). One example of such a genuinely multi-dimensional feature is the low Mach number limit.

The strong influence of the numerical diffusion in the low Mach number limit wipes out all the features of the flow. This problem is dealt with in the present paper using a novel approach. 

The usual approach to devise numerical methods for multiple spatial dimensions is to extend a given one-dimensional method via dimensional splitting. This approach is successful in that it mostly yields a stable numerical method, but it is obvious that it neglects subtle balances that might exist between the different directions. The dimensionally split approach was so far also applied when achieving low Mach number compliance of a numerical method. The ne\-ces\-sary modifications were applied to the one-dimensional scheme and the new scheme was used in a dimensionally split manner. If the resulting method was stable, then this strategy has been found to successfully yield a low Mach number compliant method. However, stability can easily be lost, because the low Mach number fix can be understood as a reduction of the stabilizing numerical diffusion. Also, it is conceptually unsatisfying to modify the one-dimensional scheme, if a nontrivial low Mach number limit and the difficulties it entails only appear in multiple dimensions and there is nothing ``wrong'' with the one-dimensional scheme in first place. Finally, the low Mach number fixes always involve an arbitrary function, which is not necessarily easy to specify generally: using the local Mach number, for example, is not sufficient for an initially static shock tube setups.

The alternative strategy to achieve low Mach number compliance presented in this paper does not suffer from these issues. The one-dimensional scheme is left untouched and it is only ensured that its multi-dimensional extension is done in a particular, low Mach number compliant way. This yields numerical methods that do not suffer from a loss of stability and do not involve any free parameters or functions. The low Mach number fix also does not need to be ``turned off'' by hand once the local Mach number approaches unity: the suggested strategy yields an all-speed scheme automatically.

The question to which this paper is also contributing is whether a multi-dimensional phenomenon necessitates a truly multi-dimensional numerical me\-thod. The proposed view is that the usage of a truly multi-dimensional method opens up many possibilities while it does not require virtually any extra computational time and only little bookkeeping effort.

The new strategy is exemplified in this paper on two schemes: a Lagrange Projection method and a relaxation solver. Both schemes are analyzed and extended to multiple dimensions in an all-speed fashion. For the first of them it has already been shown that a simple low Mach number fix achieves low Mach number compliance, so this example serves to demonstrate that the new strategy is at least as good. For the second example it has been previously found that the usual low Mach number fix (which modifies the one-dimensional scheme already) leads to instability. The multi-dimensional strategy of achieving low Mach number compliance is demonstrated to lead to a stable all-speed scheme also in this example. This scheme is simpler than the Lagrange Projection method and therefore might be more attractive for practical implementation.

The multi-dimensional schemes, and in particular the truly multi-dimensional divergence discretizations first arose in the context of involution preserving schemes for linear acoustics. In future, it might be interesting to investigate whether the enlarged flexibility of the multi-dimensional all-speed extension as described here leads to further structure preservation properties of the resulting scheme.

Achieving structure preservation is independent from the extension of a method to higher order. Although the latter generally comes at an increased computational or memory cost, clearly, certain applications might profit from combining these two approaches of diffusion reduction. 

Another aspect of the difficulties that numerical methods face in the low Mach number limit is the decrease of the time step for explicit time integration, which can be dealt with via implicit or semi-implicit time integration. Future work shall be devoted to finding combinations of these approaches with the multi-dimensional all-speed extension as well as an analysis of the conditions necessary to show a discrete entropy instability.

\section*{Acknowledgements}
 The author was supported by the German Academic Exchange Service (DAAD) with funds from the German Federal Ministry of Education and Research (BMBF) and the European Union (FP7-PEOPLE-2013-COFUND -- grant agreement no. 605728) as well as by the Deutsche Forschungsgemeinschaft (DFG) through project 429491391 (BA 6878/1-1).

\bibliographystyle{elsarticle-num}

\begin{thebibliography}{10}
\expandafter\ifx\csname url\endcsname\relax
  \def\url#1{\texttt{#1}}\fi
\expandafter\ifx\csname urlprefix\endcsname\relax\def\urlprefix{URL }\fi
\expandafter\ifx\csname href\endcsname\relax
  \def\href#1#2{#2} \def\path#1{#1}\fi

\bibitem{ebin77}
D.~G. Ebin, The motion of slightly compressible fluids viewed as a motion with
  strong constraining force, Annals of mathematics (1977) 141--200.

\bibitem{klainerman81}
S.~Klainerman, A.~Majda, Singular limits of quasilinear hyperbolic systems with
  large parameters and the incompressible limit of compressible fluids,
  Communications on Pure and Applied Mathematics 34~(4) (1981) 481--524.

\bibitem{metivier01}
G.~M{\'e}tivier, S.~Schochet, The incompressible limit of the non-isentropic
  {E}uler equations, Archive for rational mechanics and analysis 158~(1) (2001)
  61--90.

\bibitem{courant28}
R.~Courant, K.~Friedrichs, H.~Lewy, {\"U}ber die partiellen
  differenzengleichungen der mathematischen {P}hysik, Mathematische Annalen
  100~(1) (1928) 32--74.

\bibitem{vanleer06}
B.~Van~Leer, Upwind and high-resolution methods for compressible flow: From
  donor cell to residual-distribution schemes, Communications in Computational
  Physics 1~(192-206) (2006) 138.

\bibitem{lax64}
P.~D. Lax, B.~Wendroff, Difference schemes for hyperbolic equations with high
  order of acuracy, Comm. Pure Appl. Math. 17 (1964) 381--398.

\bibitem{tadmor2003entropy}
E.~Tadmor, Entropy stability theory for difference approximations of nonlinear
  conservation laws and related time-dependent problems, Acta Numerica 12
  (2003) 451--512.

\bibitem{vanleer79}
B.~Van~Leer, Towards the ultimate conservative difference scheme. v. {A}
  second-order sequel to {G}odunov's method, Journal of computational Physics
  32~(1) (1979) 101--136.

\bibitem{colella84}
P.~Colella, P.~R. Woodward, The piecewise parabolic method ({PPM}) for
  gas-dynamical simulations, Journal of computational physics 54~(1) (1984)
  174--201.

\bibitem{jiang96}
G.-S. Jiang, C.-W. Shu, Efficient implementation of weighted {ENO} schemes,
  Journal of computational physics 126~(1) (1996) 202--228.

\bibitem{cockburn98}
B.~Cockburn, C.-W. Shu, The {R}unge-{K}utta discontinuous {G}alerkin method for
  conservation laws {V}: {M}ultidimensional systems, Journal of Computational
  Physics 141~(2) (1998) 199--224.

\bibitem{abgrall93}
R.~Abgrall, \href{https://hal.inria.fr/inria-00074814}{A genuinely
  multidimensional {R}iemann solver}, hal:inria-00074814.
\newline\urlprefix\url{https://hal.inria.fr/inria-00074814}

\bibitem{leveque97}
R.~J. LeVeque, Wave propagation algorithms for multidimensional hyperbolic
  systems, Journal of Computational Physics 131~(2) (1997) 327--353.

\bibitem{gilquin96}
H.~Gilquin, J.~Laurens, C.~Rosier, Multi-dimensional {R}iemann problems for
  linear hyperbolic systems, ESAIM: Mathematical Modelling and Numerical
  Analysis 30~(5) (1996) 527--548.

\bibitem{barsukow17}
W.~Barsukow, C.~Klingenberg, Exact solution and a truly multidimensional
  {G}odunov scheme for the acoustic equations, submitted, preprint available as
  arXiv:2004.04217.

\bibitem{colella90}
P.~Colella, Multidimensional upwind methods for hyperbolic conservation laws,
  Journal of Computational Physics 87~(1) (1990) 171--200.

\bibitem{zheng12}
Y.~Zheng, Systems of Conservation Laws: Two-Dimensional Riemann Problems,
  Springer Science \& Business Media, 2001.

\bibitem{morton01}
K.~W. Morton, P.~L. Roe, Vorticity-preserving {L}ax-{W}endroff-type schemes for
  the system wave equation, SIAM Journal on Scientific Computing 23~(1) (2001)
  170--192.

\bibitem{jeltsch06}
R.~Jeltsch, M.~Torrilhon, On curl-preserving finite volume discretizations for
  shallow water equations, BIT Numerical Mathematics 46~(1) (2006) 35--53.

\bibitem{mishra09preprint}
S.~Mishra, E.~Tadmor, Constraint preserving schemes using potential-based
  fluxes {II}. genuinely multi-dimensional central schemes for systems of
  conservation laws, ETH preprint~(2009-32).

\bibitem{barsukow17a}
W.~Barsukow, Stationarity preserving schemes for multi-dimensional linear
  systems, Mathematics of Computation 88~(318) (2019) 1621--1645.

\bibitem{guillard04}
H.~Guillard, A.~Murrone, On the behavior of upwind schemes in the low {M}ach
  number limit: {II}. {G}odunov type schemes, Computers \& fluids 33~(4) (2004)
  655--675.

\bibitem{rieper09}
F.~Rieper, G.~Bader, The influence of cell geometry on the accuracy of upwind
  schemes in the low {M}ach number regime, Journal of Computational Physics
  228~(8) (2009) 2918--2933.

\bibitem{dellacherierieper10}
S.~Dellacherie, P.~Omnes, F.~Rieper, The influence of cell geometry on the
  {G}odunov scheme applied to the linear wave equation, Journal of
  Computational Physics 229~(14) (2010) 5315--5338.

\bibitem{guillard09}
H.~Guillard, On the behavior of upwind schemes in the low {M}ach number limit.
  {IV}: {P}0 approximation on triangular and tetrahedral cells, Computers \&
  Fluids 38~(10) (2009) 1969--1972.

\bibitem{rieper10}
F.~Rieper, On the dissipation mechanism of upwind-schemes in the low {M}ach
  number regime: {A} comparison between {R}oe and {HLL}, Journal of
  Computational Physics 229~(2) (2010) 221--232.

\bibitem{haack12}
J.~Haack, S.~Jin, J.-G. Liu, An all-speed asymptotic-preserving method for the
  isentropic {E}uler and {N}avier-{S}tokes equations, preprint.

\bibitem{cordier12}
F.~Cordier, P.~Degond, A.~Kumbaro, An asymptotic-preserving all-speed scheme
  for the {E}uler and {N}avier-{S}tokes equations, Journal of Computational
  Physics 231~(17) (2012) 5685--5704.

\bibitem{li08}
X.-s. Li, C.-w. Gu, An all-speed {R}oe-type scheme and its asymptotic analysis
  of low {M}ach number behaviour, Journal of Computational Physics 227~(10)
  (2008) 5144--5159.

\bibitem{thornber08}
B.~Thornber, D.~Drikakis, Numerical dissipation of upwind schemes in low {M}ach
  flow, International journal for numerical methods in fluids 56~(8) (2008)
  1535--1541.

\bibitem{dellacherie10}
S.~Dellacherie, Analysis of {G}odunov type schemes applied to the compressible
  {E}uler system at low {M}ach number, Journal of Computational Physics 229~(4)
  (2010) 978--1016.

\bibitem{rieper11}
F.~Rieper, A low-{M}ach number fix for {R}oe’s approximate {R}iemann solver,
  Journal of Computational Physics 230~(13) (2011) 5263--5287.

\bibitem{li13}
X.-s. Li, C.-w. Gu, Mechanism of {R}oe-type schemes for all-speed flows and its
  application, Computers \& Fluids 86 (2013) 56--70.

\bibitem{chalons16}
C.~Chalons, M.~Girardin, S.~Kokh, An all-regime lagrange-projection like scheme
  for the gas dynamics equations on unstructured meshes, Communications in
  Computational Physics 20~(1) (2016) 188--233.

\bibitem{birken16}
K.~O{\ss}wald, A.~Siegmund, P.~Birken, V.~Hannemann, A.~Meister, L2roe: a low
  dissipation version of {R}oe's approximate {R}iemann solver for low {M}ach
  numbers, International Journal for Numerical Methods in Fluids 81~(2) (2016)
  71--86.

\bibitem{barsukow16}
W.~Barsukow, P.~V. Edelmann, C.~Klingenberg, F.~Miczek, F.~K. R{\"o}pke, A
  numerical scheme for the compressible low-{M}ach number regime of ideal fluid
  dynamics, Journal of Scientific Computing 72~(2) (2017) 623--646.

\bibitem{dellacherie16}
S.~Dellacherie, J.~Jung, P.~Omnes, P.-A. Raviart, Construction of modified
  {G}odunov-type schemes accurate at any {M}ach number for the compressible
  {E}uler system, Mathematical Models and Methods in Applied Sciences 26~(13)
  (2016) 2525--2615.

\bibitem{birken05}
P.~Birken, A.~Meister, Stability of preconditioned finite volume schemes at low
  {M}ach numbers, BIT Numerical Mathematics 45~(3) (2005) 463--480.

\bibitem{klein95}
R.~Klein, Semi-implicit extension of a {G}odunov-type scheme based on low
  {M}ach number asymptotics {I}: {O}ne-dimensional flow, Journal of
  Computational Physics 121~(2) (1995) 213--237.

\bibitem{lung14}
T.~Lung, P.~Roe, Toward a reduction of mesh imprinting, International Journal
  for Numerical Methods in Fluids 76~(7) (2014) 450--470.

\bibitem{barsukow18hypproceeding}
W.~Barsukow, Stationary states of finite volume discretizations of
  multi-dimensional linear hyperbolic systems, in: XVII International
  Conference on Hyperbolic Problems, Vol.~10, AIMS Series on Applied
  Mathematics, 2020, pp. 296--303.

\bibitem{guillard99}
H.~Guillard, C.~Viozat, On the behaviour of upwind schemes in the low {M}ach
  number limit, Computers \& fluids 28~(1) (1999) 63--86.

\bibitem{barsukow18thesis}
W.~Barsukow, Low {M}ach number finite volume methods for the acoustic and
  {E}uler equations, Doctoral thesis, University of Wuerzburg (2018).

\bibitem{turkel87}
E.~Turkel, Preconditioned methods for solving the incompressible and low speed
  compressible equations, Journal of computational physics 72~(2) (1987)
  277--298.

\bibitem{bispen17}
G.~Bispen, M.~Lukacova-Medvidova, L.~Yelash, Asymptotic preserving {IMEX}
  finite volume schemes for low {M}ach number {E}uler equations with
  gravitation, Journal of Computational Physics 335 (2017) 222--248.

\bibitem{boscarino19}
S.~Boscarino, J.-M. Qiu, G.~Russo, T.~Xiong,
  \href{http://www.sciencedirect.com/science/article/pii/S0021999119303109}{A
  high order semi-implicit imex weno scheme for the all-mach isentropic euler
  system}, Journal of Computational Physics 392 (2019) 594 -- 618.
\newblock \href {http://dx.doi.org/https://doi.org/10.1016/j.jcp.2019.04.057}
  {\path{doi:https://doi.org/10.1016/j.jcp.2019.04.057}}.
\newline\urlprefix\url{http://www.sciencedirect.com/science/article/pii/S0021999119303109}

\bibitem{thomann19}
A.~Thomann, M.~Zenk, G.~Puppo, C.~Klingenberg, An all speed second order imex
  relaxation scheme for the euler equations, arXiv preprint arXiv:1907.08398.

\bibitem{despres17}
B.~Despr{\'e}s, Numerical methods for {E}ulerian and {L}agrangian conservation
  laws, Birkh{\"a}user, 2017.

\bibitem{roe17}
P.~Roe, Multidimensional upwinding, Handbook of Numerical Analysis 18 (2017)
  53--80.

\bibitem{leveque02}
R.~J. LeVeque, Finite volume methods for hyperbolic problems, Vol.~31,
  {C}ambridge {U}niversity {P}ress, 2002.

\bibitem{sidilkover02}
D.~Sidilkover, Factorizable schemes for the equations of fluid flow, Applied
  numerical mathematics 41~(3) (2002) 423--436.

\bibitem{chalons10}
C.~Chalons, F.~Coquel, E.~Godlewski, P.-A. Raviart, N.~Seguin, Godunov-type
  schemes for hyperbolic systems with parameter-dependent source: the case of
  euler system with friction, Mathematical Models and Methods in Applied
  Sciences 20~(11) (2010) 2109--2166.

\bibitem{girardin14}
M.~Girardin, Asymptotic preserving and all-regime lagrange-projection like
  numerical schemes: application to two-phase flows in low mach regime, Ph.D.
  thesis (2014).

\bibitem{jin95}
S.~Jin, Z.~Xin, The relaxation schemes for systems of conservation laws in
  arbitrary space dimensions, Communications on pure and applied mathematics
  48~(3) (1995) 235--276.

\bibitem{bouchut04}
F.~Bouchut, Nonlinear stability of finite Volume Methods for hyperbolic
  conservation laws and Well-Balanced schemes for sources, Springer Science \&
  Business Media, 2004.

\bibitem{holden10}
H.~Holden, K.~H. Karlsen, K.-A. Lie, Splitting methods for partial differential
  equations with rough solutions: {A}nalysis and {MATLAB} programs, Vol.~11,
  European Mathematical Society, 2010.

\bibitem{gresho90}
P.~M. Gresho, S.~T. Chan, On the theory of semi-implicit projection methods for
  viscous incompressible flow and its implementation via a finite element
  method that also introduces a nearly consistent mass matrix. {P}art 2:
  {I}mplementation, International Journal for Numerical Methods in Fluids
  11~(5) (1990) 621--659.

\bibitem{miczek15}
F.~Miczek, F.~R{\"o}pke, P.~Edelmann, New numerical solver for flows at various
  {M}ach numbers, Astronomy \& Astrophysics 576 (2015) A50.

\end{thebibliography}

\appendix

\section{Linear stability of the multi-dimensional acoustic scheme} \label{ap:stability}

Recall (see \cite{barsukow18thesis,barsukow17a}) that a stationarity preserving scheme for the linear $m \times m$ system\footnote{Indices never denote derivatives here.}
\begin{align}
 \del_t q + J_x \del_x q + J_y \del_y q &= 0 & q &: \mathbb R^+_0 \times \mathbb R^d \to \mathbb R^m\\
 && J_x, J_y &: \mathcal M^{m \times m}
\end{align}
reads
\begin{align}
 \del_t q &+ J_x \frac{ \{\{ [q]_{i\pm1} \}\}_{j\pm\frac12} }{8 \Delta x} + J_y \frac{  [\{\{q\}\}_{i\pm\frac12} ]_{j\pm1} }{8 \Delta y} \nonumber\\
 &- \frac12 \sign J_x \left( J_x \frac{ \{\{ [[q]]_{i\pm\frac12} \}\}_{j\pm\frac12} }{4 \Delta x} + J_y \frac{  [[q]_{i\pm1} ]_{j\pm1} }{4 \Delta y} \right ) \label{eq:appgenscheme} \\
 &- \frac12 \sign J_y \left( J_x \frac{ [ [q]_{i\pm1} ]_{j\pm1} }{4 \Delta x} + J_y \frac{  [[\{\{q\}\}_{i\pm\frac12} ]]_{j\pm\frac12} }{4 \Delta y}  \right ) = 0 \nonumber
\end{align}

Here, both $|J_x|$ and $\sign J_x$ are defined on the eigenvalues.

This scheme reduces to the scheme \eqref{eq:multidschemeacousticu}--\eqref{eq:multidschemeacousticp} for the special case of the acoustic equations \eqref{eq:acoustic2dv}--\eqref{eq:acoustic2dp}.

The discrete Fourier transform replaces $q_{ij}$ by
\begin{align}
 \hat q(t) \exp(\ii k_x \Delta x \cdot i + \ii k_y \Delta y \cdot j)
\end{align}
Define the translation factors
\begin{align}
 t_x &:= \exp(\ii k_x \Delta x) & t_y &:= \exp(\ii k_y \Delta y)
\end{align}

Then the Fourier transform of the scheme \eqref{eq:appgenscheme} can be written as
\begin{align}
 \del_t \hat q + f \mathcal D \hat q &= 0
\end{align}
with
\begin{align}
 \mathcal D &= J_x \frac{(t_x-1)(t_x+1)}{2 t_x \Delta x} \frac{(1 + t_y)^2}{4 t_y}   + J_y \frac{(t_y-1)(t_y+1)}{2 t_y \Delta y} \frac{(1 + t_x)^2}{4 t_x} \\
 f &= \id - \sign J_x \frac{t_x-1}{t_x+1} - \sign J_y \frac{t_y-1}{t_y+1}
\end{align}

Thus the amplification matrix for explicit time integration is
\begin{align}
 \mathcal A = \id - \Delta t f \mathcal D
\end{align}

Von Neumann stability requires the eigenvalues of $\mathcal A$ not to be larger than 1 in absolute value. This then puts a condition on $\Delta t$.

\begin{theorem}
 For $\Delta y = \Delta x$, the scheme \eqref{eq:appgenscheme} for the acoustic equations \eqref{eq:acoustic2dv}--\eqref{eq:acoustic2dp} is von Neumann stable if
 \begin{align}
 c \Delta t \leq \Delta x
 \end{align}
 i.e. for $\textsc{CFL}\leq 1$.
\end{theorem}
\begin{proof}
The first part of the stability proof involves showing that actually, for $\Delta y = \Delta x$ $f$ and $\mathcal D$ commute and thus are simultaneously diagonalizable:

\begin{align}
 f \mathcal D - \mathcal D f &= 
 (- \sign J_x J_y - \sign J_y  J_x + J_x \sign J_y + J_y \sign J_x) \nonumber\\&\phantom{mmmmmmmmmm}\cdot \frac{(t_y-1)(t_y+1)}{2 t_y \Delta x} \frac{(1 + t_x)(t_x-1)}{4 t_x}
\end{align}

Taking
\begin{align}
 J_x &= \left ( \begin{array}{ccc} 0&0&c\\0&0&0\\c&0&0 \end{array} \right ) & J_y &= \left ( \begin{array}{ccc} 0&0&0\\0&0&c\\0&c&0 \end{array} \right )
\end{align}
one shows upon direct computation that
\begin{align}
 J_x \sign J_y  - \sign J_x J_y &= 0\\
 J_y \sign J_x - \sign J_y  J_x &= 0
\end{align}
This might point to some algebraic structure of the equations that is yet to be understood.

\newcommand{\diag}{\mathrm{diag}}

Thus the eigenvalue of $f\mathcal D$ is the product of the eigenvalues of $f$ and $\mathcal D$. These are
\begin{align}
 f &\cong \diag\left(1, 1 + \frac{\sqrt{2} \sqrt{s}}{(1+t_x)(1+t_y)}, 1 - \frac{\sqrt{2} \sqrt{s}}{(1+t_x)(1+t_y)} \right)\\
 \mathcal D &\cong \diag\left(0, - c \frac{(1+t_x)(1+t_y) \sqrt{s}}{4 \sqrt{2} \Delta x t_xt_y}, c \frac{(1+t_x)(1+t_y) \sqrt{s} }{4 \sqrt{2} \Delta x t_xt_y} \right )
\end{align}
with $s = 1 + t_x^2 + t_y^2 + t_x^2  t_y^2 - 4 t_x t_y$.

Call $\lambda$ an eigenvalue of the evolution matrix $f \mathscr D$. Then
\begin{align}
 \alpha = 1 - \Delta t \lambda
\end{align}
is an eigenvalue of the amplification matrix. Its absolute value squared is
\begin{align}
 |\alpha|^2 = (1 - \Delta t \lambda)(1 - \Delta t \lambda^*) = 1 - \Delta t (\lambda + \lambda^*) + \Delta t^2 \lambda \lambda^* 
\end{align}
For stability $|\alpha|^2 \leq 1$, which amounts to
\begin{align}
 \Delta t \leq \frac{\lambda + \lambda^*}{\lambda \lambda^*}
\end{align}

This computation can be easily performed for the above eigenvalues and after inserting the definitions of $t_x$ and $t_y$ yields
\begin{align}
 \Delta t \leq \frac{\Delta x}{c} \cdot \frac{4}{3 + \cos \beta_x + \cos \beta_y - \cos \beta_x \cos \beta_y} =: \frac{\Delta x}{c} f(\beta_x, \beta_y)
\end{align}
with $\beta_x = k_x \Delta x$, $\beta_y = k_y \Delta y$.

\begin{figure}
 \centering
 \includegraphics[width=0.6\textwidth]{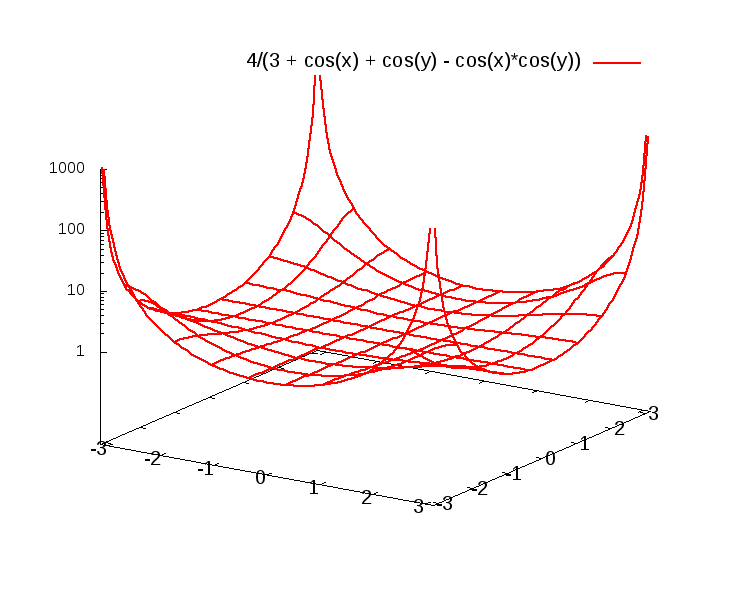}
 \caption{The function $f(\beta_x, \beta_y)$ in the parameter range $[-\pi,\pi]^2$.}
 \label{fig:stabilityfunction}
\end{figure}

A plot of the function $f$ is shown in Figure \ref{fig:stabilityfunction}. The derivative is easily found:

\begin{align}
 \nabla_\beta f(\beta_x, \beta_y) = \frac{4}{(3 + \cos \beta_x + \cos \beta_y - \cos \beta_x \cos \beta_y)^2} \vecc{%
 \sin \beta_x (1-\cos \beta_y)}{%
 \sin \beta_y (1-\cos \beta_x)}
\end{align}
The derivative vanishes for $\sin \beta_x = 0$ or $\sin \beta_y = 0$. The value of $f$ in all such points is 1, which proves
\begin{align}
 c \Delta t \leq \Delta x
\end{align}
as the bound for stability if $\Delta y = \Delta x$.
\end{proof}

Note that, e.g. in two dimensions, a dimensionally split scheme generally is only stable under the more restrictive stability condition $\textsc{CFL} < \frac12$.

\section{On the interpretation of results obtained by discrete asymptotic analysis for time-explicit numerical schemes} \label{ap:asymptotic}

Consider the following model equation
\begin{align}
 \frac{q^{n+1} - q^n}{\Delta t} + \frac{1}{\epsilon} (q^n - a) = 0 \label{eq:modelequationasymptotic}
\end{align}
with $a \in \mathbb R$. This is a forward Euler discretization of the ODE $\frac{\dd q}{\dd t} + \frac{1}{\epsilon} (q - a) = 0$. As $\epsilon \to 0$, formal asymptotic analysis of \eqref{eq:modelequationasymptotic} gives
\begin{align}
 q(t) &= q^{(0)} + \epsilon q^{(1)} + \epsilon^2 q^{(2)} + \mathcal O(\epsilon^3)\\
 (q^n)^{(0)}  &= a \label{eq:modelequationlimit}
\end{align}

It is a common source of irritation that the asymptotic analysis seemingly yields a statement about the timestep $n$ rather than $n+1$, and that \eqref{eq:modelequationlimit} might be in contradiction with the data at timestep $n$. However, the asymptotic analysis is not a proof, and the statement thus obtained has to be taken \emph{cum grano salis}. The correct interpretation is subject of this section. The model equation \eqref{eq:modelequationasymptotic} admits the exact solution
\begin{align}
 q^n = \left( 1 - \frac{\Delta t}{\epsilon}\right )^n (q^0 - a) + a
\end{align}
where $q^0$ denotes the initial value. For reasons of stability, $\Delta t$ cannot be independent of $\epsilon$, as is usually the case for a time-explicit integrator. Choose $\Delta t = \epsilon \tau$. Then, stability requires $1 - \tau > -1$, i.e. $0 < \tau < 2$. Integrating up to a final time of $T$ means $n = T/\Delta t = \frac{T}{\epsilon \tau}$:
\begin{align}
 q^n = \left( 1 - \tau \right )^{\frac{T}{\epsilon \tau}} (q^0 - a) + a
\end{align}

Consider now the case $q^0 \neq a$. The result \eqref{eq:modelequationlimit} of the asymptotic analysis requires an interpretation, because it clearly is wrong when applied to $n=0$ and taken literally. Upon the evolution in time, the difference $q^0-a$ decays exponentially. The time $T_{1/2}$ for it to reduce by one half is
\begin{align}
  T_{1/2} =  \epsilon \tau \frac{\log \frac12}{\log \left( 1 - \tau \right )} \in \mathcal O(\epsilon)
\end{align}
Thus, as $\epsilon \to 0$, the limit solution $q^0 = a$ is reached exponentially quickly, with a decay time $\mathcal O(\epsilon)$. The statement \eqref{eq:modelequationlimit} of the asymptotic analysis should thus be understood for $n>0$ only.

A possible remedy seems to choose an implicit time integration
\begin{align}
 \frac{q^{n+1} - q^n}{\Delta t} + \frac{1}{\epsilon} (q^{n+1} - a) = 0 \label{eq:modelequationasymptoticimplicit}
\end{align}
because then, as $\epsilon \to 0$, formal asymptotic analysis gives the more pleasing
\begin{align}
 (q^{n+1})^{(0)}  &= a 
\end{align}
Indeed, this argument is sometimes used to support the need for implicit time integration in the context of low Mach number compliant schemes for Euler. However, the exact solution to \eqref{eq:modelequationasymptoticimplicit} is
\begin{align}
 q^n = \left( 1 + \frac{\Delta t}{\epsilon}\right )^{-n} (q^0 - a) + a
\end{align}
which shows, apart from the improved stability, the same behaviour, and $q^{n}$ approaches $a$ exponentially in time with a decay rate proportional to $\epsilon$. Thus, the advantage of implicit time integration is only the absence of a CFL condition, not any fundamentally different asymptotic behaviour.

\end{document}